\documentclass[11pt]{article}

\usepackage{amsmath}
\usepackage{amssymb}
\usepackage{amsthm}
\usepackage{latexsym}
\usepackage{color}
\usepackage{graphicx}
\usepackage{color}
\usepackage[T1]{fontenc}
\usepackage[english]{babel}
\usepackage{geometry}
\usepackage{url}
\DeclareSymbolFont{calletters}{OMS}{cmsy}{m}{n}
\DeclareSymbolFontAlphabet{\mathcal}{calletters}

\def\be{\begin{eqnarray}}
\def\ee{\end{eqnarray}}

\def\b*{\begin{eqnarray*}}
\def\e*{\end{eqnarray*}}

\newtheorem{Theorem}{Theorem}
\newtheorem{Definition}{Definition}
\newtheorem{Proposition}{Proposition}
\newtheorem{Assumption}{Assumption}
\newtheorem{Lemma}{Lemma}
\newtheorem{Remark}{Remark}

\makeatletter \@addtoreset{equation}{section}

\@addtoreset{Definition}{section}

\@addtoreset{Theorem}{section}

\@addtoreset{Proposition}{section}

\@addtoreset{Property}{section}

\@addtoreset{Assumption}{section}

\@addtoreset{Corollary}{section}

\@addtoreset{Lemma}{section}

\@addtoreset{Remark}{section}

\@addtoreset{Example}{section}

\newcommand{\rmi}{{\rm (i)$\>\>$}}
\newcommand{\rmii}{{\rm (ii)$\>\>$}}
\newcommand{\rmiii}{{\rm (iii)$\>\>$}}
\newcommand{\rmiv}{{\rm (iv)$\>\>$}}
\newcommand{\rmv}{{\rm (v)$\>\>$}}

\newcommand{\No}[1]{\left\|#1\right\|}     
\newcommand{\abs}[1]{\left|#1\right|}     

\def \C{\mathbb{C}}
\def \D{\mathbb{D}}
\def \E{\mathbb{E}}
\def \F{\mathbb{F}}
\def \H{\mathbb{H}}
\def \L{\mathbb{L}}
\def \M{\mathbb{M}}
\def \N{\mathbb{N}}
\def \P{\mathbb{P}}
\def \Q{\mathbb{Q}}
\def \R{\mathbb{R}}
\def \S{\mathbb{S}}
\def \X{\mathbb{X}}

\def\Bc{{\cal B}}

\def\Ec{{\cal E}}
\def\Fc{{\cal F}}

\def\Hc{{\cal H}}

\def\Jc{{\cal J}}
\def\Kc{{\cal K}}

\def\Mc{{\cal M}}

\def\Pc{{\cal P}}

\def\Tc{{\cal T}}
\def\Uc{{\cal U}}

\def\Xc{{\cal X}}
\def\Yc{{\cal Y}}
\def\Zc{{\cal Z}}

\def\Fb{{\bar F}}

\def\Pb{{\overline \P}}

\def\Ych{{\widehat \Yc}}
\def\Zch{{\widehat \Zc}}

\def\einf{{\rm ess \, inf}}
\def\esup{{\rm ess \, sup}}

\def\x{\times}

\def\={\;=\;}
\def\.{\;.}

\def\eps{\varepsilon}

\def\reff#1{{\rm(\ref{#1})}}

\def \i{1\!\mbox{\rm I}}
\def\1{{\bf 1}}
\def \ep{\hbox{ }\hfill{ ${\cal t}$~\hspace{-5.1mm}~${\cal u}$   } }

\def \proof{{\noindent \bf Proof. }}
\def \ep{\hbox{ }\hfill$\Box$}

\def\Om{\Omega}
\def\om{\omega}

\def\crochetX{\langle X \rangle}
\def\ah{\widehat{a}}
\def\Ycb{\overline{\Yc}}
\def\Bf{\mathfrak{B}}
\def\Sp{\mathbb S_d^{\geq 0}}
\def\ox{\otimes}
\def \w{\mathsf{w}}
\def\Hlr{\langle H \rangle}

\def\Omb{\overline{\Om}}
\def\omb{\bar \om}
\def\Fb{\overline{\F}}
\def\Fcb{\overline{\Fc}}
\def\Ycb{\overline{\Yc}}
\def\Zcb{\overline{\Zc}}
\def\Mcb{\overline{\Mc}}
\def\Psib{\overline{\Psi}}
\def\0{\mathbf{0}}
\def \Ec{\mathcal{E}}

\def\normeL2#1{\left\|{#1}\right\|_{L^2}}

 \title{Stochastic control for a class of nonlinear kernels and applications
 }

\author{
 Dylan {\sc Possama\"{i}} \footnote{Universit\'e Paris--Dauphine, PSL Research University, CNRS, CEREMADE, 75016 Paris, France, possa-mai@ceremade.dauphine.fr}
 \and Xiaolu {\sc Tan}\footnote{Universit\'e Paris--Dauphine, PSL Research University, CNRS, CEREMADE, 75016 Paris, France, tan@ce-remade.dauphine.fr}
      \and Chao {\sc Zhou}\footnote{Department of Mathematics, National University of Singapore, Singapore, matzc@nus.edu.sg}}
             \date{\today}

\begin{document}

 \maketitle

 \begin{abstract}
 We consider a stochastic control problem for a class of nonlinear kernels. More precisely, our problem of interest consists in the optimisation, over a set of possibly non--dominated probability measures, of solutions of backward stochastic differential equations (BSDEs). Since BSDEs are nonlinear generalisations of the traditional (linear) expectations, this problem can be understood as stochastic control of a family of nonlinear expectations, or equivalently of nonlinear  kernels. Our first main contribution is to prove a dynamic programming principle for this control problem in an abstract setting, which we then use to provide a semi--martingale characterisation of the value function. We next explore several applications of our results. We first obtain a wellposedness result for second order BSDEs (as introduced in \cite{soner2012wellposedness}) which does not require any regularity assumption on the terminal condition and the generator. Then we prove a nonlinear optional decomposition in a robust setting, extending recent results of \cite{nutz2015robust}, which we then use to obtain a super--hedging duality in uncertain, incomplete and nonlinear financial markets. Finally, we relate, under additional regularity assumptions, the value function to a viscosity solution of an appropriate path--dependent partial differential equation (PPDE).

\vspace{5mm}

\vspace{5mm}


\end{abstract}

\section{Introduction}
The dynamic programming principle (DPP for short) has been a major tool in the control theory, since the latter took off in the 1970's. Informally speaking, this principle simply states that a global optimisation problem can be split into a series of local optimisation problems. Although such a principle is extremely intuitive, its rigorous justification has proved to be a surprisingly difficult issue. Hence, for stochastic control problems, the dynamic programming principle is generally based on the stability of the controls with respect to conditioning and concatenation, together with a measurable selection argument, which, roughly speaking, allow to prove the measurability of the associated value function, as well as constructing almost optimal controls through "pasting". This is exactly the approach followed by Bertsekas and Shreve \cite{bertsekas1978stochastic}, and Dellacherie \cite{dellacherie1985quelques} for discrete time stochastic control problems. In continuous time, a comprehensive study of the dynamic programming principle remained more elusive. Thus, El Karoui, in \cite{el1981aspects}, established the dynamic programming principle for the optimal stopping problem in a continuous time setting, using crucially the strong stability properties of stopping times, as well as the fact that the measurable selection argument can be avoided in this context, since an essential supremum over stopping times can be approximated by a supremum over a countable family of random variables. Later, for general controlled Markov processes (in continuous time) problems, El Karoui, Huu Nguyen and Jeanblanc \cite{el1987compactification} provided a framework to derive the dynamic programming principle using the measurable selection theorem, by interpreting the controls as probability measures on the canonical trajectory space (see e.g. Theorems 6.2, 6.3 and 6.4 of \cite{el1987compactification}). Another commonly used approach to derive the DPP was to bypass the measurable selection argument by proving, under additional assumptions, {\it a priori} regularity of the value function. This was the strategy adopted, among others, by Fleming and Soner \cite{fleming2006controlled}, and in the so--called weak DPP of Bouchard and Touzi \cite{bouchard2011weak}, which has then been extended by Bouchard and Nutz \cite{bouchard2012weak,bouchard2015stochastic} and Bouchard, Moreau and Nutz \cite{bouchard2014stochastic} to optimal control problems with state constraints as well as to differential games (see also Dumitrescu, Quenez and Sulem \cite{dumitrescu2015weak} for a combined stopping--control problem on BSDEs).
	One of the main motivations of this weak DPP is that it is generally enough to characterise the value function as a viscosity solution of the associated Hamilton--Jacobi--Bellman partial differential equation (PDE).
	Let us also mention the so--called stochastic Perron's method, which has been developed by Bayraktar and S\^irbu, see {\it e.g.} \cite{bayraktar2012stochastic,bayraktar2013stochastic}, which allows, for Markov problems, to obtain the viscosity solution characterization of the value function without using the DPP, and then to prove the latter {\it a posteriori}. Recently, motivated by the emerging theory of robust finance, Nutz et al. \cite{neufeld2013superreplication,nutz2013constructing} gave a framework which allowed to prove the dynamic programming principle for sub--linear expectations (or equivalently a non--Markovian stochastic control problem),
	where the essential arguments are close to those in \cite{el1987compactification}, though the presentation is more modern, pedagogic and accessible.
	The problem in continuous--time has also been studied by El Karoui and Tan \cite{karoui2013capacities,karoui2013capacities2}, in a more general context than the previous references, but still based on the same arguments as in \cite{el1987compactification} and \cite{neufeld2013superreplication}.

\vspace{0.5em}
\noindent However, all the above works consider only what needs to be qualified as the sub--linear case. Indeed, the control problems considered consists generically in the maximisation of a family of expectations over the set of controls. Nonetheless, so--called nonlinear expectations on a given probability space (that is to say operators acting on random variables which preserve all the properties of expectations but linearity) have now a long history, be it from the capacity theory, used in economics to axiomatise preferences of economic agents which do not satisfy the usual axiom's of von Neumann and Morgenstern, or from the seminal $g-$expectations (or BSDEs) introduced by Peng \cite{peng1997backward}. Before pursuing, let us just recall that in the simple setting of a probability space carrying a Brownian motion $W$, with its (completed) natural filtration $\F$, finding the solution of a {\rm BSDE} with generator $g$ and terminal condition $\xi\in\Fc_T$ amounts to finding a pair of $\F-$progressively measurable processes $(Y,Z)$ such that
$$Y_t=\xi-\int_t^Tg_s(Y_s,Z_s)ds-\int_t^TZ_s\cdot dW_s,\ t\in[0,T],\ a.s.$$
This theory is particularly attractive from the point of view of stochastic control, since it is constructed to be filtration (or time) consistent, that is to say that its conditional version satisfies a tower property similar to that of linear expectations, which is itself a kind of dynamic programming principle. Furthermore, it has been proved by Coquet {\it et al.} \cite{coquet2002filtration} that essentially all filtration consistent nonlinear expectations satisfying appropriate domination properties could be represented with BSDEs (we refer the reader to \cite{hu2008representation} and \cite{cohen2012representing} for more recent extensions of this result). Our first contribution in this paper, in Section \ref{sec:1}, is therefore to generalise the measurable selection argument to derive the dynamic programming principle in the context of optimal stochastic control of nonlinear expectations (or kernels) which can be represented by BSDEs (which as mentioned above is not such a stringent assumption). We emphasise that such an extension is certainly not straightforward. Indeed, in the context of linear expectations, there is a very well established theory studying how the measurability properties of a given map are impacted by its integration with respect to a so--called stochastic kernel (roughly speaking one can see this as a regular version of a conditional expectation in our context, see for instance \cite[Chapter 7]{bertsekas1978stochastic}). For instance, integrating a Borel map with respect to a Borel stochastic kernel preserves the Borel measurability. However, in the context of BSDEs, one has to integrate with respect to nonlinear stochastic kernels, for which, as far as we know, no such theory of measurability exists. Moreover, we also obtain a semi--martingale decomposition for the value function of our control problem. This is the object of Section \ref{sec:2}.

\vspace{0.5em}
\noindent Let us now explain where our motivation for studying this problem stems from. The problem of studying a controlled system of BSDEs is not new. For instance, it was shown by El Karoui and Quenez \cite{el1991programmation,el1995dynamic} and Hamad\`ene and Lepeltier \cite{hamadene1995backward} (see also \cite{el1997backward} and the references therein) that a stochastic control problem with control on the drift only could be represented via a controlled family of BSDEs (which can actually be itself represented by a unique {\rm BSDE} with convex generator). More recently, motivated by obtaining probabilistic representations for fully nonlinear PDEs, Soner, Touzi and Zhang \cite{soner2012wellposedness,soner2013dual} (see also the earlier works \cite{cheridito2007second} and \cite{soner2011martingale}) introduced a notion of second--order BSDEs (2BSDEs for short), whose solutions could actually be written as a supremum, over a family of non--dominated probability measures (unlike in \cite{el1995dynamic} where the family is dominated), of standard BSDEs. Therefore the 2BSDEs fall precisely in the class of problem that we want to study, that is stochastic control of nonlinear kernels. The authors of \cite{soner2012wellposedness,soner2013dual} managed to obtain the dynamic programming principle, but under very strong continuity assumptions w.r.t. $\omega$ on the terminal condition and the generator of the BSDEs, and obtained a semi--martingale decomposition of the value function of the corresponding stochastic control problem, which ensured wellposedness of the associated {\rm 2BSDE}. Again, these regularity assumptions are made to obtain the continuity of the value function {\it a priori}, which allows to avoid completely the use of the measurable selection theorem. Since then, the {\rm 2BSDE} theory has been extended by allowing more general generators, filtrations and constraints (see \cite{kazi2015second,kazi2015second2,matoussi2014second,matoussi2013second,possamai2013second1,possamai2013second}), but no progress has been made concerning the regularity assumptions. However, the 2BSDEs (see for instance \cite{matoussi2015robust}) have proved to provide a particularly nice framework to study the so--called robust problems in finance, which were introduced by \cite{avellaneda1995pricing,lyons1995uncertain} and in a more rigorous setting by \cite{denis2006theoretical}. However, the regularity assumptions put strong limitations to the range of the potential applications of the theory. 

\vspace{0.5em}
\noindent We also would like to mention a related theory introduced by Peng \cite{peng2012nonlinear}, and developed around the notion of $G-$expectations, which lead to the $G-$BSDEs (see \cite{hu2014backward,hu2014comparison}). Instead of working on a fixed probability space carrying different probability measures corresponding to the controls, they work directly on a sub--linear expectation space in which the canonical process already incorporates the different measures, without having to refer to a probabilistic setting. Although their method of proof is different, since they mainly use PDE arguments to construct a solution in the Markovian case and then a closure argument, the final objects are extremely close to 2BSDEs, with similar restrictions in terms of regularity. Moreover, the PDE approach they use is unlikely to be compatible with a theory without any regularity, since the PDEs they consider need at the very least to have a continuous solution. On the other hand, there is more hope for the probabilistic approach of the 2BSDEs, since, as shown in \cite{nutz2013constructing} in the case of linear expectations (that is when the generator of the BSDEs is $0$), everything can be well defined by assuming only that the terminal condition is (Borel) measurable. 

\vspace{0.5em}
\noindent There is a third theory which shares deep links with 2BSDEs, namely that of viscosity solutions of fully nonlinear path dependent PDEs (PPDEs for short), which has been introduced recently by Ekren, Keller, Touzi and Zhang \cite{ekren2014viscosity,ekren2016viscosity,ekren2012viscosity}. Indeed, they showed that the solution of a {\rm 2BSDE}, with a generator and a terminal condition uniformly continuous (in $\omega$), was nothing else than the viscosity solution of a particular PPDE, making the previous theory of 2BSDEs a special case of the theory of PPDEs. The second contribution of our paper is therefore that we show (a suitable version of) the value function for which we have obtained the dynamic programming principle provides a solution to a {\rm 2BSDE} without requiring any regularity assumption, a case which cannot be covered by the PPDE theory. This takes care of the existence problem, while we tackle, as usual, the uniqueness problem through {\it a priori} $\L^p$ estimates on the solution, for any $p>1$. We emphasise that in the very general setting that we consider, the classical method of proof fails (in particular since the filtration we work with is not quasi--left continuous in general), and the estimates follow from a general result that we prove in our accompanying paper \cite{bouchard2015unified}. In particular, our wellposedness results contains directly as a special case the theory of BSDEs, which was not immediate for the 2BSDEs of \cite{soner2012wellposedness}, or the $G-${\rm BSDE} (indeed one first has to prove a Lusin type result showing that the closure of the space of uniformly continuous random variables for the $\L^p(\P)$ norm for a fixed $\P$ is actually the whole space). Moreover, the class of probability measures that we can consider is much more general than the ones considered in the previous literature, even allowing for degeneracy of the diffusion coefficient. This is the object of Section \ref{sec:3}. 

\vspace{0.5em}
\noindent The rest of the paper is mainly concerned with applications of the previous theory. First, in Section \ref{sec:4}, we use our previous results to obtain a nonlinear and robust generalisation of the so--called optional decomposition for super--martingales (see for instance \cite{el1995dynamic,kramkov1996optional} and the other references given in Section \ref{sec:4} for more details), which is new in the literature. This allows us to introduce, under an additional assumption stating that the family of measures is roughly speaking rich enough, a new notion of solutions, which we coined saturated 2BSDEs. This new formulation has the advantage that it allows us to get rid of the orthogonal martingales which generically appear in the definition of a {\rm 2BSDE} (see Definitions \ref{def:1} and \ref{def:2} for more details). This is particularly important in some applications, see for instance \cite{cvitanic2015dynamic}. We then give a duality result for the robust pricing of contingent claims in nonlinear and incomplete financial markets. Finally, in Section \ref{sec:5}, we recall in our context the link between 2BSDEs and PPDEs when we work under additional regularity assumptions. Compared to \cite{ekren2016viscosity}, our result can accommodate degenerate diffusions.

\vspace{0.5em}
\noindent To conclude this introduction, we really want to insist on the fact that our new results have much more far--reaching applications, and are not a mere mathematical extension. Indeed, in the paper \cite{cvitanic2015dynamic}, the wellposedness theory of 2BSDEs we have obtained is used crucially to solve general Principal--Agent problems in contracting theory, when the agent controls both the drift and the volatility of the corresponding output process (we refer the reader to the excellent monograph \cite{cvitanic2012contract} for more details on contract theory), a problem which could not be treated with the techniques prevailing in the previous literature. Such a result has potential applications in many fields, ranging from economics (see for instance \cite{cvitanic2014moral,mastrolia2015moral}) to energy management (see \cite{aid2017principal}).

\subsection*{Notations}
	Throughout this paper, we fix the constants $p > \kappa > 1$.
	Let $\mathbb{N}^*:=\mathbb{N}\setminus\{0\}$ and let $\mathbb{R}_+^*$ be the set of real positive numbers. 
	For every $d-$dimensional vector $b$ with $d\in \mathbb{N}^*$, we denote by $b^{1},\ldots,b^{d}$ its coordinates and for $\alpha,\beta \in \R^d$ we denote by $\alpha\cdot \beta$ the usual inner product, with associated norm $\No{\cdot}$, which we simplify to $|\cdot|$ when $d$ is equal to $1$. We also let ${\bf 1}_d$ be the vector whose coordinates are all equal to $1$. 
	For any $(\ell,c)\in\mathbb N^*\times\mathbb N^*$, $\mathcal M_{\ell,c}(\mathbb R)$ will denote the space of $\ell\times c$ matrices with real entries. Elements of the matrix $M\in\mathcal M_{\ell,c}$ will be denoted by $(M^{i,j})_{1\leq i\leq \ell,\; 1\leq j\leq c}$, and the transpose of $M$ will be denoted by $M^\top$. When $\ell=c$, we let $\mathcal M_{\ell}(\mathbb R):=\mathcal M_{\ell,\ell}(\mathbb R)$. We also identify $\mathcal M_{\ell,1}(\R)$ and $\R^\ell$. The identity matrix in $\mathcal M_\ell(\R)$ will be denoted by $I_\ell$.
	Let $\Sp$ denote the set of all symmetric positive semi-definite $d\x d$ matrices.
	We fix a map $\psi : \Sp \longrightarrow \Mc_d(\R)$ which is (Borel) measurable and satisfies $\psi(a) \psi(a)^\top = a$ for all $a \in \Sp$, 
	and denote $a^{\frac 12} := \psi(a)$.
	{Finally, we denote by $a^{\oplus}$ the Moore--Penrose pseudoinverse of $a \in \Sp$.
	In particular, the map $a \longmapsto a^{\oplus} = \lim_{\delta \searrow 0} (a^{\top} a + \delta I_d)^{-1} a^{\top}$ is Borel measurable.
	}

\subsection{Probabilistic framework}

	Let us first present our probabilistic framework with a canonical space and the associated notations.

\subsubsection{Canonical space}

\vspace{0.5em}
	Let $d \in \N^{\ast}$,
	we denote by $\Om := C\left( \left[ 0,T \right], \mathbb{R}^d \right)$ the canonical space of all $\R^d-$valued continuous paths $\om$ on $[0,T]$ such that $\om_0 = 0$,
	equipped with the canonical process $X$, {\it i.e.} $X_t(\om) := \om_t,$ for all $\om  \in \Om$.
	Denote 
	by $\F = (\Fc_t)_{0 \le t \le T}$ the canonical filtration generated by $X$,
	and by $ \mathbb{F}_+ = ( \mathcal{F}^+_t)_{ 0 \le t \le T } $
	the right limit of  $ \mathbb{F} $ with $\Fc_t^+ := \Fc_{t+} \mathrel{\mathop:}= \cap_{ s> t} \Fc_s $ for all $t \in [0, T)$ and $\Fc^+_T := \Fc_T$.
	We equip $\Om$ with the uniform convergence norm $\No{\om}_{\infty} := \sup_{0 \le t \le T} \No{\om_t}$, so that the Borel $\sigma-$field of $\Om$ coincides with $\Fc_T$.
	Let $\P_0$ denote the Wiener measure on $\Om$ under which $X$ is a Brownian motion.

	\vspace{0.5em}
	
	\noindent Let $\M_1$ denote the collection of all probability measures on $(\Om, \Fc_T)$.
	Notice that $\M_1$ is a Polish space equipped with the weak convergence topology.
	We denote by $\Bf$ its Borel $\sigma-$field.
	Then for any $\P\in \M_1$, denote by $\Fc_t^{\P}$ the completed $\sigma-$field of $\Fc_t$ under $\P$.
	Denote also the completed filtration by $\F^{\P} = \big(\Fc_t^{\P} \big)_{t \in [0,T]}$
	and $\F^{\P}_+$ the right limit of $\F^{\P}$, so that $\F^{\P}_+$ satisfies the usual conditions.
	Moreover,  for $\Pc \subset \M_1$, we introduce the universally completed filtration 
	$\F^{U} := \big(\Fc^{U}_t \big)_{0 \le t \le T}$, $\F^{\Pc} := \big(\Fc^{\Pc}_t \big)_{0 \le t \le T}$,
	and $\F^{\Pc+} := \big(\Fc^{\Pc+}_{t} \big)_{0 \le t \le T}$, defined as follows
	$$
		 \mathcal F^{U}_t
		:=
		\bigcap_{\P \in \M_1} ~\Fc_t^{\P},\ 
		\Fc^{\Pc}_t := \bigcap_{\P \in \Pc} ~\Fc_t^{\P},\ t\in[0,T],\ 
		\Fc^{\Pc+}_t := \Fc^{\Pc}_{t+},\ t\in[0,T), \ \text{and}\ \Fc^{\Pc+}_T:=\Fc^{\Pc}_T.
	$$
	We also introduce an enlarged canonical space $\Omb := \Om \x \Om'$, where $\Om'$ is identical to $\Om$.
	By abuse of notation, we denote by $(X, B)$ its canonical process, {\it i.e.} $X_t(\omb) := \om_t$, $B_t(\omb) := \omega'_t$ for all $\omb := (\om, \omega') \in \Omb$,
	by $\Fb = (\Fcb_t)_{0 \le t \le T}$ the canonical filtration generated by $(X,B)$,
	and by $\Fb^X = (\Fcb^X_t)_{0 \le t \le T}$ the filtration generated by $X$.
	Similarly, we denote the corresponding right--continuous filtrations by $\Fb^X_+$ and $\Fb_+$,
	and the augmented filtration by $\Fb^{X, \Pb}_+$ and $\Fb^{\Pb}_+$,
	given a probability measure $\Pb$ on $\Omb$.

\subsubsection{Semi--martingale measures}

	We say that a probability measure $\P$ on $(\Om, \Fc_T)$ is a semi--martingale measure if $X$ is a semi-martingale under $\P$.
	Then on the canonical space $\Om$, there is some $\F-$progressively measurable non-decreasing process
	(see {\it e.g.} Karandikar \cite{karandikar1995pathwise}, Bichteler \cite{bichteler1981stochastic}, or Neufeld and Nutz \cite[Proposition 6.6]{neufeld2014measurability}),
	denoted by $\crochetX = (\crochetX_t)_{0 \le t \le T}$, which coincides with the quadratic variation of $X$ under each semi--martingale measure $\P$.
	Denote further
	$$
		\ah_t
		:=
		\limsup_{\eps \searrow 0}~ \frac{\crochetX_t - \crochetX_{t-\eps}}{\eps}.
	$$
	For every $t \in [0, T]$, let $\Pc^W_t$ denote the collection of all probability measures $\P$ on $(\Om, \Fc_T)$ such that
	\begin{itemize}
		\item $(X_s)_{s \in [t,T]}$ is a $(\P, \F)-$semi--martingale admitting the canonical decomposition 
		(see e.g. \cite[Theorem I.4.18]{jacod2013limit})
		$$X_s=\int_t^sb^\P_rdr + X^{c,\P}_s, \ s\in [t, T], \ \P-a.s., $$ 
		where $b^\P$ is a $\F^{\P}-$predictable $\R^d-$valued process,
		and $X^{c,\P}$ is the continuous local martingale part of $X$ under $\P$.

		\item $\big(\crochetX_s \big)_{s \in [t,T]}$ is absolutely continuous in $s$ with respect to the Lebesgue measure, and $\ah$ takes values in $\Sp$, $\P-a.s.$
	\end{itemize}

\noindent Given a random variable or process $\lambda$ defined on $\Om$, we can naturally define its extension on $\Omb$ (which, abusing notations slightly, we still denote by $\lambda$) by
	\begin{equation} \label{eq:extension_def}
		\lambda(\omb) := \lambda(\om), ~\forall \omb := (\om, \om') \in \Omb.
	\end{equation}
	In particular, the process $\ah$ can be extended on $\Omb$.
	Given a probability measure $\P \in \Pc^W_t$, we define a probability measure $\Pb$ on the enlarged canonical space $\Omb$ by $\Pb := \P \otimes \P_0$, so that
	$X$ in $(\Omb, \Fcb_T, \Pb, \Fb)$ is a semi--martingale with  the same triplet of characteristics as $X$ in $(\Om, \Fc_T, \P, \F)$,
	$B$ is a $\Fb-$Brownian motion, and $X$ is independent of $B$.
	Then for every $\P \in \Pc^W_t$, there is some $\R^d-$valued, $\Fb-$Brownian motion $W^{\P} = (W^{\P}_r)_{t \le r \le s}$ such that  (see e.g. Theorem 4.5.2 of \cite{stroock2007multidimensional})
	\begin{equation} \label{eq:XW}
		{X_s = \int_t^s b_r^{\P} dr + \int_t^s \ah^{ 1/2}_r dW_r^{\P},
		~s \in [t, T], ~\Pb-a.s.,}
	\end{equation}
	where we extend the definition of $b^{\P}$ and $\ah$ on $\Omb$ as in \eqref{eq:extension_def}, and where we recall that $\widehat a^{1/2}$ has been defined in the Notations above.

\vspace{0.5em}
	\noindent Notice that when $\ah_r$ is non--degenerate $\P-a.s.$, for all $r \in [t, T]$, we can construct the Brownian motion $W^{\P}$ on $\Om$ as follows
	$${W_t^\P:=\int_0^t\widehat a_s^{-1/2}dX_s^{c,\P},\ t\in[0,T],\ \P-a.s.,}$$
	 and do not need to consider the above enlarged space equipped with an independent Brownian motion to construct $W^{\P}$.
	 
	\begin{Remark}[{On the choice of $\ah^{1/2}$}]
		{\rm The measurable map $a \longmapsto a^{ 1/2}$ is fixed throughout the paper.
		A first choice is to take $a^{ 1/2}$ as the unique non--negative symmetric square root of $a$ $($see {\it e.g.} Lemma 5.2.1 of {\rm\cite{stroock2007multidimensional}}$)$.
		One can also use the Cholesky decomposition to obtain $a^{ 1/2}$ as a lower triangular matrix.
		Finally, when $d = m+n$ for $m, n \in \N^{\ast}$, and $\ah$ has the specific structure of Remark \ref{rem:sde} below, one can take $a^{ 1/2}$ in the following way
		\be \label{eq:def_a12}
			a = \left( \begin{array}{cc}
			\sigma \sigma^\top & \sigma \\
			\sigma^\top & I_n 
			\end{array}
			\right)
			~\mbox{and}~
			a^{1/2} = \left( \begin{array}{cc}
			\sigma & 0 \\
			I_n & 0
			\end{array}
			\right),
			~\mbox{for some}~ \sigma \in \Mc_{m,n}.
		\ee}
	\end{Remark} 
	
\subsubsection{Conditioning and concatenation of probability measures}
	
	We also recall that for every probability measure $\P$ on $\Om$ and $\F-$stopping time $\tau$ taking value in $[0,T]$, there exists a regular conditional probability distribution (r.c.p.d. for short) $(\P^{\tau}_{\om})_{\om \in \Om}$
	 (see {\it e.g.} Stroock and Varadhan \cite{stroock2007multidimensional}), satisfying
	 
	 \vspace{0.5em}
		\noindent (i) For every $\om \in \Om$, $\mathbb{P}^{\tau}_{\om}$ is a probability measure on $(\Om,  \Fc_T)$.

	\vspace{0.5em}
		\noindent (ii) For every $ E \in \mathcal{F}_T $, the mapping $ \om \longmapsto \mathbb{P}^{\tau}_{\om}(E) $ is 
	$\mathcal{F}_\tau-$measurable.

	\vspace{0.5em}
		\noindent (iii) The family $ (\mathbb{P}^{\tau}_{\om})_{\om \in \Om} $ is a version of the conditional probability measure of $ \mathbb{P} $ on $ \mathcal{F}_{\tau}$, {\it i.e.}, for every integrable $ \mathcal{F}_T-$measurable random variable $ \xi $ we have
	$ \mathbb{E}^{\mathbb{P} } [ \xi | \mathcal{F}_{\tau}](\omega)=\mathbb{E}^{ \mathbb{P}^{\tau}_{\omega}} \big[\xi \big],$ for $\P-a.e.$ $\om\in\Om$.

	\vspace{0.5em}
		\noindent (iv) For every $\om \in \Om$, $ \mathbb{P}^{\tau}_{\om} (\Omega_{\tau}^{\om})=1$, 
	where $\Omega_{\tau}^{\om}\mathrel{\mathop:}= \big \lbrace \overline\omega \in \Omega : \ \overline\omega(s)=\om(s), \ 0\le s \le \tau(\om) \big \rbrace.$

	\vspace{0.5em}
	\noindent Furthermore, given some $\P$ and a family $(\Q_{\om})_{\om \in \Om}$ such that 
	$\om \longmapsto \Q_{\om}$ is $\Fc_{\tau}-$measurable and $\Q_{\om} (\Om_{\tau}^{\om}) = 1$ for all $\om \in \Om$, 
	one can then define a concatenated probability measure $\P \otimes_{\tau} \Q_{\cdot}$ by
	\b*
		\P \otimes_{\tau} \Q_{\cdot} \big[ A \big]
		:=
		\int_{\Om} \Q_{\om} \big[A \big] ~\P(d \om),
		\ \forall A \in \Fc_T.
	\e*



\subsection{Functional spaces and norms}
\label{subsec:space_norms}
We now give the spaces and norms which will be needed in the rest of the paper. 
	We are given a fixed family $\left (\Pc(t,\om) \right)_{(t,\om) \in [0,T] \x \Om}$
	of sets of probability measures on $(\Om, \Fc_T)$, where $\Pc(t,\om) \subset \Pc^W_t$.
	Fix some $t\in[0,T]$ and some $\omega\in\Omega$. In what follows, $\X:=(\mathcal X_s)_{t\leq s\leq T}$ will denote an arbitrary filtration on $\Om$, $\mathcal X$ an arbitrary $\sigma-$algebra on $\Om$, 
	and $\P$ an arbitrary element in $\Pc(t,\omega)$.
	Denote also by $\X_{\P}$ the $\P-$augmented filtration associated to $\X$.

\vspace{0.5em}

			\noindent $\bullet$ For $p \geq  1$, $\L^{p}_{t,\omega}(\Xc)$ (resp. $\L^p_{t,\omega}(\Xc,\P)$) denotes the space of all $\mathcal X-$measurable scalar random variable $\xi$ with
$$\No{\xi}_{\L^{p}_{t,\omega}}^p:=\underset{\mathbb{P} \in \mathcal{P}(t,\omega)}{\sup}\mathbb E^{\mathbb P}\left[|\xi|^p\right]<+\infty,\ \left(\text{resp. }\No{\xi}_{\L^{p}_{t,\omega}(\P)}^p:=\mathbb E^{\mathbb P}\left[|\xi|^p\right]<+\infty\right).$$

			\noindent $\bullet$ $\mathbb H^{p}_{t,\omega}(\X)$ (resp. $\mathbb H^p_{t,\omega}(\X,\P)$) denotes the space of all $\mathbb X-$predictable $\mathbb R^d-$valued processes $Z$, which are defined $\widehat a_s ds-a.e.$ on $[t,T]$, with
\begin{align*}
&\No{Z}_{\mathbb H^{p}_{t,\omega}}^p:=\underset{\mathbb{P} \in \mathcal{P}(t,\omega)}{\sup}\mathbb E^{\mathbb P}\left[\left(\int_t^T\No{\big(\ah_s^{1/2}\big)^\top Z_s}^2ds\right)^{\frac p2}\right]<+\infty, \\
& \left(\text{resp. } \No{Z}_{\mathbb H^{p}_{t,\omega}(\P)}^p:=\mathbb E^{\mathbb P}\left[\left(\int_t^T\No{ \big(\ah_s^{1/2}\big)^\top Z_s}^2ds\right)^{\frac p2}\right]<+\infty\right).
\end{align*}

			\noindent $\bullet$ $\mathbb M^p_{t,\omega}(\X,\P)$ denotes the space of all $(\mathbb X,\P)-$optional martingales $M$ with $\P-a.s.$ c\`adl\`ag paths on $[t,T]$, with $M_t=0$, $\P-a.s.$, and 
$$\No{M}_{\mathbb M^p_{t,\omega}(\P)}^p:=\E^\P\left[\left[M\right]_T^{\frac p2}\right]<+\infty.$$
Furthermore, we will say that a family $(M^\P)_{\P\in\mathcal P(t,\omega)}$ belongs to $\mathbb M^p_{t,\omega}((\X_\P)_{\P\in\Pc(t,\omega)})$ if, for any $\P\in\Pc(t,\omega)$, $M^\P\in\mathbb M^p_{t,\omega}(\X_\P,\P)$ and 
$$\underset{\mathbb P\in\mathcal P(t,\omega)}{\sup}\No{M^\P}_{\mathbb M^p_{t,\omega}(\P)}<+\infty.$$

			\noindent $\bullet$ $\mathbb I^p_{t,\omega}(\X,\P)$ (resp. $\mathbb I^{o,p}_{t,\omega}(\X,\P)$) denotes the space of all $\mathbb X-$predictable (resp. $\mathbb X-$optional) processes $K$ with $\P-a.s.$ c\`adl\`ag and non-decreasing paths on $[t,T]$, with $K_t=0$, $\P-a.s.$, and
$$\No{K}_{\mathbb I^p_{t,\omega}(\P)}^p:=\E^\P\left[K_T^{p}\right]<+\infty
	~~\mbox{(resp.}~ \No{K}_{\mathbb I^{o,p}_{t,\omega}(\P)}^p :=\E^\P\left[K_T^{p}\right] <+\infty \mbox{)}.
$$
We will say that a family $(K^\P)_{\P\in\mathcal P(t,\omega)}$ belongs to $\mathbb I^{p}_{t,\omega}((\X_\P)_{\P\in\Pc(t,\omega)})$ (resp. $\mathbb I^{o,p}_{t,\omega}((\X_\P)_{\P\in\Pc(t,\omega)})$) if, for any $\P\in\mathcal P(t,\omega)$, $K^\P\in\mathbb I^{p}_{t,\omega}(\X_\P,\P)$ (resp. $K^\P\in\mathbb I^{o,p}_{t,\omega}(\X_\P,\P)$) and 
$$
	\underset{\mathbb P\in\mathcal P(t,\omega)}{\sup}\No{K^\P}_{\mathbb I^p_{t,\omega}(\P)}<+\infty
	~\left(\mbox{resp.}~
	\underset{\mathbb P\in\mathcal P(t,\omega)}{\sup}\No{K^\P}_{\mathbb I^{o,p}_{t,\omega}(\P)}<+\infty
	\right).
$$

			\noindent $\bullet$ $\mathbb D^{p}_{t,\omega}(\X)$ (resp. $\mathbb D^p_{t,\omega}(\X,\P)$) denotes the space of all $\mathbb X-$progressively measurable $\mathbb R-$valued processes $Y$ with $\mathcal P(t,\omega)-q.s.$ (resp. $\P-a.s.$) c\`adl\`ag paths on $[t,T]$, with
$$\No{Y}_{\mathbb D^{p}_{t,\omega}}^p:=\underset{\mathbb{P} \in \mathcal{P}(t,\omega)}{\sup}\mathbb E^{\mathbb P}\left[\underset{t\leq s\leq T}{\sup}|Y_s|^p\right]<+\infty, \; \Bigg(\text{resp. }\No{Y}_{\mathbb D_{t,\omega}^{p}(\P)}^p:=\mathbb E^{\mathbb P}\left[\underset{t\leq s\leq T}{\sup}|Y_s|^p\right]<+\infty\Bigg).$$

			\noindent $\bullet$ Similarly, given a probability measure $\Pb$ and a filtration $\overline \X$ on the enlarged canonical space $\Omb$, we denote the corresponding spaces by
$\mathbb D^p_{t,\omega}(\overline \X,\Pb)$,
$\mathbb H^p_{t,\omega}(\overline \X,\Pb)$,
$\mathbb M^p_{t,\omega}(\overline \X,\Pb)$,... Furthermore, when $t=0$, there is no longer any dependence on $\omega$, since $\omega_0=0$, so that we simplify the notations by suppressing the $\omega-$dependence and write $\H^p_0(\X)$, $\H^p_0(\X,\P)$,... Similar notations are used on the enlarged canonical space.

\subsection{Assumptions}

	Let us  provide here a class of conditions which will be assumed throughout the paper.
	We shall consider a random variable $\xi : \Om \longrightarrow \R$ and a generator function 
	\begin{equation*}
		f: (t, \om, y, z, a, b)  \in [0,T] \x \Om \x \R \x \R^d \x \Sp \x \R^d
		 \longrightarrow  
		\R.
	\end{equation*}
	Define for simplicity
	$$
		\widehat{f}^\P_s(y,z):= f(s, X_{\cdot\wedge s},y,z,\widehat{a}_s, b^\P_s)
		\mbox{ and } 
		\widehat{f}^{\P,0}_s:= f(s, X_{\cdot\wedge s},0,0,\widehat{a}_s, b^\P_s).
	$$
	Recall that we are given a family $\left (\Pc(t,\om) \right)_{(t,\om) \in [0,T] \x \Om}$
	of sets of probability measures on $(\Om, \Fc_T)$, where $\Pc(t,\om) \subset \Pc^W_t$ 
	for all $(t,\om) \in [0,T] \x \Om$.
	Denote also $\Pc_t := \cup_{\om \in \Om} \Pc(t,\om)$.
	We make the following assumption on $\xi$, $f$ and the family $\left (\Pc(t,\om) \right)_{(t,\om) \in [0,T] \x \Om}$.

	\begin{Assumption} 
	\label{assum:main}

	$\mathrm{(i)}$ The random variable $\xi$ is $\Fc_T-$measurable,
		the generator function $f$ is jointly Borel measurable and such that for every $(t,\omega,y,y',z,z',a,b)\in[0,T]\times\Omega\times\R\times\R\times\R^d\times\R^d\times\Sp\times\R^d$,
				$$\abs{f(t,\omega,y,z,a,b)-f(t,\omega,y',z',a,b)}\leq C\left(\abs{y-y'}+\No{z-z'}\right),$$
		and for every fixed $(y,z,a,b)$, the map $(t,\om) \longmapsto f(t,\om, y,z,a,b)$ is $\F-$progressively measurable.

\vspace{0.5em}
		\noindent$\mathrm{(ii)}$ For the fixed constant $p > 1$, one has for every $(t,\omega) \in[0,T] \x \Omega$,
			\begin{equation} \label{eq:integrability}
				\sup_{\P \in \Pc(t,\omega)} 
				\E^{\P} \left[ 
					|\xi|^p
					+ 
					\int_t^T \big| f(s, X_{\cdot \wedge s}, 0, 0, \ah_s, b_s^\P) \big|^p ds 
				\right] 
				<
				+ \infty.
			\end{equation}

\vspace{0.5em}	
	\noindent	$\mathrm{(iii)}$ For every $(t,\om) \in [0,T] \x \Om$, one has $\Pc(t, \om) = \Pc(t, \om_{\cdot \wedge t})$ and
		$\P(\Om_t^{\om}) = 1$ whenever $\P \in \Pc(t, \om)$.
		The graph $[[ \Pc ]]$ of $\Pc$, defined by $[[ \Pc ]] := \{ (t, \om, \P): \P \in \Pc(t,\om) \}$, is  analytic in $[0,T] \x \Om \x \M_1$.

\vspace{0.5em}	
	\noindent	$\mathrm{(iv)}$ $\Pc$ is stable under conditioning, {\it i.e.} for every $(t,\om) \in [0,T] \x \Om$ and every $\P \in \Pc(t,\om)$ together with an $\F-$stopping time $\tau$ taking values in $[t, T]$,
		there is a family of r.c.p.d. $(\P_{\w})_{\w \in \Om}$ such that $\P_{\w} \in \Pc(\tau(\w), \w)$, for $\P-a.e.$ $\w \in \Om$.

\vspace{0.5em}	
	\noindent	$\mathrm{(v)}$ $\Pc$ is stable under concatenation, {\it i.e.} for every $(t,\om) \in [0,T] \x \Om$ and $\P \in \Pc(t,\om)$ together with a $\F-$stopping time $\tau$ taking values in $[t, T]$,
		let $(\Q_{\w})_{\w \in \Om}$ be a family of probability measures such that 
		$\Q_{\w} \in \Pc(\tau(\w), \w)$ for all $\w \in \Om$ and $\w \longmapsto \Q_{\w}$ is $\Fc_{\tau}-$measurable,
		then the concatenated probability measure $\P \ox_{\tau} \Q_{\cdot} \in \Pc(t,\om)$.
	\end{Assumption}

	\begin{Remark}{\rm
		Conditions \rmi and \rmii of Assumption \ref{assum:main} are very standard conditions
		to obtain existence and uniqueness to the standard {\rm {\rm BSDE}} with generator $f$ and terminal condition $\xi$.	The only change here is that \eqref{eq:integrability} takes into account that we are working with a whole family of measures, hence the supremum.	}
	\end{Remark}
	
	\begin{Remark} \label{remark:assum_meas_selec}{\rm
		Conditions $\mathrm{(iii)}$--\rmv of Assumption \ref{assum:main} will be essentially used to prove the dynamic programming principle for our control problem on the nonlinear kernels, and are the generic conditions used in such a setting $($see for instance {\rm\cite{nutz2013constructing}}, and the extension given in {\rm\cite{possamai2013robust}}$)$.
		Notice also that when $t = 0$, we have $\Pc_0 = \Pc(0, \om)$ since $\om_0 = 0$ for any $\om \in \Om$.

	\vspace{0.5em}	
	\noindent	
		In particular, let us consider the case where the sets $\Pc(t,\omega)$ are induced by controlled diffusion processes.
		Let $U$ be some $($non--empty$)$ Polish space,
		$\Uc$ denote the collection of all $U-$valued and $\F-$progressively measurable processes,
		$(\mu, \sigma) : [0,T] \x \Om \x U \longrightarrow \R^d \x \S^d$ be the drift and volatility coefficient function which are assumed to be such that for some constant $L > 0$,
		$\|(\mu, \sigma)(t, 0, u)\|  \le L$ and
		\begin{equation*}
			\No{ (\mu, \sigma) (t, \om, u) - (\mu, \sigma) (t', \om', u) }
			\le
			L \left( \sqrt{ |t - t'|} + \No{\om_{t \wedge \cdot} - \om'_{t ' \wedge \cdot} } \right).
		\end{equation*}
		Recall that the canonical process $X$ on the canonical space $\Om$ is a standard Brownian motion under the Wiener measure $\P_0$.
		We denote by $S^{t,\om, \nu}$ the unique (strong) solution to the {\rm SDE}
		\begin{equation*}
			S^{t, \om, \nu}_s 
			=
			\om_t 
			+ 
			\int_t^s \mu(r, S^{t, \om, \nu}, \nu_r) dr 
			+ 
			\int_t^s \sigma(r, S^{t, \om, \nu}, \nu_r) dX_r, 
			~s \in [t, T],
			~\P_0-a.s.,
		\end{equation*}
		with initial condition $S^{t,\om, \nu}_s = \om_s$ for all $s \in [0,t]$ and $\nu \in \Uc$.
		Then the collection $\Pc^{\Uc}(t,\om)$ of sets of measures defined by
		\begin{equation*}
			\Pc^{\Uc}(t, \om) 
			:=
			\big\{
				\P_0 \circ \big( S^{t, \om, \nu}\big)^{-1},\; \nu \in \Uc
			\big\}
		\end{equation*}
		satisfies Assumption \ref{assum:main} $\mathrm{(iii)-(v)}$
		$($see Theorem 3.5 and Lemma 3.6 in {\rm\cite{karoui2013capacities2}} or Theorem 2.4 and Proposition 2.2 in {\rm\cite{neufeld2013superreplication}} in a simpler context$)$.
		One can find more examples of $\Pc(t,\om)$ satisfying Assumption \ref{assum:main} $\mathrm{(iii)-(v)}$ in \cite{karoui2013capacities2}, which are induced from the weak/relaxed formulation of the controlled {\rm SDEs}.}
	\end{Remark}

\section{Stochastic control for a class of nonlinear stochastic kernels}\label{sec:1}

	For every $(t,\om) \in [0,T] \x \Om$ and $\P \in \Pc(t, \om)$, we consider the following {\rm BSDE}
	\begin{equation} \label{eq:BSDE}
		\Yc_s
		=
		\xi 
		- 
		\int_s^T f \left(r, X_{\cdot \wedge r}, \Yc_r,  (\ah_r^{1/2})^\top \Zc_r, \ah_r, b_r^\P \right) dr
		-
		 \left( \int_s^T \Zc_r \cdot  dX^{c,\P}_r \right)^\P
		-
		\int_s^T d\Mc_r,
		~\P-a.s.,
	\end{equation}
	where $( \int_s^T \Zc_r \cdot  dX^{c,\P}_r )^\P$ denotes the stochastic integrable of $\Zc$ w.r.t. $X^{c,\P}$ under the probability $\P$.
	Following El Karoui and Huang \cite{elkaroui1997general}, we define a solution to {\rm BSDE} \eqref{eq:BSDE} as a triple 
	$\left(\Yc^{\P}_s, \Zc^{\P}_s, \Mc^{\P}_s \right)_{s \in [t, T]}\in \mathbb D^p_{t,\omega}(\F^\P_+,\P)\times \mathbb H^p_{t,\omega}(\F^\P_+,\P)\times \mathbb M^p_{t,\omega}(\F^\P_+,\P)$ satisfying the equality \eqref{eq:BSDE} under $\P$.
	In particular, under the integrability condition of $\xi$ and $f$ in Assumption \ref{assum:main}, one has existence and uniqueness of the solution to {\rm BSDE} \eqref{eq:BSDE} (see Lemma \ref{lemma:deuxbsde} below).

\vspace{0.5em}
	\noindent One may also consider a stopping time $\tau$ and $\Fc^U_{\tau}-$measurable random variable $\zeta$, and consider the {\rm BSDE}, with terminal time $\tau$ and terminal condition $\zeta$,
	\begin{equation}\label{eq:bsde}
		\Yc_t
		=
		\zeta
			- 
		\int_t^{\tau} f \left(s, X_{\cdot \wedge s}, \Yc_s,  (\ah_s^{1/2})^\top \Zc_s, \ah_s, b^\P_s\right) ds
		-
		 \left( \int_t^{\tau}  \Zc_s \cdot dX^{c,\P}_s \right)^\P
		-
		\int_t^{\tau}  d\Mc_s.
	\end{equation}
	We will denote by $\Yc^{\P}_{\cdot}(\tau, \zeta)$ the $Y$ part of its solution whenever the above {\rm BSDE} \eqref{eq:bsde} is well--posed in the above sense, and set $\Yc^{\P}_{\cdot}(\tau, \zeta) = \infty$ otherwise.

\subsection{Control on a class of nonlinear stochastic kernels and the dynamic programming principle}

	Remember that under different $\P \in \Pc(t, \om)$, the law of the generating process $X$ is different,
	and hence the solution $\Yc^{\P}_t$ will also be different.
	We then define, for every $(t,\om) \in [0,T] \x \Om$,
	\begin{equation} \label{eq:def_Y}
		\Ych_t(\om) := \sup_{\P \in \Pc(t,\om)} \E^{\P} \Big[ \Yc^{\P}_t \Big].
	\end{equation}
	Notice that $\Yc_t^{\P}$ is only $\Fc^{\P}_{t+}$-measurable, we hence consider its expectation before taking the supremum.
	Moreover, since $\Pc(t,\om)$ depends only on $\om_{t \wedge \cdot}$, it follows that $\Ych_t(\om)$ also depends  only on $\om_{t \wedge \cdot}$.

	\vspace{0.5em}

	\noindent Our first main result is the following dynamic programming principle.
	\begin{Theorem}\label{theo:main_dpp}
		Suppose that Assumption \ref{assum:main} holds true.
		Then for all $(t,\om)\in [0,T] \x \Om$, one has $\Ych_t(\om) = \Ych_t(\om_{t \wedge \cdot})$,
		and $(t, \om) \longmapsto \Ych_t(\om)$ is $\Bc([0,T]) \otimes \Fc_T-$universally measurable.
		Consequently, for all $(t,\om) \in [0,T] \x \Om$ and $\F-$stopping time $\tau$ taking values in $[t,T]$,
		the variable $\Ych_{\tau}$ is $\Fc^U_{\tau}-$measurable. 
		Moreover, one has for any $\P\in\Pc(t,\omega)$,
		\begin{equation*}
			\mathbb E^\P\Big[ \big| \widehat \Yc_\tau \big|^p\Big] < +\infty
			~\mbox{and}~
			\Ych_t(\om) 
			=
			\sup_{\P \in \Pc(t, \om)} \E^{\P} \left[ \Yc_t^{\P} \big(\tau, \Ych_{\tau} \big) \right].
		\end{equation*}
		\end{Theorem}

	\begin{Remark}\label{rem:sde}{\rm 
		In some contexts, the sets $\Pc(t,\om)$ are defined as the collections of probability measures induced by a family of controlled diffusion processes (recall Remark \ref{remark:assum_meas_selec}).
		For example, let $\C_1$ $($resp. $\C_2)$ denote the canonical space of all continuous paths $\om^1$ in $\C([0,T], \R^n)$ $($resp. $\om^2$ in $\C([0,T], \R^{m}))$ such that $\om^1_0 = 0$ $($resp. $\om^2_0 = 0)$, with canonical process $B$, canonical filtration $\F^1$, and let $\P_0$ be the corresponding Wiener measure.
		Let $U$ be a Polish space, $(\mu, \sigma): [0,T] \x \C_1 \x U \longrightarrow \R^n \x \Mc_{n,m}$ be the coefficient functions,
		then, given $(t,\om^1) \in [0,T] \x \C_1$, we denote by $\Jc(t, \om^1)$ the collection of all terms
		$$
			\alpha := \big( \Om^{\alpha}, \Fc^{\alpha}, \P^{\alpha}, \F^{\alpha} = (\Fc^{\alpha}_t)_{t \ge 0}, W^{\alpha}, (\nu^{\alpha}_t)_{t \ge 0}, X^{\alpha} \big),
		$$
		where $\big( \Om^{\alpha}, \Fc^{\alpha}, \P^{\alpha}, \F^{\alpha} \big)$ is a filtered probability space,
		$W^{\alpha}$ is a $\F^{\alpha}-$Brownian motion,
		$\nu^{\alpha}$ is a $U-$valued $\F^{\alpha}-$predictable process
		and $X^{\alpha}$ solves the SDE $($under some appropriate additional conditions on $\mu$ and $\sigma)$, with initial condition $X^{\alpha}_s = \om^1_s$ for all $s \in [0,t]$,
		\begin{align*}
			X_s^{\alpha} ~=~ \om^1_t + \int_t^s \mu(r, X^{\alpha}_{r \wedge \cdot}, \nu^{\alpha}_r) dr
			+ \int_t^s \sigma(r, X^{\alpha}_{r \wedge \cdot}, \nu^{\alpha}_r) d W^{\alpha}_r,
			~s \in [t, T],
			~\P^{\alpha}-a.s.
		\end{align*}
In this case, one can let $d = m+n$ so that $\Om = \C_1 \x \C_2$ 
		and define $\Pc(t,\om)$ for $\om = (\om^1, \om^2)$ as the collection of all probability measures induced by $(X^{\alpha}, B^{\alpha})_{\alpha \in \Jc(t, \om^1)}$.
		
		\vspace{0.5em}\noindent
		Then, with the choice of $\ah^{1/2}$ as in \eqref{eq:def_a12}, one can recover $\sigma$ from it directly, which may be useful for some applications.
		Moreover, notice that  $\Pc(t, \om)$ depends only on $(t, \om^1)$ for $\om = (\om^1, \om^2)$, then the value $\Ych_t(\om)$ in \eqref{eq:def_Y} depends also only on $(t, \om^1)$.}
	\end{Remark}

\subsection{Proof of Theorem \ref{theo:main_dpp}}

\subsubsection{An equivalent formulation on the enlarged canonical space}

	We would like formulate {\rm BSDE} \eqref{eq:BSDE} on the enlarged canonical space in an equivalent way.
	Remember that $\Omb := \Om \x \Om'$ and for a probability measure $\P$ on $\Om$, we define $\Pb := \P \otimes \P_0$.
	Then a $\P-$null event on $\Om$ becomes a $\Pb-$null event on $\Omb$ if it is considered in the enlarged space.
	Let $\pi: \Om \x \Om' \longrightarrow \Om$ be the projection operator defined by $\pi(\om, \om') := \om$, for any $(\om,\om')\in\Omb$.

\begin{Lemma} \label{lemm:enlargement}
	Let  $A \subseteq \Om$ be a subset in $\Om$.
	Then saying that $A$ is a $\P-$null set is equivalent to saying that $\{ \omb ~: \pi(\omb) \in A\}$ is a $\Pb := \P \otimes \P_0-$null set.
\end{Lemma}
\proof
	For $A \subseteq \Om$, denote $\overline A := \{ \omb ~: \pi(\omb) \in A \} = A \x \Om'$.
	Then by the definition of the product measure, it is clear that
	$$
		\P(A) = 0 \Longleftrightarrow \P \otimes \P_0(\overline A) = 0,
	$$
	which concludes the proof.
	\qed

\vspace{0.5em}
	\noindent We now consider two BSDEs on the enlarged canonical space, w.r.t. two different filtrations. The first one is the following {\rm BSDE} on $(\Omb, \Fcb^X_T, \Pb)$ w.r.t the filtration $\Fb^{X,\Pb}$
	\begin{equation} \label{eq:BSDE2}
		\Ycb_s
		=
		\xi(X_{\cdot})
		- 
		\int_s^T f \left(r, X_{\cdot \wedge r}, \Ycb_r, (\ah_r^{1/2})^\top \Zcb_r, \ah_r, b_r^\P \right) dr
		-
		 \left( \int_s^T  \Zcb_r \cdot dX^{c,\P}_r \right)^\Pb
		-
		\int_s^T d\Mcb_r,
	\end{equation}
	a solution being a triple $(\Ycb^{\Pb}_s, \Zcb^{\Pb}_s, \Mcb^{\Pb}_s )_{s \in [t, T]}\in \mathbb D^p_{t,\omega}(\Fb^{X,\Pb}_+,\Pb)\times \mathbb H^p_{t,\omega}(\Fb^{X,\Pb}_+,\Pb)\times \mathbb M^p_{t,\omega}(\Fb^{X,\Pb}_+,\Pb)$ satisfying \eqref{eq:BSDE2}.
	Notice that in the enlargement, the Brownian motion $B$ is independent of $X$, so that the above {\rm BSDE} \eqref{eq:BSDE2} is equivalent to {\rm BSDE} \eqref{eq:BSDE} (see Lemma \ref{lemma:deuxbsde} below for a precise statement and justification).

	\vspace{0.5em}
	
	\noindent We then introduce a second {\rm BSDE} on the enlarged space $(\Omb, \Fcb_T, \Pb)$, w.r.t. the filtration $\Fb$,
	\begin{equation} \label{eq:BSDE3}
		{\widetilde \Yc_s
		=
		\xi(X_{\cdot})
		- 
		\int_s^T f \left(r, X_{\cdot \wedge r}, \widetilde \Yc_r, (\widehat{a}_r^{1/2})^\top \widetilde \Zc_r, \ah_r, b_r^\P \right) dr
		-
		 \left( \int_s^T  \widetilde\Zc_r \cdot \widehat{a}_r^{1/2} dW^{\P}_r \right)^\Pb
		-
		\int_s^T d \widetilde \Mc_r,}
	\end{equation}	
a solution being a triple $( \widetilde\Yc^{\Pb}_s, \widetilde \Zc^{\Pb}_s, \widetilde \Mc^{\Pb}_s )_{s \in [t, T]}\in \mathbb D^p_{t,\omega}(\Fb^{\Pb}_+,\Pb)\times \mathbb H^p_{t,\omega}(\Fb^{\Pb}_+,\Pb)\times \mathbb M^p_{t,\omega}(\Fb^{\Pb}_+,\Pb)$ satisfying \eqref{eq:BSDE3}.

	\begin{Lemma}\label{lemma:deuxbsde}
		Let $(t,\omega) \in [0,T]\times\Omega$, $\P \in \Pc(t,\omega)$ and $\Pb := \P \otimes \P_0$, 
		then each of the three {\rm BSDEs} \eqref{eq:BSDE}, \eqref{eq:BSDE2} and \eqref{eq:BSDE3} has a unique solution, denoted respectively by $\left(\Yc, \Zc, \Mc \right)$, $\left(\Ycb, \Zcb, \Mcb\right)$ and $(\widetilde \Yc, \widetilde \Zc, \widetilde\Mc )$.
		Moreover, their solution coincide in the sense that there is some functional 
		$$\Psi: = (\Psi^Y, \Psi^Z, \Psi^M) : [t,T] \x \Om \longrightarrow \R \x \R^d \x \R,$$ 
		such that $\Psi^Y$ and $\Psi^M$ are $\F_+-$progressively measurable and $\P-a.s.$ c\`adl\`ag, $\Psi^Z$ is $\F-$predic-table, and
		$$
			\Yc_s = \Psi^Y_s,
			~ \Zc_r =  \Psi^Z_r, \ \ah_rdr-a.e. \ \text{on }[t,s],
			~\Mc_s = \Psi^M_s,
			~\mbox{for all}~ s \in [t,T],~ \P-a.s.,
		$$
		\vspace{-0.5em}
		$$
			\Ycb_s = \widetilde \Yc_s = \Psi^Y_s(X_\cdot),
			~
			 \Zcb_r =  \widetilde \Zc_r 
			= \Psi^Z_r(X_\cdot), \ \ah_rdr-a.e. \ \text{on }[t,s],
			~
			\Mcb_s = \widetilde \Mc_s = \Psi^M_s(X_\cdot),
		$$
		for all $s \in [t,T]$, $\Pb-a.s.$
	\end{Lemma}
	\proof \rmi The existence and uniqueness of a solution to \eqref{eq:BSDE3} is a classical result,
	we can for example refer to Theorem 4.1 of \cite{bouchard2015unified}.
	Then it is enough to show that the three BSDEs share the same solution in the sense given in the statement.
	Without loss of generality, we assume in the following $t = 0$.

	\vspace{0.5em}
	
	\rmii We next show that \eqref{eq:BSDE2} and \eqref{eq:BSDE3} admit the same solution in $(\Omb, \Fcb^{\Pb}_T, \Pb)$.
	Notice that a solution to \eqref{eq:BSDE2} is clearly a solution to \eqref{eq:BSDE3} by \eqref{eq:XW}.
	We then show that a solution to \eqref{eq:BSDE3} is also a solution to \eqref{eq:BSDE2}.
	
	\vspace{0.5em}
	
	\noindent Let $\zeta: \Omb \longrightarrow \R$ be a $\Fcb^{X,\Pb}_T-$measurable random variable, which admits a unique martingale representation
	\begin{equation} \label{eq:mart_repres}
		\zeta = \E^{\Pb} [ \zeta ] + \int_0^T \Zcb^{\zeta}_s \cdot dX^{c,\P}_s + \int_0^T \Mcb^{\zeta}_s,
	\end{equation}
	w.r.t. the filtration $\Fb^{X,\Pb}_+$.
	Since $B$ is independent of $X$ in the enlarged space, and since $X$ admits the same semi-martingale triplet of characteristics in both space, the above representation \eqref{eq:mart_repres} w.r.t. $\Fb^{X,\Pb}_+$ is the same as the one w.r.t. $\Fb^{\Pb}_+$, which are all unique up to a $\Pb-$evanescent set.
	Remember now that the solution of {\rm BSDE} \eqref{eq:BSDE3} is actually obtained as an iteration of the above martingale representation (see e.g. Section \ref{subsec:BSDE_construction} below). Therefore, a solution to \eqref{eq:BSDE3} is clearly a solution to \eqref{eq:BSDE2}.
	
	\vspace{0.5em}

	\rmiii We now show that a solution $(\Ycb, \Zcb, \Mcb)$ to \eqref{eq:BSDE2} induces a solution to \eqref{eq:BSDE}.
	Notice that $\Ycb$ and $\Mcb$ are $\Fb^{X,\Pb}_+-$optional, and $\Zcb$ is $\Fb^{X,\P}_+-$predictable,
	then (see e.g. Lemma 2.4 of \cite{soner2013dual} and Theorem IV.78 and Remark IV.74 of \cite{dellacherie1978probabilities}) there exists a functional $(\Psib^Y, \Psib^Z, \Psib^M) : [0,T] \x \Omb \longrightarrow \R \x \R^d \x \R$ such that
	$\Psib^Y$ and $\Psib^M$ are $\Fb^{X}_+-$progressively mesurable and $\Pb-a.s.$ c\`adl\`ag, $\Psib^Z$ is $\Fb^X-$predictable,
	and $\Ycb_t = \Psib^X_t$, $\Zcb_t = \Psib^Z_t$ and $\Mcb_t = \Psib^M_t$, for all $t \in [0,T]$, $\Pb-a.s.$
	Define 
	$$
		(\Psib^{Y,0} (\om), \Psib^{Z,0}(\om), \Psib^{M,0}(\om)) 
		~:=~
		(\Psib^{Y} (\om, \0), \Psib^{Z}(\om, \0), \Psib^{M}(\om, \0)),
	$$
	where $\0$ denotes the path taking value $0$ for all $t \in [0,T]$.
	
	\vspace{0.5em}
	\noindent Since $(\Psib^Y, \Psib^Z, \Psib^M)$ are $\Fb^X-$progressively measurable,
	the functions $(\Psib^{Y,0}, \Psib^{Z,0}, \Psib^{M,0})$ are then $\Fb-$progressively measurable, and it is easy to see that they provide a version of a solution to \eqref{eq:BSDE} in $(\Om, \Fc^{\P}_T, \P)$.

	\vspace{0.5em}
	
	{\rm (iv)}  Finally, let $(\Yc, \Zc, \Mc)$ be a solution to \eqref{eq:BSDE}, then there is a function 
	$(\Psi^Y, \Psi^Z, \Psi^M) : [0,T] \x \Om \longrightarrow \R \x \R^d \x \R$ such that
	$\Psi^Y$ and $\Psi^M$ are $\F_+-$progressively measurable and $\P-a.s.$ c\`adl\`ag, $\Psi^Z$ is $\F-$predictable,
	and $\Yc_t = \Psi_t$, $\Zc_t = \Psi^Z_t$ and $\Mc_t = \Psi^M_t$, for all $t \in [0,T]$, $\P-a.s.$
	Since $\Pb := \P \otimes \P_0$, it is easy to see that $(\Psi^Y, \Psi^Z, \Psi^M)$ is the required functional in the lemma.
	\qed

\vspace{0.5em}

\noindent The main interest of Lemma \ref{lemma:deuxbsde} is that it allows us, when studying the {\rm BSDE} \reff{eq:BSDE}, to equivalently work with the {\rm BSDE} \reff{eq:BSDE3}, in which the Brownian motion $W^\P$ appears explicitly. This will be particularly important for us when using linearisation arguments. Indeed, in such type of arguments, one usually introduces a new probability measure equivalent to $\P$. But if we use formulation \reff{eq:BSDE}, then one must make the inverse of $\widehat a$ appear explicitly in the Radon--Nykodym density of the new probability measure. Since such an inverse is not always defined in our setting, we therefore take advantage of the enlarged space formulation to bypass this problem.

\subsubsection{An iterative construction of the solution to BSDE \eqref{eq:BSDE}}
\label{subsec:BSDE_construction}

	In preparation of the proof of the dynamic programming principle for the control problem in Theorem \ref{theo:main_dpp}, 
	let us first recall the classical construction of the $\Yc^{\P}$ part of the solution to the {\rm BSDE} \eqref{eq:BSDE} 
	under some probability $\P \in \Pc(t,\omega)$ using Picard iterations. 
	Let us first  define for any $m\geq 0$
	$$\xi^m:=(\xi \wedge m) \vee (-m),\ \ f^m(t,\omega,y,z,a,b):=(f(t,\omega,y,z,a,b)\wedge m)\vee(-m).$$
	\begin{itemize}
		\item[$\mathrm{(i)}$] First, let $\Yc^{\P,0,m}_s \equiv 0$ and $\Zc^{\P,0,m}_s \equiv 0$, for all $s \in [t,T]$.
		
		\item[$\mathrm{(ii)}$] Given a family of $\F_+-$progressively measurable processes
		$ \big(\Yc^{\P,n,m}_s, \Zc^{\P,n,m}_s \big)_{s \in [t,T]}$, we let
		\begin{equation} \label{eq:Ycb}
			\Ycb^{\P, n+1,m}_s
			:=
			\E^{\P}
			\left[ \left.\xi^m - \int_s^T f^m(r, X_{\cdot \wedge r}, \Yc^{\P,n,m}_r, (\widehat a_r^{1/2})^\top \Zc^{\P,n,m}_r, \ah_r, b^\P_r) dr  \right| \Fc_s \right],\ \P-a.s.
		\end{equation}
		
		\item[$\mathrm{(iii)}$] Let $\Yc^{\P, n+1,m}$ be a right--continuous modification of $\Ycb^{\P,n+1,m}$ defined by
		\begin{equation} \label{eq:lim_yb}
			\Yc^{\P, n+1,m}_s 
			:=
			\limsup_{\Q \ni r \downarrow s} \Ycb^{\P, n+1,m}_r,\ \P-a.s.
		\end{equation}
		
		\item[$\mathrm{(iv)}$] Notice that $\Yc^{\P, n+1,m}$ is a semi--martingale under $\P$.
		 Let $\langle \Yc^{\P, n+1,m}, X  \rangle^{\P}$ be the predictable quadratic covariation of the process $\Yc^{\P, n+1,m}$ and $X$ under $\P$.
		Define
		\be \label{eq:def_zc}
			{\Zc^{\P, n+1,m}_s
			:=
			\ah_s^{\oplus}  \bigg( \limsup_{\Q \ni \eps \downarrow 0} 
			\frac{\langle \Yc^{\P, n+1,m}, X \rangle^{\P}_s - \langle \Yc^{\P, n+1,m}, X \rangle^{\P}_{s-\eps}}{\eps} \bigg)
			.} 
		\ee
		
		\item[$\mathrm{(v)}$] Notice that the sequence $(\Yc^{\P,n,m})_{n \ge 0}$ is a Cauchy sequence for the norm 
		$$\No{(Y,Z)}^2_{\alpha}:=\E^\P\left[\int_0^Te^{\alpha s}\abs{Y_s}^2ds\right]^2+\E^\P\left[\int_0^Te^{\alpha s}\No{(\widehat a_s^{1/2})^\top Z_s}^2ds\right]^2,$$
		for $\alpha$ large enough. Indeed, this is a consequence of the classical estimates for BSDEs, for which we refer to Section $4$ of \cite{bouchard2015unified}\footnote{Notice here that the results of \cite{bouchard2015unified} apply for BSDEs of the form \eqref{eq:BSDE3} in the enlarged space. However, by Lemma \ref{lemma:deuxbsde}, this implies the same convergence result in the original space.}.
		Then by taking some suitable sub-sequence $(n^{\P,m}_k)_{k \ge 1}$, 
		we can define 
				\begin{equation*}
			\Yc^{\P,m}_s := \limsup_{k \to \infty} \Yc^{\P, n^{\P,m}_k,m}_s.
		\end{equation*}
		
		\item[$\mathrm{(vi)}$] Finally, we can again use the estimates given in \cite{bouchard2015unified} (see again Section $4$) to show that the sequence $(\Yc^{\P,m})_{m\geq 0}$ is a Cauchy sequence in $\mathbb D^p_0(\overline{\F}_+^\P,\P)$, so that by taking once more a suitable subsequence $(m^{\P}_k)_{k \ge 1}$, 
		we can define the solution to the {\rm BSDE} as
				\be \label{eq:YcP}
			\Yc^{\P}_s := \limsup_{k \to \infty} \Yc^{\P, m^{\P}_k}_s.
		\ee

	\end{itemize}

\subsubsection{On the measurability issues of the iteration}

	Here we show that the iteration in Section \ref{subsec:BSDE_construction} can be taken in a measurable way w.r.t. the reference probability measure $\P$,
	which allows us to use the measurable selection theorem to derive the dynamic programming principle.

	\begin{Lemma} \label{lemm:QVariation}
		Let $\Pc$ be a measurable set in $\M_1$,
		$(\P, \om, t) \longmapsto H^{\P}_t(\om)$ be a measurable function such that for all $\P \in \Pc$,
		$H^{\P}$ is right-continuous, $\F_+-$adapted and a $(\P,\F^{\P}_+)-$semi-martingale.
		Then there is a measurable function $(\P, \om, t) \longmapsto \Hlr^{\P}_t(\om)$
		such that for all $\P \in \Pc$, $\Hlr^{\P}$ is right--continuous, $\F_+-$adapted and $\F^{\P}_+-$predictable, and
		\b*
			\Hlr^{\P}_{\cdot} 
			~\mbox{is the predictable quadratic variation of the semi--martingale}
			~H^{\P}
			~\mbox{under}~ \P.
		\e*
	\end{Lemma}
	\proof
	\rmi For every $n\geq 1$, we define the following sequence of random times
	\begin{equation}\label{stoptimes}
	\begin{cases}
		\displaystyle\tau_0^{\P,n}(\omega):=0,\ \omega\in\Omega,\\[0.5em]
		\displaystyle\tau_{i+1}^{\P,n}(\omega):=\inf\left\{t\geq \tau_i^{\P, n} (\omega),\ \abs{H^{\P}_t(\omega)-H^{\P}_{\tau^{\P, n}_i}(\omega)}\geq2^{-n}\right\}\wedge T,\ \omega\in\Omega,\ i\geq 0.
	\end{cases}
	\end{equation}
	We notice that the $\tau^{\P, n}_i$ are all $\F_+-$stopping times since the $H^{\P}$ are right--continuous
	and $\F_+-$adapted.
	We then define
	\begin{equation}\label{quadvar}
		[H^{\P}]_{\cdot} (\om) 
		:=
		\underset{n\rightarrow+\infty}{\limsup}
		~\sum_{i\geq 0}
			\left(H^{\P}_{\tau^{\P,n}_{i+1}\wedge\cdot}(\omega)
				- H^{\P}_{\tau^{\P,n}_{i}\wedge \cdot}(\omega)
			\right)^2.
	\end{equation}
	It is clear that $(\P, \om, t) \longmapsto [H^{\P}]_t(\om)$ is a measurable function,
	and for all $\P \in \Pc$, $[H^{\P}]$ is non-decreasing, $\F_+-$adapted and $\F^{\P}_+-$optional.
	Then, it follows by Karandikar \cite{karandikar1995pathwise} that $[H^{\P}]$ coincides with the quadratic variation of the semi--martingale $H^{\P}$ under $\P$.
	Moreover, by taking its right limit {over rational time instants}, we can choose $[H^{\P}]$ to be right continuous.

	\vspace{0.5em}

	\rmii Finally, using Proposition 5.1 of Neufeld and Nutz \cite{neufeld2014measurability}, we can then construct a process $\Hlr^{\P}_t(\om)$ satisfying the required conditions.
	\qed

	\vspace{0.5em}

\noindent 	Notice that the construction above can also be carried out for the predictable quadratic covariation 
	$\langle H^{\P,1},H^{\P,2} \rangle^{\P}$, by defining it through the polarisation identity
	\be \label{eq:covariation}
		\langle H^{\P,1} , H^{\P,2} \rangle^{\P} 
		:=
		\frac14 \left(\langle H^{\P,1} + H^{\P,2} \rangle^{\P}
		- \langle H^{\P,1} - H^{\P,2}  \rangle^{\P} \right),
	\ee
	for all measurable functions $H^{\P,1}_t(\om)$ and $H^{\P,2}_t(\om)$ satisfying the conditions in Lemma \ref{lemm:QVariation}. We now show that the iteration in Section \ref{subsec:BSDE_construction} can be taken in a measurable way w.r.t. $\P$,
	which provides a key step for the proof of Theorem \ref{theo:main_dpp}.

	\begin{Lemma}
		Let $m, n > 0$ be fixed, and let
		$ (s,\om, \P) \longmapsto (\Yc^{\P,n,m}_s(\om),  \Zc^{\P,n,m}_s(\om) )$
		be a measurable map such that for every $\P \in \Pc_t$, 
		$\Yc^{\P,n,m}$ is right--continuous, $\F_+-$adapted and $\F^{\P}_+-$optional,
		$ \Zc^{\P,n,m}$ is $\F_+-$adapted and $\F^{\P}_+-$predictable.
		We can choose a measurable map
		$ (s,\om, \P) \longmapsto \big(\Yc^{\P,n,m}_s(\om), \Zc^{\P,n,m}_s(\om) \big)$
		s.t. for every $\P \in \Pc_t$, 
		$\Yc^{\P,n+1,m}$ is right--continuous, $\F_+-$adapted and $\F^{\P}_+-$optional,
		and $ \Zc^{\P,n+1,m}$ is $\F_+-$adapted and $\F^{\P}_+-$predictable.
	\end{Lemma}
	\proof 		
		\rmi First, using Lemma 3.1 of Neufeld and Nutz \cite{neufeld2014measurability},
		there is a version of $( \Ycb^{\P, n+1,m} )$ defined by \eqref{eq:Ycb},
		such that $(\P, \om) \longmapsto \Ycb^{\P, n+1,m}_s$ is {$\Bf \ox \Fc_s-$}measurable
		for every $s \in [t, T]$.

	\vspace{0.5em}
		\rmii Next, we notice that the measurability is not lost by taking the limit along a countable sequence.
		Then with the above version of $( \Ycb^{\P, n+1,m} )$, 
		it is clear that the family $(\Yc^{\P,n+1,m}_s(\om))$ 
		defined by \eqref{eq:lim_yb} is measurable in $(s, \om, \P)$,
		and for all $\P \in \Pc_t$, $\Yc^{\P, n+1,m}$ is $\F_+-$adapted and $\F^{\P}_+-$optional.

	\vspace{0.5em}
		\rmiii Then using Lemma \ref{lemm:QVariation} as well as the definition of the quadratic covariation in \eqref{eq:covariation},
		it follows that there is a measurable function 
		$$(s, \om, \P) \longmapsto \langle \Yc^{\P, n+1,m}, X  \rangle^{\P}_s(\om), $$ 
		such that for every $\P \in \Pc_t$, 
		$\langle \Yc^{\P, n+1,m}, X  \rangle^{\P}$ is right--continuous, $\F_+-$adapted 
		and coincides with the predictable quadratic covariation of $\Yc^{\P, n+1,m}$ and $X$ under $\P$.
		
	\vspace{0.5em}
		\rmiv Finally, with the above version of $\left(\langle \Yc^{\P, n+1,m}, X  \rangle^{\P}\right)$, 
		it is clear that the family $( \Zc^{\P, n+1,m}_s(\om))$ defined by \eqref{eq:def_zc} is measurable in 
		$(s, \om, \P)$ and for every $\P \in \Pc_t$,
		$\Zc^{\P, n+1,m}$ is $\F_+-$adapted and $\F_+^{\P}-$predictable.
	\qed

	\begin{Lemma}
		For every $\P \in \Pc_t$, there is a right--continuous, 
		$\F^{\P}_+-$martingale $\Mc^{\P, n+1,m}$ orthogonal to $X$ under $\P$,
		such that $\P-a.s.$
		\begin{align} \label{eq:YcPn1}
			\Yc^{\P,n+1,m}_t
			=&\ 
			\xi^m
			-
			\int_t^T f^m(s, X_{\cdot \wedge s}, \Yc^{\P,n,m}_s, (\widehat a_s^{1/2})^\top \Zc^{\P,n,m}_s, \ah_s, b^\P_s) ds
			-\int_t^T d\Mc^{\P,n+1,m}_s
			 \nonumber \\
			&-
			\left(\int_t^T \Zc^{\P,n+1,m}_s \cdot dX^{c,\P}_s\right)^\P .
		\end{align}
	\end{Lemma}
	\proof Using Doob's upcrossing inequality, 
	the limit $\lim_{r \downarrow s} \Ycb^{\P, n+1,m}_r$ exists $\P-$almost surely, for every $\P \in \Pc_t$.
	In other words, $\Yc^{\P, n+1,m}$ is a version of the right--continuous modification of the semi--martingale $\Ycb^{\P, n+1,m}$.
	Using  the martingale representation, it follows that there is a right--continuous, 
	$\F^{\P}_+-$martingale $\Mc^{\P, n+1,m}$ orthogonal to $X$ under $\P$, and an $\F^{\P}_+-$predictable process $\Zch^{\P, n+1, m}$ such that
\begin{align*}
		\Yc^{\P,n+1,m}_t
			=&\
			\xi^m
			-
			\int_t^T f^m(s, X_{\cdot \wedge s}, \Yc^{\P,n,m}_s, (\widehat a_s^{1/2})^\top \Zc^{\P,n,m}_s, \ah_s, b^\P_s) ds
			-\int_t^T d\Mc^{\P,n+1,m}_s
			 \nonumber \\
			&-
			\left(\int_t^T \Zch^{\P,n+1,m}_s \cdot dX^{c,\P}_s\right)^\P .
	\end{align*}
	In particular, $\Zch^{\P, n+1, m}$ satisfies 
	$$\langle \Yc^{\P, n+1, m}, X^{c, \P} \rangle_t = \int_0^t \ah_s \Zch^{\P, n+1, m}_s ds, \;\P-a.s.$$
	Besides, by the definition of $\Zc^{\P, n+1,m}$ in \eqref{eq:def_zc}, one has
	$$\int_0^t \ah_s \Zc^{\P, n+1, m}_s ds = \langle \Yc^{\P, n+1, m}, X^{c, \P} \rangle_t = \int_0^t \ah_s \Zch^{\P, n+1, m}_s ds, ~\P-\mbox{a.s.}$$
	It follows that 
\begin{align*}
		&\int_0^t  \big\| (\ah_s^{ 1/2})^\top (\Zc^{\P, n+1, m}_s  - \Zch^{\P, n+1, m}_s)  \big\|^2 ds  \\
		&=
		\int_0^t  (\Zc^{\P, n+1, m}_s  - \Zch^{\P, n+1, m}_s)^\top \ah_s (\Zc^{\P, n+1, m}_s  - \Zch^{\P, n+1, m}_s) ds 
		= 0, ~\P-\mbox{a.s}.
\end{align*}
	Hence \eqref{eq:YcPn1} holds true.
	\qed

	\begin{Lemma} \label{lemma:measurability_BSDE}
		There are families of subsequences $(n^{\P,m}_k, ~ k \ge 1)$ and 
		$(m^{\P}_i, ~ i \ge 1)$ such that 
		the limit 
		$ \Yc^{\P}_s(\om) = \lim_{i \to \infty} \lim_{k \to \infty} \Yc^{\P,n^{\P,m}_k, m^{\P}_i}_s$ exists for all $s \in [t,T]$, 
		$\P-$almost surely, 
		for every $\P \in \Pc_t$, and $(s, \om, \P) \longmapsto \Yc^{\P}_s(\om)$ is a measurable function.
		Moreover, $\Yc^{\P}$ provides a solution to the {\rm BSDE} \eqref{eq:BSDE} for every $\P \in \Pc_t$.
	\end{Lemma}
	\proof 
	By the conditions in \eqref{eq:integrability}, $(\Yc^{\P,n,m}, \Zc^{\P,n,m})_{n \ge 1}$ provides a Picard iteration under the $(\P, \beta)-$norm, for $\beta > 0$ large enough (see e.g. Section 4 of \cite{bouchard2015unified}\footnote{\label{foot}Again, we remind the reader that one should first apply the result of \cite{bouchard2015unified} to the corresponding Picard iteration of \eqref{eq:BSDE3} in the enlarged space and then use Lemma \ref{lemma:deuxbsde} to go back to the original space.}), 
	defined by
		$$
			||{\varphi}||^{ 2}_{\P, \beta} 
			:= \E ^{\P} \left[ \sup_{t\leq s\leq T}e^{\beta s}|\varphi_s|^2 \right].
		$$
	Hence, $\Yc^{\P,n,m}$ converges (under the $(\P, \beta)-$norm) to some process $\Yc^{\P,m}$ as $n \longrightarrow \infty$, which solves the {\rm BSDE} \eqref{eq:BSDE} with the truncated terminal condition $\xi^m$ and truncated generator $f^m$.
	Moreover, by the estimates in Section 4 of \cite{bouchard2015unified} (see again Footnote \ref{foot}), $(\Yc^{\P,m})_{m \ge 1}$ is a Cauchy sequence
	in $\D^p_{t,\omega}(\overline{\F}^\P_+,\P)$.
	Then using Lemma 3.2 of \cite{neufeld2014measurability}, we can find two families of subsequences $(n^{\P,m}_k, k\ge 1, \P \in \Pc_t)$ 
	and $(m^{\P}_i, i\ge 1, \P \in \Pc_t)$ satisfying the required properties.
	\qed

\subsubsection{End of the proof of Theorem \ref{theo:main_dpp}}

	Now we can complete the proof of the dynamic programming in Theorem \ref{theo:main_dpp}.
	Let us first provide a tower property for the {\rm BSDE} \eqref{eq:BSDE}.

	\begin{Lemma}\label{lemma:dppbsde}
		Let $(t,\om) \in [0,T] \x \Om$, $\P \in \Pc(t, \om)$,  $\tau$ be an $\F-$stopping time taking values in $[t,T]$
		and $\big(\Yc^{\P}, \Zc^{\P}, \Mc^{\P} \big)$ be a solution to the {\rm BSDE} \eqref{eq:BSDE} under $\P$.
		Then one has
		\begin{equation*}
			\Yc^{\P}_t(T, \xi) = 
			\Yc^{\P}_t \big(\tau, \Yc^{\P}_{\tau} \big)
			=
			\Yc^{\P}_t \big(\tau, \E^{\P} \big[ \Yc^{\P}_{\tau} \big|\Fc_{\tau} \big]\big),\ \P-a.s.	
		\end{equation*}
	\end{Lemma}
	\proof \rmi Given a solution $\big(\Yc^{\P}, \Zc^{\P}, \Mc^{\P} \big)$ to the {\rm BSDE} \eqref{eq:BSDE} under $\P$
	w.r.t the filtration $\F^{\P}_+ = (\Fc^{\P}_{s+})_{t \le s \le T}$,
		one then has
		\begin{equation*}
			\Yc^{\P}_t
			=
			\Yc^{\P}_{\tau} 
			-
			\int_t^{\tau} f \left(s, X_{\cdot \wedge s}, \Yc^{\P}_s, (\widehat a_s^{1/2})^\top \Zc^{\P}_s, \ah_s, b^\P_s\right) ds
			-
			 \left( \int_t^{\tau}  \Zc^{\P}_s \cdot dX_s^{c,\P} \right)^\P
			-
			\int_t^{\tau}  d \Mc^{\P}_s,
			\ \P-a.s.,
		\end{equation*}
		By taking conditional expectation w.r.t. $\Fc_{\tau}^{\P}$ under $\P$, it follows that, $\P-a.s.$,
		\begin{equation*}
			\Yc^{\P}_t
			=
			\E^{\P} \big[ \Yc^{\P}_{\tau}  \big | \Fc^{\P}_{\tau} \big]
			+
			\int_t^{\tau} f \left(s, X_{\cdot \wedge s}, \Yc^{\P}_s, (\widehat a_s^{1/2})^\top \Zc^{\P}_s, \ah_s, b^\P_s\right) ds
			-
			 \left( \int_t^{\tau}  \Zc^{\P}_s \cdot dX_s^{c,\P} \right)^\P
			-
			\int_t^{\tau}  d\widetilde \Mc^{\P}_s,
		\end{equation*}		
		where $\widetilde \Mc^{\P}_{\tau} := \E^{\P} \big[ \Mc^{\P}_{\tau}  \big | \Fc^{\P}_{\tau} \big]$, and $\widetilde \Mc^{\P}_s := \Mc^{\P}_s$ when $s < \tau$.
		By identification, we deduce that $\widetilde \Mc^{\P}_{\tau}=\Mc^{\P}_{\tau} + \E^{\P} \big[ \Yc^{\P}_{\tau}  \big | \Fc_{\tau} \big] - \Yc^{\P}_{\tau}$.
		Moreover, it is clear that $\widetilde \Mc^{\P}\in \mathbb M^p_t(\F^\P_+,\P)$ 
		and $\widetilde \Mc^{\P}$ is orthogonal to the continuous martingale $X$ under $\P$.
		
		\vspace{0.5em}
		\rmii Let us now consider the {\rm BSDE} with generator $f$ and terminal condition 
		$\E^{\P} \big[ \Yc^{\P}_{\tau}  \big | \Fc^{\P}_{\tau} \big]$, on $[t, \tau]$.
		By uniqueness of the solution to {\rm BSDE}, it follows that
		$$
			\Yc^{\P}_t(T, \xi) = 
			\Yc^{\P}_t \big(\tau, \Yc^{\P}_{\tau} \big)
			=
			\Yc^{\P}_t \big(\tau, \E^{\P} \big[ \Yc^{\P}_{\tau} \big|\Fc_{\tau} \big]\big),\ \P-a.s.	
		$$	
		\qed
	\vspace{0.5em}
	
\proof[Proof of Theorem \ref{theo:main_dpp}] 
	\rmi First, by the item \rmiii of Assumption \ref{assum:main}, it is clear that $\Ych_t(\om) = \Ych_t(\om_{t \wedge \cdot})$.
	Moreover, since $(t, \om, \P) \longmapsto \Yc^{\P}_t (\om)$ is a Borel measurable map from 
	$[0,T] \x \Om \x \M_1$ to $\R$ by Lemma \ref{lemma:measurability_BSDE},
	and the graph $[[\Pc]]$ is also a Borel measurable in $[0,T] \x \Om \x \M_1$ by Assumption \ref{assum:main},
	it follows by the measurable selection theorem that
	$(t,\om) \longmapsto \Ych_t(\om)$ is $\Bc([0,T]) \otimes \Fc_T-$universally measurable
	(or more precisely upper semi--analytic, see e.g. Proposition 7.47 of Bertsekas and Shreve \cite{bertsekas1978stochastic} or Theorem III.82 (p. 161) of Dellacherie and Meyer \cite{dellacherie1978probabilities}.


\vspace{0.5em}
	\rmii Now, using the measurable selection argument, the DPP as well as the integrability of $\Ych_{\tau}$ is a direct consequence 
	of the comparison principle and the stability of {\rm BSDE} \eqref{eq:BSDE}.
	
	First, for every $\P \in \Pc(t,\om)$ and $\eps > 0$, using the measurable selection theorem (see e.g. Proposition 7.50 of \cite{bertsekas1978stochastic} or Theorem III.82 in \cite{dellacherie1978probabilities}),
	one can choose a family of probability measures $(\Q^{\eps}_{\w})_{\w \in \Om}$ such that
	$\w \longmapsto \Q^{\eps}_{\w}$ is $\Fc_{\tau}-$measurable,
	and for $\P-a.e.$ $\w \in \Om$,
	\be \label{eq:DPP_interm_ineq}
		\Q^{\eps}_{\w} \in \Pc(\tau(\w), \w)
		~~\mbox{and}~~
		\Ych_{\tau(\w)}(\w) - \eps
		\le \E^{\Q^{\eps}_{\w}} \big[ \Yc_{\tau(\w)}^{\Q^{\eps}_{\w}}(T, \xi) \big]
		\le \Ych_{\tau(\w)}(\w).
	\ee
	We can then define the concatenated probability $\P^{\eps} := \P \otimes_{\tau} \Q^{\eps}_{\cdot}$ so that,
	by Assumption \ref{assum:main} (v), $\P^{\eps} \in \Pc(t, \om)$. Notice that $\P$ and $\P^{\eps}$ coincide on $\Fc_{\tau}$ and hence $\E^{\P^{\eps}}\big[\Yc^{\P^{\eps}}_{\tau} \big| \Fc_{\tau} \big] \in \L^p_{t,\om}(\Fc_\tau, \P)$.
	It follows then from the inequality in \eqref{eq:DPP_interm_ineq} that $\E^{\P}\big[ \big| \Ych_{\tau} \big|^p \big] < \infty$.

	\vspace{0.5em}
	\noindent Further, using the stability of the solution to {\rm BSDE} \eqref{eq:BSDE} in Lemma \ref{lemma:estimbsde} (together with Lemma \ref{lemma:deuxbsde}),
	it follows that
	$$
		\Ych_t(\om)
		\ge
		\E^{\P^{\eps}} \big[ \Yc^{\P^{\eps}}_t \big]
		=
		\E^{\P^{\eps}} \big[ \Yc^{\P^{\eps}}_t (\tau, \Yc^{\P^{\eps}}_{\tau}) \big]
		=
		\E^{\P} \big[ \Yc^{\P}_t \big(\tau, \E^{\P^{\eps}}\big[\Yc^{\P^{\eps}}_{\tau} \big| \Fc_{\tau} \big] \big) \big]
		\ge
		\E^{\P} \big[ \Yc^{\P}_t (\tau, \Ych_\tau) \big] - C \eps,
	$$
	for some constant $C > 0$ independent of $\eps$.
	By the arbitrariness of $\eps > 0$ as well as that of $\P \in \Pc(t,\om)$, one can conclude that
	$$
		\Ych_t(\om)
		\ge
		\sup_{\P \in \Pc(t,\om)} \E^{\P} \big[ \Yc^{\P}_t (\tau, \Ych_\tau) \big].
	$$
	Finally, for every $\P \in \Pc(t, \om)$, we have
	\b*
		\Yc^{\P}_t(T, \xi) = 
		\Yc^{\P}_t \big(\tau, \Yc^{\P}_{\tau} \big)
		=
		\Yc^{\P}_t \big(\tau, \E^{\P} \big[ \Yc^{\P}_{\tau} \big|\Fc_{\tau} \big]\big),
		\ \P-a.s.
	\e*
	It follows by the comparison principle of the {\rm BSDE} \eqref{eq:BSDE} (see Lemma \ref{lemma:comp} in Appendix together with Lemma \ref{lemma:deuxbsde})  that
	\b*
		\Ych_t(\om) 
		:= 
		\sup_{\P \in \Pc(t,\om)} \E^{\P} \big[ \Yc^{\P}_t (T, \xi) \big]
		=
		\sup_{\P \in \Pc(t,\om)} \E^{\P} \Big[ \Yc^{\P}_t \big(\tau, \E^{\P} \big[ \Yc^{\P}_{\tau} \big|\Fc_{\tau} \big]\big) \Big]
		\le
		\sup_{\P \in \Pc(t,\om)} \E^{\P} \big[ \Yc^{\P}_t(\tau, \Ych_{\tau}) \big].
	\e*
	We hence conclude the proof.
	\qed


\subsubsection{Further discussions}

	Notice that the essential arguments to prove the measurability of $\Ych_t(\om)$ is to construct the solution of the {\rm BSDE} in a measurable way with respect to different probabilities.
	Then the dynamic programming principle follows directly from the measurable selection theorem together with the comparison and stability of the {\rm BSDE}.
	This general approach is not limited to BSDEs with Lipschitz generators. Indeed, the solution of any BSDEs that can be approximated by a countable sequence of Lipschitz BSDEs inherits directly the measurability property. More precisely, we have the following proposition which also applies to specific super--solutions (see Section 2.3 in \cite{el1997backward} for a precise definition) of the BSDEs.
	
	\begin{Proposition} \label{others}
	Let $\mathcal Y^\P$ be the first component of the $($minimal$)$ super--solution of a {\rm BSDE} with possibly non--Lipschitz generator. We have
	\begin{enumerate}
	\item[$\mathrm{(i)}$] If there is a family $(\mathcal Y^{\P, n})$, which corresponds to the first component of a family of Lipschitz BSDEs, and a family of subsequence $(n^{\P}_k)_{k \ge 1}$ such that, $\P \longmapsto n^{\P}_k$ is $($Borel$)$ measurable,
	and $\mathcal Y^\P = \lim_{k \to \infty} \mathcal Y^{\mathbb P, n_k^{\P}}$.
	{Then $(s,\omega,\P) \longmapsto \E^{\P} \big[ \Yc^{\P}_s \big]$ is a measurable map,} and $(t, \om) \longmapsto \Ych_t(\om)$ is $\Bc([0,T]) \otimes \Fc_T-$universally measurable.

	\item[$\mathrm{(ii)}$] Furthermore, if the $($possibly non--Lipschitz$)$ {\rm BSDE} for $\cal Y^\P$ admits the comparison principle and the stability result w.r.t. its terminal condition, then
	for all $(t,\om) \in [0,T] \x \Om$ and $\F-$stopping times $\tau$ taking value in $[t,T]$, we have
		\begin{equation*}
			\Ych_t(\om) 
			=
			\sup_{\P \in \Pc(t, \om)} \E^{\P} \left[ \Yc_t^{\P} \big(\tau, \Ych_{\tau} \big) \right].
		\end{equation*}
	\end{enumerate}
	  
	\end{Proposition}

\noindent In particular, this result can be applied to BSDEs with linear growth \cite{lepeltier1997backward}, to BSDEs with general growth in $y$ \cite{pardoux1999bsdes}, to quadratic BSDEs \cite{barrieu2013monotone,kobylanski1997resultats,kobylanski2000backward}, to BSDEs with unbounded horizon \cite{darling1997backwards}, to reflected BSDEs \cite{el1997reflected}, to constrained BSDEs \cite{cvitanic1998backward,peng1999monotonic} (for the point $\mathrm{(i)}$ only),...


	\begin{Remark}
		{\rm In Assumption \ref{assum:main}, the terminal condition $\xi : \Om \longrightarrow \R$ is assumed to be Borel measurable,
		which is more restrictive comparing to the results in the context of controlled diffusion/jump process problems (where $\xi$ is only assumed to be upper semi--analytic, see {\it e.g.} \cite{neufeld2013superreplication} or \cite{karoui2013capacities}).
		This Borel measurability condition is however crucial in our {\rm BSDE} context.
		For example, when $f(t,\om, y, z, a, b) = |z|$, we know that the solution of the {\rm BSDE} \eqref{eq:BSDE} is given by $\inf_{\tilde \P \in \tilde \Pc} \E^{\tilde \P}[ \xi]$ for some family $\tilde \Pc$ of probability measure equivalent to $\P$.
		However, as is well known, the upper--semianalytic property is stable by taking supremum but not by taking infimum.}
	\end{Remark}

\section{Path regularisation of the value function}\label{sec:2}

	In this section, we will characterise a c\`adl\`ag modification of the value function $\Ych$ defined in \eqref{eq:def_Y} 
	as a semi--martingale under any $ \P\in\Pc_0$ and provide its semi--martingale decomposition. 
	In particular, this c\`adl\`ag modification will allow us to construct a solution to the second order {\rm BSDE} defined in Section \ref{sec:3}.
	
	\vspace{0.5em}
	\noindent 
	Recall that by Theorem \ref{theo:main_dpp},
	one has, under Assumption \ref{assum:main}, for any $\F-$stopping times $\tau\geq \sigma$
\begin{equation}\label{super}
			\Ych_{\sigma(\omega)}(\om) 
			=
			\sup_{\P \in \Pc(\sigma(\omega), \om)} \E^{\P} \left[ \Yc_{\sigma(\omega)}^{\P} \big(\tau, \Ych_{\tau} \big) \right].
		\end{equation}
	Moreover, we also have
\begin{equation}\label{super2}
			\Ych_{\sigma(\omega)}(\om) 
			=
			\sup_{\P \in \Pc(\sigma(\omega), \om)} \E^{\P\otimes\P_0} \left[ \widetilde \Yc_{\sigma(\omega)}^{\P\otimes\P_0} \big(\tau, \Ych_{\tau} \big) \right],
		\end{equation}
		where $\widetilde \Yc_{\sigma(\omega)}^{\P\otimes\P_0} \big(\tau, \Ych_{\tau} \big)$ is the equivalent of $\Yc_{\sigma(\omega)}^{\P} \big(\tau, \Ych_{\tau} \big)$ but defined on the enlarged space, recall \reff{eq:BSDE3} and Lemma \ref{lemma:deuxbsde}.
		
		\vspace{0.5em}

	\noindent The following technical lemma can be formally obtained by simply taking conditional expectations of the corresponding BSDEs. 
	However, this raises subtle problems about negligible sets and conditional probability measures. 
	We therefore refer the reader to \cite{claisse2016pseudo} for the precise details.

\begin{Lemma}\label{lemma:eqbsde}
	Let Assumption \ref{assum:main} hold true.
	Then for any $\P\in\Pc_0$, for any $\F-$stopping times $0\leq \sigma\leq \tau\leq T$, we have
\begin{align*}
&{\E^{\P^{\sigma(\omega)}_\omega\otimes\P_0}\Big[\widetilde \Yc^{\P^{\sigma(\omega)}_\omega\otimes\P_0}_{\sigma(\omega)}(\tau,\Ych_\tau)\Big]=\E^{\P\otimes\P_0}\Big[\left.\widetilde \Yc^{\P\otimes\P_0}_{\sigma}(\tau,\Ych_\tau)\right|\overline\Fc_\sigma\Big](\omega,\omega'),\ \text{$\P\otimes\P_0-a.e.$ $(\omega,\omega')\in\overline\Omega$},}\\
&\E^{\P^{\sigma(\omega)}_\omega}\Big[ \Yc^{\P^{\sigma(\omega)}_\omega}_{\sigma(\omega)}(\tau,\Ych_\tau)\Big]=\E^{\P}\Big[\left. \Yc^{\P}_{\sigma}(\tau,\Ych_\tau)\right|\Fc_\sigma\Big](\omega),\ \text{for $\P-a.e.$ $\omega\in\Omega$}.
\end{align*}
\end{Lemma}

\noindent Let us next remark the following immediate consequences of  \eqref{super} and \eqref{super2}
\begin{align}\label{eq:supermart}
&\Ych_{\sigma(\omega)}\geq \E^{\P} \Big[ \Yc_{\sigma(\omega)}^{\P} \big(\tau, \Ych_{\tau} \big) \Big],\ \text{for any $\P\in\Pc(\sigma(\omega),\omega)$},\\
&{{\Ych_{\sigma(\omega)}\geq \E^{\P\otimes\P_0} \left[ \widetilde\Yc_{\sigma(\omega)}^{\P\otimes\P_0} \big(\tau, \Ych_{\tau} \big) \right],\ \text{for any $\P\in\Pc(\sigma(\omega),\omega)$}.}}
\end{align}
With these inequalities, we can prove a downcrossing inequality for $\Ych$, which ensures that $\Ych$ admits right-- and left--limits outside a $\Pc_0-$polar set. Recall that 
$$\widehat{f}^\P_s(y,z):= f(s, X_{\cdot\wedge s},y,z,\widehat{a}_s, b^\P_s),\ \ \widehat{f}^{\P,0}_s:= f(s, X_{\cdot\wedge s},0,0,\widehat{a}_s, b^\P_s).$$
Let $J:=(\tau_n)_{n\in\N}$ be a countable family of $\F-$stopping times taking values in $[0,T]$ such that for any $(i,j)\in\N^2$, one has either $\tau_i\leq \tau_j $ for every $\omega\in\Omega$, or $\tau_i\geq \tau_j $ for every $\omega\in\Omega$. Let $a<b$ and $J_n\subset J$ be a finite subset $(J_n=\{0\leq \tau_1\leq \cdots \leq \tau_n\leq T\})$. We denote by $ D^b_a(\Ych,J_n)$ the number of downcrossings of the process $(\Ych_{\tau_k})_{1\leq k\leq n}$ from $b$ to $a$. We then define 
\begin{equation*}
D^b_a(\Ych,J):=\sup \{D^b_a(\Ych,J_n): J_n\subset J, \text{ and } J_n \text{ is a finite set}\}.
\end{equation*}  
The following lemma follows very closely the related result proved in Lemma A.1 of \cite{bouchard2016}. However, since $\Ych$ is not exactly an { $\Ec^{\widehat f^\P}-$super--martingale in their terminology, }we give a short proof.
\begin{Lemma}\label{downcrossing}
Fix some $\P\in\Pc_0$ and let Assumption \ref{assum:main} hold. Denote by $L$ the Lipschitz constant of the generator $f$. Then, for all $a<b$, there exists a probability measure $\overline \Q$, equivalent to {$\P\otimes\P_0$}, such that
\begin{align*}
\E^{\overline \Q}\left[ D^b_a(\Ych,J)\right]\leq \frac{e^{LT}}{b-a} \E^{\overline\Q}&\left[e^{LT}(\Ych_0\wedge b-a) -e^{-LT}(\Ych_T\wedge b-a)^++e^{LT}(\Ych_T\wedge b-a)^-\right.\\
&\left.\hspace{0.5em}+e^{LT}\int^T_0\abs{\widehat{f}^\P_s(a,0)}ds\right].
\end{align*}
Moreover, outside a $\Pc_0-$polar set, we have
$$\underset{r\in\mathbb Q\cap(t,T],r\downarrow t}{\overline\lim}\Ych_r(\om)=\underset{r\in\mathbb Q\cap(t,T],r\downarrow t}{\lim}\Ych_r(\om),\ \text{and }\underset{r\in\mathbb Q\cap(t,T],r\uparrow t}{\overline\lim}\Ych_r(\om)=\underset{r\in\mathbb Q\cap(t,T],r\uparrow t}{\lim}\Ych_r(\om).$$

\end{Lemma}

\noindent \proof Without loss of generality, we can always suppose that $0$ and $ T$ belong to $J$ and $b>a=0$. Let $J_n=\{\tau_0,\tau_1,\cdots,\tau_n\}$ with $0=\tau_0<\tau_1<\cdots<\tau_n=T$. We then consider for any $i=1,\dots,n$, and any $\omega\in\Omega$, the following {\rm BSDE} {in the enlarged space under $\overline\P^{\tau_{i-1}(\omega)}_\omega:=\P^{\tau_{i-1}(\omega)}_\omega\otimes\P_0$} on $[\tau_{i-1},\tau_i]$
\begin{align*}
\widetilde{\mathcal Y}^{i,\overline\P^{\tau_{i-1}(\omega)}_\omega}_t:=&\ \Ych_{\tau_i}-\int_t^{\tau_i}\left(\widehat f^{\P^{\tau_{i-1}(\omega)}_\omega,0}_s+\lambda^{i}_s\widetilde{\mathcal Y}^{i,\overline\P^{\tau_{i-1}(\omega)}_\omega}_s+\eta^{i}_s\cdot (\widehat{a}_s^{1/2})^\top \widetilde{\mathcal Z}^{i,\overline\P^{\tau_{i-1}(\omega)}_\omega}_s\right)ds\\
&-\int_t^{\tau_i}\widetilde{\mathcal Z}^{i,\overline\P^{\tau_{i-1}(\omega)}_\omega}_s\cdot 
{\widehat{a}_s^{1/2}dW^{\P^{\tau_{i-1}(\omega)}_\omega}_s } -\int_t^{\tau_i}d\widetilde{\mathcal M}^{i,\overline\P^{\tau_{i-1}(\omega)}_\omega}_{s},\text{ } \overline\P^{\tau_{i-1}(\omega)}_\omega-a.s.,   
\end{align*}
where $\lambda^{i}$ and $\eta^{i}$ are the two bounded processes appearing in the linearisation of $\widehat f$ (recall Assumption \ref{assum:main}(i)). We consider then the linear {\rm BSDE}, {also on the enlarged space}
{\begin{align*}\label{linearbsde}
\overline{\mathcal Y}^{i,\overline\P^{\tau_{i-1}(\omega)}_\omega}_t:=&\ \Ych_{\tau_i}-\int_t^{\tau_i}\left(\abs{\widehat f^{\P^{\tau_{i-1}(\omega)}_\omega,0}_s}+\lambda^{i}_s\overline{\mathcal Y}^{i,\overline\P^{\tau_{i-1}(\omega)}_\omega}_s+\eta^{i}_s\cdot (\widehat{a}_s^{1/2})^\top \overline{\mathcal Z}^{i,\overline\P^{\tau_{i-1}(\omega)}_\omega}_s\right)ds\\
&-\int_t^{\tau_i}\overline{\mathcal Z}^{i,\overline\P^{\tau_{i-1}(\omega)}_\omega}_s\cdot 
{\widehat{a}_s^{1/2}dW^{\P^{\tau_{i-1}(\omega)}_\omega}_s } -\int_t^{\tau_i}d\overline{\mathcal M}^{i,\overline\P^{\tau_{i-1}(\omega)}_\omega}_{s},\text{ }\overline\P^{\tau_{i-1}(\omega)}_\omega-a.s.  
\end{align*}}
It is immediate that
\begin{align*}
\overline{\mathcal Y}^{i,\overline\P^{\tau_{i-1}(\omega)}_\omega}_{\tau_{i-1}}&=\E^{\overline\P^{\tau_{i-1}(\omega)}_\omega}\left[\left.L_{\tau_i} \left(\Ych_{\tau_i}e^{\int^{\tau_i}_{\tau_{i-1}}\lambda^{i}_sds}- \int^{\tau_i}_{\tau_{i-1}} e^{\int^{s}_{\tau_{i-1}}\lambda^{i}_rdr}\abs{\widehat{f}^{\P^{\tau_{i-1}(\omega)}_\omega,0}_s}ds\right)\right|\overline\Fc^+_{\tau_{i-1}}\right],
\end{align*}
where\footnote{As usual in stochastic analysis, for any local martingale $M$, we denote by $\mathcal E(M)$ its Dol\'eans--Dade exponential.}
$$L_t:= { \Ec\left(\int_{\tau_{i-1}}^{t}\eta^{i}_s\cdot dW_s^{\overline\P^{\tau_{i-1}(\omega)}_\omega}\right),\ t\in[\tau_{i-1},\tau_i].}$$
By Assumption \ref{assum:main}(iv), for $\P-a.e.$ $\omega\in\Omega$, $\P^{\tau_{i-1}(\omega)}_{\omega}\in\Pc(\tau_{i-1}(\omega),\omega)$. Hence, by the comparison principle for BSDEs, recalled in Lemma \ref{lemma:comp} below, and {\reff{super2}}, it is clear that 
\begin{align*}
\E^{\overline \P^{\tau_{i-1}(\omega)}_\omega}\left[L_{\tau_i} \left(\Ych_{\tau_i}e^{\int^{\tau_i}_{\tau_{i-1}}\lambda^{i}_sds}- \int^{\tau_i}_{\tau_{i-1}} e^{\int^{s}_{\tau_{i-1}}\lambda^{i}_rdr}\abs{\widehat{f}^{\P^{\tau_{i-1}(\omega)}_\omega,0}_s}ds\right)\right]\leq \Ych_{\tau_{i-1}}(\omega).
\end{align*}
But, by definition of the r.p.c.d., this implies that
$$\E^{\overline \Q}\left[\left.\Ych_{\tau_i}e^{\int^{\tau_i}_{\tau_{i-1}}\lambda^{i}_sds}- \int^{\tau_i}_{\tau_{i-1}} e^{\int^{s}_{\tau_{i-1}}\lambda^{i}_rdr}\abs{\widehat{f}^{\P^{\tau_{i-1}(\cdot)}_\cdot,0}_s}ds\right|\overline\Fc_{\tau_{i-1}}\right]\leq \Ych_{\tau_{i-1}},\ \P\otimes\P_0-a.s.,$$
where the probability measure $\overline \Q$ is equivalent to {$\P\otimes\P_0$} and defined by
$$\frac{d \overline\Q}{d\P\otimes\P_0}:= { \Ec\left(\int_{\tau_{i-1}}^{t}\eta^{i}_s\cdot dW^{\P}_s\right),\ t\in[\tau_{i-1},\tau_i].}$$
Let $\lambda_s:=\sum^n_{i=1} \lambda^{i}_s{\bf 1}_{\left.\left[\tau_{i-1},\tau_i\right.\right)}(s)$, then one has that the discrete process $ (V_{\tau_i})_{0\leq i\leq n}$ defined by
$$V_{\tau_i}:=\Ych_{\tau_i}e^{\int^{\tau_i}_{0}\lambda_sds}- \int^{\tau_i}_{0} e^{\int^{s}_{0}\lambda_rdr}\abs{\widehat{f}^{\P,0}_s}ds, $$
is a $\overline \Q-$super--martingale relative to {$\overline\F$}. Then, the control on the down--crossings can be obtained exactly as in the proof of Lemma A.1 in \cite{bouchard2016}. Indeed, it is enough to observe that the original down--crossing inequality for super--martingales (see e.g. \cite[p.446]{doob1984classical}) does not require the filtration to satisfy the usual assumptions. We now prove the second part of the lemma. We define the set 
$$\Sigma:=\{\om\in\Omega\  \text{s.t.} \ \Ych_\cdot(\om) \text{ has no right$-$ or left$-$limits along the rationals at some point}\}. $$
We claim that $\Sigma$ is a $\Pc_0-$polar set. Indeed, suppose that there exists $\P\in\Pc_0$ satisfying $\P(\Sigma)>0$. Then, $\Sigma$ is non--empty and for any $\om\in\Sigma$, the path $\Ych_\cdot(\om)$ has, e.g., no right$-$limit along the rationals at some point $t\in[0,T]$. We can therefore find two rational numbers $a $, $b$ such that 
$$\underset{r\in\mathbb Q\cap(t,T],r\downarrow t}{\underline\lim}\Ych_r(\om) <a<b<\underset{r\in\mathbb Q\cap(t,T],r\downarrow t}{\overline\lim}\Ych_r(\om),$$
and the number of down--crossings $D^b_a(\Ych,J)(\om)$ of the path $\Ych_\cdot(\om) $ on the interval $[a,b]$ is equal to $+\infty$. However, the down--crossing inequality proved above shows that $D^b_a(\Ych,J)$ is $\overline\Q-a.s.$, and thus $\P-a.s.$ (see Lemma \ref{lemm:enlargement}), finite, for any pair $(a,b)$. This implies a contradiction since we assumed that $\P(\Sigma)>0$. Therefore, outside the $\Pc_0-$polar set $\Sigma$, $\Ych$ admits both right$-$ and left$-$limits along the rationals. 
\ep 

\vspace{0.5em}
\noindent We next define for all $(t,\omega)\in[0,T)\times\Omega$
\begin{equation} \label{eq:def_Ychplus}
\Ych_t^+:=\underset{r\in\mathbb Q\cap(t,T],r\downarrow t}{\overline \lim}\Ych_r,
\end{equation}
and $\Ych^+_T:=\Yc_T$. By Lemma \ref{downcrossing}, 
$$\Ych_t^+=\underset{r\in\mathbb Q\cap(t,T],r\downarrow t}{\lim}\Ych_r,\text{ outside a $\Pc_0-$polar set},$$ 
and we deduce that $\Ych^+$ is c\`ad outside a $\Pc_0-$polar set. Hence, since for any $t\in[0,T]$, $\Ych_t^+$ is by definition $\mathcal F^{\Pc_0+}_t-$measurable, we deduce that $\Ych^+$ is actually {$\F^{\Pc_0+}-$}optional. Our next result extends \reff{eq:supermart} to $\Ych^+$.

\begin{Lemma}\label{lemma:truc}
	Let Assumption \ref{assum:main} hold true.
	Then for any $0\leq s\leq t\leq T$ and $\P\in\Pc_0$, we have $\widehat{\mathcal Y}^+_t\in\mathbb L^p_0(\Fc^{\Pc_0+}_t,\P)$ and
	$$\Ych^+_s\geq \Yc^\P_s(t,\Ych^+_t),\ \P-a.s.$$
\end{Lemma}

\proof { \rmi Let $\P\in\Pc_0$, $(s_n)_{n \ge 1}$ be a sequence of rational numbers such that $s_n \downarrow t$,
	we notice that $\Ych_{s_n} \longrightarrow \Ych_t^+$, $\P-a.s.$
	Moreover, using the same conditioning argument and the inequality in \eqref{eq:DPP_interm_ineq},
	together with the estimates of Lemma \ref{lemma:estimbsde} (recall Lemma \ref{lemma:deuxbsde} and condition \eqref{eq:integrability}),
	we see that $(\Ych_{s_n})_{n \ge 1}$ are uniformly bounded in $\L^p_0(\Fc_T^{\P+}, \P)$.
	It follows that $\widehat{\mathcal Y}^+_t\in\mathbb L^p_0(\Fc^{\P+}_T,\P)$.
}
\vspace{0.5em}

	\rmii Fix some $(s,t,\omega)\in[0,T]\times[s,T]\times\Omega$ and some $\P\in\Pc_0$. Let $r^1_n\in\mathbb Q\cap(s,T],r^1_n\downarrow s $ and $r^2_n\in\mathbb Q\cap(t,T],r^2_n\downarrow t $. By \reff{eq:supermart}, we have for any $m,n\geq 1$ and $\widetilde \P\in \Pc(r^1_n,\omega)$
$$\Ych_{r^1_n}(\omega)\geq \E^{\widetilde \P} \left[\Yc_{r^1_n}^{\widetilde\P} \big(r^2_m, \Ych_{r^2_m} \big)\right].$$
In particular, thanks to Assumption \ref{assum:main}(iv), for $\P-a.e.$ $\omega\in\Omega$, we have
\begin{equation}\label{eq:eqeq}
\Ych_{r^1_n}(\omega)\geq \E^{\P^{r^1_n}_\omega} \left[\Yc_{r^1_n}^{\P^{r^1_n}_\omega} \big(r^2_m, \Ych_{r^2_m} \big)\right]=\E^\P\left[\left.\Yc^{\P}_{r^1_n}\big(r^2_m, \Ych_{r^2_m} \big)\right|\Fc_{r_n^1}\right](\omega),
\end{equation}
where we have used Lemma \ref{lemma:eqbsde}. By definition, we have 
$$\underset{n\rightarrow+\infty}{\lim}\Ych_{r^1_n}=\Ych^+_s,\ \P-a.s.$$
Next, we want to show that
\be \label{eq:cvg_cad_modif}
	\E^\P\left[\left.\Yc^{\P}_{r^1_n}\big(r^2_m, \Ych_{r^2_m} \big)\right|\Fc_{r_n^1}\right]\underset{n\rightarrow+\infty}{\longrightarrow}\Yc^{\P}_{s}\big(r^2_m, \Ych_{r^2_m} \big),\ \text{for the norm $\No{\cdot}_{\L^{1}_{s,\omega}}$}.
\ee
Indeed, we have
\begin{align*}
&\E^\P\left[\abs{\E^\P\left[\left.\Yc^{\P}_{r^1_n}\big(r^2_m, \Ych_{r^2_m} \big)\right|\Fc_{r_n^1}\right]-\Yc^{\P}_{s}\big(r^2_m, \Ych_{r^2_m} \big)}\right]\\
&=\E^\P\left[\abs{\E^\P\left[\left.\Yc^{\P}_{r^1_n}\big(r^2_m, \Ych_{r^2_m} \big)-\Yc^{\P}_{s}\big(r^2_m, \Ych_{r^2_m} \big)\right|\Fc_{r_n^1}\right]}\right]\\
&\leq \E^\P\left[\E^\P\left[\left.\abs{\Yc^{\P}_{r^1_n}\big(r^2_m, \Ych_{r^2_m} \big)-\Yc^{\P}_{s}\big(r^2_m, \Ych_{r^2_m} \big)}\right|\Fc_{r_n^1}\right]\right]\\
&= \E^\P\left[\abs{\Yc^{\P}_{r^1_n}\big(r^2_m, \Ych_{r^2_m} \big)-\Yc^{\P}_{s}\big(r^2_m, \Ych_{r^2_m} \big)}\right].
\end{align*} 
Since $\Yc^{\P}_{r^1_n}\big(r^2_m, \Ych_{r^2_m} \big)$ is c\`adl\`ag, we know that $\Yc^{\P}_{r^1_n}\big(r^2_m, \Ych_{r^2_m} \big)$ converges, $\P-a.s.$, to $\Yc^{\P}_{s}\big(r^2_m, \Ych_{r^2_m} \big)$, as $n$ goes to $+\infty$. Moreover, by the estimates of Lemma \ref{lemma:estimbsde} (together with Lemma \ref{lemma:deuxbsde}),
the quantity in the expectation above is uniformly bounded in $\L^p(\F^{\P,+},\P)$, and therefore forms a uniformly integrable family by de la Vall\'ee--Poussin criterion (since $p>1$). 
Therefore the desired convergence in \eqref{eq:cvg_cad_modif} is a simple consequence of the dominated convergence theorem. 

\vspace{0.5em}
\noindent
Hence, taking a subsequence if necessary, we have that the right--hand side of \reff{eq:eqeq} goes $\P-a.s$ to $\Yc^{\P}_{s}\big(r^2_m, \Ych_{r^2_m} \big)$ as $n$ goes to $+\infty$, so that we have
$$\Ych^+_s\geq \Yc^{\P}_{s}\big(r^2_m, \Ych_{r^2_m} \big),\ \P-a.s.$$
Next, we have by the dynamic programming for BSDEs
$$\Yc^{\P}_{s}\big(r^2_m, \Ych_{r^2_m} \big)-\Yc^{\P}_{s}\big(t, \Ych_{t}^+ \big)=\Yc^{\P}_{s}\big(r^2_m, \Ych_{r^2_m} \big)-\Yc^{\P}_{s}\big(r^2_m, \Ych_{t}^+ \big)+\Yc^{\P}_{s}\big(t, \Yc^\P_t\big(r^2_m,\Ych_t^+ \big)\big)-\Yc^{\P}_{s}\big(t, \Ych_{t}^+ \big).$$
The first difference on the right--hand side converges to $0$, $\P-a.s.$, once more thanks to the estimates of Lemma \ref{lemma:estimbsde} (together with Lemma \ref{lemma:deuxbsde}) and the definition of $\Ych^+$. As for the second difference, the same estimates show that it is controlled by $$\E^\P\left[\left.\abs{\Ych_t^+-\Yc^\P_t\big(r^2_m,\Ych_t^+ \big)}^{\tilde p}\right|\Fc_{s+}\right],$$ for some $1<\tilde p<p$. This term goes $\P-a.s.$ (at least along a subsequence) to $0$ as $m$ goes to $+\infty$ as well by Lemma \ref{lemma:bsde} (together with Lemma \ref{lemma:deuxbsde}), which ends the proof.
\ep

\vspace{0.5em}
\noindent The next lemma follows the classical proof of the optional sampling theorem for c\`adl\`ag super--martingales and extends the previous result to stopping times.
 
 \begin{Lemma}\label{lemma:optsamp}
	Let Assumption \ref{assum:main} hold true.
	Then for any $\F_+-$stopping times $0\leq \sigma\leq \tau\leq T$, for any $\P\in\Pc_0$, 
	{we have $\Ych^+_{\tau} \in \L^p_0(\Fc^{\P+}_T,\P)$ and}
	$$\Ych^+_{\sigma}\geq \Yc^\P_{\sigma}(\tau,\Ych^+_{\tau}),\ \P-a.s.$$
In particular $\Ych^+$ is c\`adl\`ag, $\Pc_0-$q.s.
\end{Lemma}

\proof { \rmi Let $\P \in \Pc_0$, $\tau$ be a $\F_+-$stopping time, we can approximate $\tau$ by a sequence of $\F-$stopping times $(\tau_n)_{n \ge 1}$ such that $\tau_n \downarrow \tau$, $\P-a.s.$
	Using again the estimates in Lemma \ref{lemma:estimbsde} as well as the same arguments in step \rmi of the proof of Lemma \ref{lemma:truc}, one can conclude that $\Ych^+_{\tau} \in \L^p_0(\Fc^{\P+}_T,\P)$.
}
\vspace{0.5em}

\rmii
Assume first that $\sigma$ takes a finite number of values $\{t_1,\dots,t_n\}$ and that $\tau$ is deterministic. Then, we have for any $\P\in\Pc_0$
$$\Ych^+_\sigma=\sum_{i=1}^n\Ych^+_{t_i}{\bf 1}_{\{\sigma=t_i\}}\geq\sum_{i=1}^n\Yc^\P_{t_i}(\tau,\Ych^+_{\tau}){\bf 1}_{\{\sigma=t_i\}}=\Yc^\P_{\sigma}(\tau,\Ych^+_{\tau}),\ \P-a.s.$$
Assume next that both $\tau$ and $\sigma$ take a finite number of values $\{t_1,\dots,t_n\}$. We have similarly
$$\Yc^\P_{\sigma}(\tau,\Ych^+_{\tau})=\sum_{i=1}^n\Yc^\P_{\sigma}(t_i,\Ych^+_{t_i}){\bf1}_{\{\tau=t_i\}}\leq \sum_{i=1}^n\Ych^+_\sigma{\bf1}_{\{\tau=t_i\}}=\Ych^+_\sigma,\ \P-a.s.$$
If $\sigma$ is general, we can always approach it from above by a decreasing sequence of $\F_+-$stopping times $(\sigma^n)_{n\geq 1}$ taking only a finite number of values. The above results imply directly that
$$\Ych^+_{\sigma^n\wedge\tau}\geq \Yc^\P_{\sigma^n\wedge\tau}(\tau,\Ych^+_{\tau}),\ \P-a.s.$$
Then, we can use the right-continuity of $\Ych^+$ and $\Yc^\P(\tau,\Ych^+_{\tau})$ to let $n$ go to $+\infty$ and obtain
$$\Ych^+_{\sigma}\geq \Yc^\P_{\sigma}(\tau,\Ych^+_{\tau}),\ \P-a.s.$$
Finally, let us take a general stopping time $\tau$. We once more approximate it by a decreasing sequence of $\F-$stopping times $(\tau^n)_{n\geq 1}$ taking only a finite number of values. We thus have
$$\Ych^+_{\sigma}\geq \Yc^\P_{\sigma}(\tau^n,\Ych^+_{\tau^n}),\ \P-a.s.$$
The term on the right-hand side converges (along a subsequence if necessary) $\P-a.s.$ to $\Yc^\P_{\sigma}(\tau,\Ych^+_{\tau})$ by Lemma \ref{lemma:bsde} (together with Lemma \ref{lemma:deuxbsde}).

\vspace{0.5em}
\noindent It remains to justify that $\Ych^+$ admits left$-$limits outside a $\Pc_0-$polar set. Fix some $\P\in\Pc_0$. Following the same arguments as in the proof of Lemma \ref{downcrossing}, we can show that for some probability measure $\overline \Q$ equivalent to {$\P\otimes\P_0$} and some bounded process $\lambda$, 
$$V_t:=\Ych_te^{\int_0^t\lambda_sds}-\int_0^te^{\int_0^s\lambda_udu}\abs{\widehat f^{\P,0}_s}ds,$$
is a right--continuous $(\overline \Q,{\overline\F_+})-$super--martingale, which is in addition uniformly integrable under $\overline \Q$ since $\Ych$ and $\widehat f^{\P,0}$ are uniformly bounded in {$\L^p(\overline\Fc_T,\P\otimes\P_0)$} and thus in {$\L^{\tilde p}(\overline\Fc_T,\overline \Q)$} for some $1<\tilde p<p$. Therefore, for any increasing sequence of {$\overline\F_+-$}stopping times $(\rho_n)_{n\geq 0}$ taking values in $[0,T]$, the sequence $(\E^{\overline \Q}[V_{\rho_n}])_{n\geq 0}$ is non--increasing and admits a limit. By Theorem VI-48 and Remark VI-50(f) of \cite{dellacherie1978probabilities}, we deduce that $V$, and thus $\Ych^+$, admit left$-$limits outside a $\overline \Q-$negligible (and thus $\P-$negligible by Lemma \ref{lemm:enlargement}) set. Moreover, the above implies that the set
$$\left\{\omega\in\Omega: \Ych^+(\omega)\text{ admits left$-$limits}\right\},$$
is the complement of a $\Pc_0-$polar set, which ends the proof.
\ep

\vspace{0.5em}
\noindent Our next result shows that $\Ych^+$ satisfies a representation formula. 
This representation is crucial to prove the existence of a solution to a second order BSDE in Section 4.4. We first define the following subset of probability measures. Given $(t,\omega)\in[0,T]\times\Omega$ and a filtration $\mathbb X = (\mathcal X_s)_{t \le s \le T}$, define 
\be \label{eq:PcPX}
		\mathcal P_{t,\omega}(s,\mathbb P,\X):=\left\{\mathbb P^{'}\in\mathcal P(t,\omega),\ \mathbb P^{'}=\mathbb P \text{ on }\mathcal X_s\right\}.
	\ee
When $t=0$, we simply denote $\mathcal P_{0,\omega}(s,\mathbb P,\X)$ by $\mathcal P_0(s,\mathbb P,\X)$.
	
\begin{Lemma}\label{lemma:rep}
	Let Assumption \ref{assum:main} hold true.
	Then for any $\F-$stopping times $0\leq \sigma\leq \tau\leq T$, for any $0\leq t\leq T$, for any $\P\in\Pc_0$, we have
$$ \Ych_\sigma=\underset{\P^{'}\in\Pc_0(\sigma,\P,\F)}{\esup^\P}\E^{\P^{'}}\left[\left.\Yc^{\P^{'}}_\sigma(\tau,\widehat\Yc_\tau)\right|\Fc_{\sigma}\right],\ \P-a.s., \ \text{and }\Ych^+_t=\underset{\P^{'}\in\Pc_0(t,\P,\F_{+})}{\esup^\P}\Yc^{\P^{'}}_{t}(T,\xi),\ \P-a.s.$$
 \end{Lemma}
 
  \proof
 We start with the first equality. By definition and Lemma \ref{lemma:eqbsde}, for any $\P^{'}\in\Pc_0(\sigma,\P,\F)$ we have
 $$\Ych_\sigma\geq \E^{\P^{'}}\left[\left.\Yc^{\P^{'}}_\sigma(\tau,\Ych_\tau)\right|\Fc_{\sigma}\right],\P^{'}-a.s.$$
 But since both sides of the inequality are $\Fc_\sigma^{U}-$measurable and $\P^{'}$ coincides with $\P$ on $\Fc_\sigma$ (and thus on $\Fc_\sigma^{U}$, by uniqueness of universal completion) the above also holds $\P-a.s.$ We deduce
 $$\Ych_\sigma\geq \underset{\P^{'}\in\Pc_0(\sigma,\P,\F)}{\esup^\P}\E^{\P^{'}}\left[\left.\Yc^{\P^{'}}_\sigma(\tau,\Ych_\tau)\right|\Fc_{\sigma}\right],\P-a.s.$$
Next, notice that by Lemmas \ref{lemma:measurability_BSDE} and \ref{lemma:dppbsde}, $(t, \om, \Q) \longmapsto \E^{\Q} \big[ \Yc_t^{\Q} (T, \xi) \big] $ is Borel measurable.
Recall that $\Ych_t(\om) =  \sup_{\P \in \Pc(t,\om)} \E^{\P}\big[ \Yc^{\P}_t(T, \xi) \big]$.
As in the proof of Theorem \ref{theo:main_dpp}, it follows by the measurable selection theorem 
(see e.g. Proposition 7.47 of \cite{bertsekas1978stochastic}) that for every $\eps > 0$, there is a family of probability measures $(\Q^{\eps}_{\w})_{\w \in \Om}$ such that $\w \longmapsto \Q^{\eps}_{\w}$ is $\Fc_{\sigma}$ measurable and for $\P-a.e.$ $\w \in \Om$,
$$
	\Ych_{\sigma(\w)} (\w) 
	\le
	\E^{\Q^{\eps}_{\w}} \big[ \Yc_{\sigma(\w)}^{\Q^{\eps}_{\w}} (T, \xi) \big] + \eps
	=
	\E^{\Q^{\eps}_{\w}} \big[ \Yc_{\sigma(\w)}^{\Q^{\eps}_{\w}} (\tau, \Yc^{\Q^{\eps}}_\tau) \big] + \eps
	\le
	\E^{\Q^{\eps}_{\w}} \big[ \Yc_{\sigma(\w)}^{\Q^{\eps}_{\w}} (\tau, \Ych_\tau) \big] + \eps,
$$
Let us now define the concatenated probability $\P^{\eps} := \P \otimes_{\sigma} \Q^{\eps}_{\cdot}$ so that
$\P^{\eps} \in \Pc_0(\sigma, \P, \F)$, it follows then by Lemma \ref{lemma:eqbsde} that
$$
	\Ych_{\sigma} 
	\le 
	 \E^{\P^{\eps}}\left[\left.\Yc^{\P^{\eps}}_\sigma(\tau,\Ych_\tau)\right|\Fc_{\sigma}\right] + \eps
	\le 
	\underset{\P^{'}\in\Pc_0(\sigma,\P,\F)}{\esup^\P} \E^{\P^{'}}\left[\left.\Yc^{\P^{'}}_\sigma(\tau,\Ych_\tau)\right|\Fc_{\sigma}\right] + \eps,~\P-a.s.
$$
We hence finish the proof of the first equality by arbitrariness of $\eps > 0$.

\vspace{0.5em}

\noindent Let us now prove the second equality. Let $r^1_n\in\mathbb Q\cap(t,T],r^1_n\downarrow t $. By the first part of the proof, we have
 $$ \Ych_{r^1_n}=\underset{\P^{'}\in\Pc_0(r_n^1,\P,\F)}{\esup^\P}\E^{\P^{'}}\left[\left.\Yc^{\P^{'}}_{r_n^1}(T,\xi)\right|\Fc_{r_n^1}\right],\ \P-a.s.$$
 Since for every $n\in\N$, $\Pc_0(r_n^1,\P,\F)\subset \Pc_0(t,\P,\F_{+})$, we deduce as above that for any $\P^{'}\in\Pc_0(t,\P,\F_{+})$ and for $n$ large enough
 $$ \Ych_{r^1_n}\geq \E^{\P^{'}}\left[\left.\Yc^{\P^{'}}_{r_n^1}(T,\xi)\right|\Fc_{r_n^1}\right],\ \P-a.s.$$
Arguing exactly as in the proof of Lemma \ref{lemma:truc}, we can let $n$ go to $+\infty$ to obtain
  $$ \Ych^+_{t}\geq \Yc^{\P^{'}}_{t}(T,\xi),\ \P-a.s.,$$
  which implies by arbitrariness of $\P^{'}$
   $$ \Ych^+_{t}\geq \underset{\P^{'}\in\Pc_0(t,\P,\F_{+})}{\esup^\P}\Yc^{\P^{'}}_{t}(T,\xi),\ \P-a.s.$$
We claim next that for any $n\in\N$, the following family is upward directed
$$\left\{\E^{\P^{'}}\left[\left.\Yc^{\P^{'}}_{r^1_n}(T,\xi)\right|\Fc_{r_n^1}\right],\ \P^{'}\in\Pc_0(r_n^1,\P,\F)\right\}.$$
Let us consider $(\P^1,\P^2)\in  \mathcal{P}_0(r_n^1,\mathbb{P},\F)\times \mathcal{P}_0(r_n^1,\mathbb{P},\F)$, and let us define the following subsets of $\Omega$
$$A_1:=\left\{\omega\in\Omega: \E^{\P^{1}}\left[\left.\Yc^{\P^{1}}_{r^1_n}(T,\xi)\right|\Fc_{r_n^1}\right](\omega)>\E^{\P^{2}}\left[\left.\Yc^{\P^{2}}_{r^1_n}(T,\xi)\right|\Fc_{r_n^1}\right](\omega)\right\},$$
 \vspace{-0.5em}
 $$A_2:=\Omega\backslash A_1.$$
 Then, $A_1,A_2\in\mathcal F_{r^n_1}^{\P}$, and we can define the following probability measure on $(\Omega,\Fc_T)$
 $$\P^{1,2}(B):=\P^1(A_1\cap B)+\P^2(A_2\cap B),\ \text{for any $B\in\Fc_T$}.$$
 By Assumption \ref{assum:main}(v), we know that $\P^{1,2}\in\mathcal P_0$, and by definition, we further have $ \P^{1,2}\in\mathcal{P}_0(r^n_1,\mathbb{P},\F)$ as well as, $\P-a.s.$,
 $$\mathbb E^{\mathbb{P}^{1,2}}\left[\left.\Yc^{\P^{1,2}}_{r^1_n}(T,\xi)\right|\Fc_{r_n^1}\right]=\mathbb E^{\mathbb{P}^{1}}\left[\left.\Yc^{\P^{1}}_{r^1_n}(T,\xi)\right|\Fc_{r_n^1}\right]\vee\mathbb E^{\mathbb{P}^{2}}\left[\left.\Yc^{\P^{2}}_{r^1_n}(T,\xi)\right|\Fc_{r_n^1}\right],$$
 which proves the claim.

According to \cite{neveu1975discrete}, we then know that there exists some sequence $(\P^m_n)_{m\geq 0}\subset\Pc_0(r_n^1,\P,\F)$ such that
$$\Ych_{r_n^1}=\underset{m\rightarrow+\infty}{\lim}\uparrow\E^{\P^m_n}\left[\left.\Yc^{\P^{m}_n}_{r_n^1}(T,\xi)\right|\Fc_{r_n^1}\right],\ \P-a.s.$$
By dominated convergence (recall that the $\Yc^\P$ are in $\D^{p}_0(\F^{\P+},\P)$, with a norm independent of $\P$, by Lemma \ref{lemma:estimbsde}), the above convergence also holds for the $\L^{\tilde p}_0(\P)-$norm, for any $1<\tilde p<p$. By the stability result of Lemma \ref{lemma:estimbsde} (together with Lemma \ref{lemma:deuxbsde}) and the monotone convergence theorem, we deduce that
\begin{align*}
\Yc_t^{\P}(r_n^1,\Ych_{r_n^1})&=\Yc_t^\P\left(r_n^1,\underset{m\rightarrow+\infty}{\lim}\uparrow\E^{\P^m_n}\left[\left.\Yc^{\P^{m}_n}_{r_n^1}(T,\xi)\right|\Fc_{r_n^1}\right]\right),\ \P-a.s.\\
&=\underset{m\rightarrow+\infty}{\lim}\ \Yc_t^\P\left(r_n^1,\E^{\P^m_n}\left[\left.\Yc^{\P^{m}_n}_{r_n^1}(T,\xi)\right|\Fc_{r_n^1}\right]\right),\ \P-a.s.\\
&=\underset{m\rightarrow+\infty}{\lim}\ \Yc_t^{\P^m_n}\left(r_n^1,\E^{\P^m_n}\left[\left.\Yc^{\P^{m}_n}_{r_n^1}(T,\xi)\right|\Fc_{r_n^1}\right]\right),\ \P-a.s.\\
&=\underset{m\rightarrow+\infty}{\lim}\ \Yc_t^{\P^m_n}\left(r_n^1,\Yc^{\P^{m}_n}_{r_n^1}(T,\xi)\right),\ \P-a.s.\\
&=\underset{m\rightarrow+\infty}{\lim}\ \Yc_t^{\P^m_n}\left(T,\xi\right),\ \P-a.s.\\
&\leq \underset{\P^{'}\in\Pc_0(t,\P,\F_+)}{\esup^\P}\Yc_t^{\P^{'}}\left(T,\xi\right),\ \P-a.s.,
\end{align*}
where we have used in the third equality the fact that $\P^m_n$ coincides with $\P$ on $\Fc_{r_n^1}$ and that $\Yc_t^\P$ is $\Fc_t^+-$measurable, Lemma \ref{lemma:dppbsde} in the fourth equality, and the dynamic programming principle for BSDEs in the fifth equality.

\vspace{0.5em}
\noindent It remains to let $n$ go to $+\infty$ and to use Lemma \ref{lemma:bsde} (together with Lemma \ref{lemma:deuxbsde}) to obtain the desired equality. 
 \ep
 
 \begin{Remark}
 {\rm Notice that we can prove with the exact same techniques the following DPP--type representation for any $\F_+-$stopping times $0\leq \tau\leq\sigma\leq T$
 $$\Ych^+_\tau=\underset{\P^{'}\in\Pc_0(\tau,\P,\F_{+})}{\esup^\P}\Yc^{\P^{'}}_{\tau}(\sigma,\widehat \Yc^+_\sigma),\ \P-a.s.$$
 However, since we did not need this result, we refrained from writing explicitly its proof.
 }
 \end{Remark}

\vspace{0.5em}
\noindent The next result shows that $\Ych^+$ is actually a semi--martingale under any $\P\in\Pc_0$, and gives its decomposition.
 \begin{Lemma}\label{lemma:semi}
 Let Assumption \ref{assum:main} hold. For any $\P\in\Pc_0$, there is $(Z^\P, M^\P, K^\P)\in \mathbb H^{p}_0(\F^{\P+},\P)\times \mathbb M^p_0(\F^{\P+},\P)\times \mathbb I^p_0(\F^{\P+}, \P)$ such that 
 $$\Ych_t^+=\xi-\int_t^T\widehat f^\P_s(\Ych^+_s, (\widehat a_s^{1/2})^\top Z^\P_s)ds-\int_t^TZ^\P_s\cdot dX^{c,\P}_s-\int_t^TdM^\P_s+\int_t^TdK_s^\P,\ t\in[0,T],\ \P-a.s.$$
Moreover, there is some $\F^{\Pc_0}-$predictable process $Z$ which aggregates the family $(Z^\P)_{\P\in\Pc_0}$.
 \end{Lemma}
 \proof
Fix some $\P\in\Pc_0$. Consider the following reflected {\rm BSDE} { on the enlarged space}. For $0\leq t\leq T$, {$\P\otimes\P_0$}$-a.s.$ 
\begin{align*}
\begin{cases}
\displaystyle \bar y_t^\P=\xi-\int_t^T\widehat f^\P_s(\bar y_s^\P, (\widehat a_s^{1/2})^\top \bar z^\P_s)ds-\int_t^T \bar z_s^\P \cdot \widehat{a}_s^{1/2}dW_s^\P-\displaystyle \int_t^Td \bar m^\P_s+ \bar  k^\P_T- \bar  k^\P_t,\\
\displaystyle \bar  y_t^\P\geq \Ych^+_t,\\ 
\displaystyle\int_0^T\left(\bar y_{t-}^\P-\Ych^+_{t-}\right)d\bar  k^\P_t=0.
\end{cases}
\end{align*}

\noindent By Theorem 3.1 in \cite{bouchard2015unified}, this reflected {\rm BSDE} is well--posed and $\bar y^\P$ is c\`adl\`ag. 
By abuse of notation, we denote $\Ych^+(\omb) := \Ych^+(\pi(\omb))$.
We claim that $\bar y^\P=\Ych^+$, {$\P\otimes\P_0-a.s.$} Indeed, we argue by contradiction, and assume without loss of generality that $\bar y^\P_0>\Ych^+_0$. For each $\varepsilon>0$, denote $\tau_\varepsilon:=\inf\;\{t: \bar y^\P_t\leq \Ych^+_t+\epsilon\} $. Then $\tau_\varepsilon $ is an {$\overline\F_{+}- $}stopping time and $\bar y^\P_{t-}\geq \Ych^+_{t-}+\epsilon>\Ych^+_{t-} $ for all $t\leq \tau_\varepsilon $. Thus {$ \bar k^\P_t=0,\ \P\otimes\P_0-a.s.,$} for $0\leq t\leq \tau_\varepsilon $ and thus
$$\bar y^\P_t=\bar y^\P_{\tau_\varepsilon}-\int_t^{\tau_\varepsilon}\widehat f^\P_s(\bar y^\P_s, (\widehat a_s^{1/2})^\top \bar z^\P_s)ds-\int_t^{\tau_\varepsilon} \bar z^\P_s\cdot \widehat{a}_s^{1/2}dW^\P_s-\int_t^{\tau_\varepsilon}d \bar m^\P_s,\ \mathbb P\otimes\P_0-a.s.$$
The same linearization argument that we used in the proof of Lemma \ref{lemma:estimbsde} implies that 
$$\bar y^\P_0\leq  \Yc^{\P\otimes\P_0}_0({\tau_\varepsilon},\Ych^+_{\tau_\varepsilon})+C\E^{\P\otimes\P_0}\left[\bar y^\P_{\tau_\varepsilon}-\Ych^+_{\tau_\varepsilon}\right]\leq  \Yc^{\P\otimes \P_0}_0({\tau_\varepsilon},\Ych^+_{\tau_\varepsilon})+C\eps,$$
 for some $C>0$. However, by Lemma \ref{lemma:optsamp}, we know that $\Yc^{\P\otimes\P_0}_0({\tau_\varepsilon},\Ych^+_{\tau_\varepsilon})\leq \Ych^+_0$, which contradicts the fact that $\bar y^\P_0>\Ych^+_0 $.  

\vspace{0.5em}
\noindent Then, by exactly the same arguments as in Lemma \ref{lemma:deuxbsde}, we can go from the enlarged space to $\Omega$ and obtain for some $(Z^\P)_{\P\in\mathcal P_0}\subset \mathbb H^p_0(\F^\P_+,\P)$, and $(M^\P,K^\P)_{\P\in\mathcal P_0} \subset \mathbb M^p_0(\F^{\P}_+,\P)\times \mathbb I^p_0(\F^{\P}_+,\P)$
$$\Ych_t^+=\xi-\int_t^T\widehat f^\P_s(\Ych^+_s, (\widehat a_s^{1/2})^\top Z^\P_s)ds-\int_t^TZ^\P_s\cdot dX^{c,\P}_s-\int_t^TdM^\P_s+\int_t^TdK_s^\P,\ t\in[0,T],\ \P-a.s.$$
Then, by Karandikar \cite{karandikar1995pathwise}, since $\Ych^+$ is a c\`adl\`ag semi--martingale, we can define a universal process denoted by $\langle \Ych^+, X \rangle$ which coincides with the quadratic co--variation of $\Ych^+$ and $X$ under each probability $\P \in \Pc_0$.
In particular, the process $\langle \Ych^+, X \rangle$ is $\Pc_0-$quasi--surely continuous and hence is $\F^{\Pc_0+}-$predictable (or equivalently  $\F^{\Pc_0}-$predictable).
Similarly to the proof of Theorem 2.4 of \cite{nutz2015robust}, we can then define a universal  $\F^{\Pc_0}-$predictable process $Z$ by
$$
	Z_t := \widehat a^{\oplus}_t \frac{d \langle \Ych^+, X \rangle_t}{dt},
$$
where we recall that $ \widehat a^{\oplus}_t$ represents the Moore--Penrose pseudoinverse of $ \widehat a_t$.
In particular, $Z$ aggregates the family $\{ Z^\mathbb P,\P\in\Pc_0\}$. 
\ep
 
 \vspace{0.5em}
\noindent We end this section with a remark, which explains that in some cases, the path regularisation that we used could be unnecessary in order to obtain a super-martingale decomposition as in Lemma \ref{lemma:semi}.

\begin{Remark}
{\rm Assume that for any $(t,\omega)\in[0,T]\times\Omega$, all the probability measures in $\Pc(t,\omega)$ satisfy the Blumenthal $0-1$ law. This would be the case for instance if we where working with the set $\overline{\Pc}_S$ defined and used in \cite{soner2012wellposedness}. Then, for any $\P\in\Pc(t,\omega)$, the filtration $\F^\P$ is right--continuous and therefore satisfies the usual conditions.
	Then we know that for any $\P$, there is a $\P-$version of $\widehat \Yc$ which is $\F^{\P}-$optional, and thus also l\`adl\`ag. Furthermore, by Lemma \ref{lemma:rep}, it verifies for any $p'<p$
$$\underset{\P\in\Pc_0}{\sup}\E^\P\left[\underset{t\in[0,T]}{{\rm essup}^\P}\abs{\widehat \Yc_t}^{p'}\right]<+\infty.$$
Moreover, by the Blumenthal $0-1$ law, \eqref{super} and \eqref{super2} rewrite
\begin{align*}
\Ych_{\sigma(\omega)}(\om) 
			&=
			\sup_{\P \in \Pc(\sigma(\omega), \om)}\Yc_{\sigma(\omega)}^{\P} \big(\tau, \Ych_{\tau} \big)\=
			\sup_{\P \in \Pc(\sigma(\omega), \om)}  \widetilde \Yc_{\sigma(\omega)}^{\P\otimes\P_0} \big(\tau, \Ych_{\tau} \big).
\end{align*}
Hence, $\widehat \Yc$ is a $\Ec^{\widehat f^\P}-$super--martingale in the terminology of \cite{bouchard2016}. We can then apply Theorem 3.1 of \cite{bouchard2016} to obtain directly the semi--martingale decomposition of Lemma \ref{lemma:semi}. The aggregation of the family $(Z^\P)_{\P\in\Pc_0}$ can still be done, but requires to use Karandikar's approach \cite{karandikar1995pathwise}, combined with the It\^o formula for l\`adl\`ag processes of \cite[p. 538]{lenglart1980tribus}. Then, one can also generalise the results on 2BSDEs of the section below. This however requires that in the definition of a {\rm 2BSDE} $($see Definition \ref{def:1}$)$, the processes $Y$ and $K$ are only l\`adl\`ag, instead of c\`adl\`ag, and $Y$ is only required to have $\P-$version which is $\F^{\P}-$optional. With this change, all our results still go through. }
\end{Remark}

\section{Application to the second order BSDEs}\label{sec:3}

	We now apply the previous dynamic programming result to establish a wellposedness result for the second order BSDEs (2BSDEs).
	Let us first formulate  a strengthened integrability condition on $\xi$ and the generator $f$.
	For $p > \kappa > 1$, $(t,\om) \in [0,T] \x \Om$ and $\zeta \in \L^{1}_{t,\omega}(\Fc_T)$, we denote
	$$
		\No{\zeta}_{\mathbb L_{t,\omega}^{p,\kappa}}^p
		:=
		\underset{\mathbb P\in\mathcal P(t,\omega)}{\sup} 
		\mathbb E^{\mathbb P}\left[
			\underset{t\leq s\leq T}{\esup}^{\mathbb P}\left( \underset{\mathbb P^{'}\in \mathcal P_{t,\omega}(s, \mathbb P,\F)}{\esup^{\mathbb P}}\mathbb E^{\mathbb P^{'}} \big[|\zeta|^{\kappa} \big|\mathcal F^+_s \big]
			\right)^{\frac{p}{\kappa}}
		\right],
	$$
	and
	$$\mathbb L_{t,\omega}^{p,\kappa}(\Fc_T):=\left\{\zeta \in \L_{t,\omega}^{p}(\Fc_T),
		\ \No{\zeta}^p_{\mathbb L_{t,\omega}^{p,\kappa}}<+\infty\right\},$$
	where, as defined in \eqref{eq:PcPX}, 
	$ 
		\mathcal P_{t,\omega}(s,\mathbb P,\F):=\left\{\mathbb P^{'}\in\mathcal P(t,\omega),\ \mathbb P^{'}=\mathbb P \text{ on }\mathcal F_s\right\}.
	$
	We will make use of the following assumption.
\begin{Assumption}\label{assum:integ}
	For the given $p > \kappa > 1$, one has $\xi \in \mathbb L_0^{p,\kappa}(\Fc_T)$ and
	$$
		\phi_f^{p,\kappa}
		:=
		\underset{\mathbb P\in\mathcal P_{0}}{\sup}\mathbb E^{\mathbb P}\left[\underset{0\leq s\leq T}{\esup}^{\mathbb P}\left(\underset{\mathbb P^{'}\in \mathcal P_0(s, \mathbb P,\F_+)}{\esup^{\mathbb P}}\mathbb E^{\mathbb P^{'}}\left[\int^T_0\big|\widehat f^{\P^{'},0}_t\big|^\kappa dt\Bigg|\Fc_s^+\right]\right)^{\frac{p}{\kappa}}\right]<+\infty.
	$$
\end{Assumption}

	\begin{Remark}{\rm
		Since the early works \cite{soner2011martingale,soner2012wellposedness}, the space $\mathbb L^{p,\kappa}_{t,\omega}(\Fc_T)$ has been recognised to be particularly suited for 2BSDE theory. This is mainly due to the fact that it is still an open problem to prove whether Doob's inequality, in the form of \cite[Proposition A.1]{possamai2013second1} can be improved or not.
		As in its current state, one knows that the norm $\No{\cdot}_{\mathbb L^{p,\kappa}_{t,\omega}}$ is dominated by the norm $\No{\cdot}_{\mathbb L^{p'}_{t,\omega}}$ for any $p'>p$. In other words, we have $\L^{p,\kappa}_{t,\omega}(\Fc_T)\subset \L^{p'}_{t,\omega}(\Fc_T)$ for any $p'>p> \kappa> 1$. As such, Assumption \ref{assum:integ} is tailor--made to obtain the {\it a priori} estimates of Theorem \ref{estimatesref} below.}
	\end{Remark}

\subsection{Definition}

We shall consider the following $2${\rm BSDE}, which verifies, $\mathcal P_0-q.s.$ 
\begin{equation}
Y_t=\xi -\! \int_t^T\widehat{f}^\P_s(Y_s, (\widehat a_s^{1/2})^\top Z_s)ds -\! \left(\int_t^T Z_s\cdot dX^{c,\P}_s\right)^\P \!-\!\int_t^T dM^{\mathbb P}_s + K^{\mathbb P}_T-K^{\mathbb P}_t,\ 0\leq t\leq T.
\label{2bsde}
\end{equation}

\begin{Definition}\label{def:1}
	We say that  
	$(Y,Z,(M^\P)_{\P\in\mathcal P_0},(K^\P)_{\P\in\mathcal P_0})
	\in 
	\mathbb D^{p}_0(\F^{\Pc_0+}) \times \mathbb H^{p}_0(\F^{\Pc_0}) \times (\mathbb M^p_0(\F^{\P+}))_{\P\in\Pc_0}$
	$\times (\mathbb I^p_0(\F^{\P}_+))_{\P\in\Pc_0}$
	is a solution to the ${\rm 2BSDE}$ \reff{2bsde} if \reff{2bsde} holds $\Pc_0-q.s.$ and if the family $(K^{\mathbb P})_{\mathbb P \in \mathcal P_{0}}$ satisfies the minimality condition
\begin{equation}
K_t^{\mathbb P}=\underset{ \mathbb{P}^{'} \in \mathcal{P}_0(t,\mathbb{P},\F_{+}) }{ \einf^{\mathbb P} }\mathbb{E}^{\mathbb P^{'}}\left[\left. K_T^{\mathbb{P}^{'}}\right|\Fc_{t}^{\P^{'}+}\right], \text{ } 0\leq t\leq T, \text{  } \mathbb P-a.s., \text{ } \forall \mathbb P \in \mathcal P_{0},
\label{2bsde.minK}
\end{equation}
where $\mathcal{P}_0(t,\mathbb{P},\F_{+})$ is defined in \eqref{eq:PcPX}.
\end{Definition}

\begin{Remark}
{\rm If we assume that $b^\P=0$, $\P-a.s.$ for any $\P\in\Pc_0$, then we have that $X^{c,\P}=X,$ $\P-a.s.$ for any $\P\in\Pc_0$. Then, we can use the general result given by Nutz \cite{nutz2012pathwise}\footnote{Notice that this result only holds under some particular set-theoretic axioms. For instance, one can assume the usual Zermelo--Fraenkel set theory, plus the axiom of choice (ZFC for short), and either add the continuum hypothesis or Martin's axiom (which is compatible with the negation of the continuum hypothesis). Actually, the required axioms must imply the existence of the so--called medial limits in the sense of Mokobodzki. As far as we know, the weakest set of axioms known to be sufficient for the existence of medial limits (see \cite[538S]{fremlin2008measure} and \cite{larson2009filter}) is ZFC plus the statement that the reals are not a union of fewer than continuum many meager sets. Moreover, ZFC alone is not sufficient, in the sense that by Corollary 3.3 of \cite{larson2009filter}, If ZFC is consistent, then so is ZFC + "there exist no medial limits".} to obtain the existence of a $\Pc_0-q.s.$ c\`adl\`ag $\F^{\Pc_0+}-$progressively measurable process, which we denote by $\int_0^\cdot Z_s\cdot dX_s$, such that
$$\int_0^\cdot Z_s\cdot dX_s=\left(\int_0^\cdot Z_s\cdot dX_s\right)^\P,\ \P-a.s.$$
Hence, we can then also find an $\F^{\Pc_0+}-$progressively measurable process $N$ which aggregates the process $M^{\mathbb P}-K^{\mathbb P} $, and which is therefore a $(\F^{\P}_+,\P)-$super--martingale for any $\P\in\Pc_0$. However, the Doob--Meyer decomposition of $N$ into a sum of a martingale and a nondecreasing process generally depends on $\mathbb P$. If furthermore the set $\Pc_0$ only contains elements satisfying the predictable martingale representation property, for instance the set $\overline{\Pc}_S$ used in \cite{soner2012wellposedness}, then we have that $M^\P=0,\ \P-a.s.$, for any $\P\in\Pc_0$, so that the above reasoning allows to aggregate the non-decreasing processes $K^\P$.}
\end{Remark}


\noindent We first state the main result of this part

\begin{Theorem}\label{th:mainex}
	Let Assumptions \ref{assum:main} and \ref{assum:integ} hold true.
	Then there exists a unique solution $(Y,Z,(M^\P)_{\P\in\mathcal P_0},(K^\P)_{\P\in\mathcal P_0}) $ to the {\rm2BSDE} \reff{2bsde}.
\end{Theorem}

\subsection{Uniqueness, stochastic control representation and comparison}
We start by proving a representation of a solution to 2BSDEs, which provides incidentally its uniqueness.
\begin{Theorem}\label{representation}
Let Assumptions \ref{assum:main} and \ref{assum:integ} hold true, 
and $(Y,Z,(M^\P)_{\P\in\mathcal P_0},(K^\P)_{\P\in\mathcal P_0})$ be a solution to the {\rm 2BSDE} \reff{2bsde}. For any $\P\in\mathcal P_0$, let $(\mathcal Y^\mathbb P,\mathcal Z^\mathbb P,\mathcal M^\mathbb P)\in \mathbb D_0^{p}(\F^{\P}_+,\P) \times \mathbb H_0^{p}(\F^{\P}_+,\P)\times \mathbb M_0^p(\F^{\P}_+,\P)$ be the solutions of the corresponding BSDEs \reff{eq:bsde}. Then, for any $\mathbb P \in \mathcal{P}_{0}$ and $0\leq t_1\leq t_2 \leq T$,

\begin{equation}\label{eq:2bsderep}
Y_{t_1}=\underset{ \mathbb{P}^{'} \in \mathcal{P}_0(t_1,\mathbb{P},\F_{+}) }{ \esup^{\mathbb P} } \mathcal Y_{t_1}^{\mathbb{P}^{'}}(t_2,Y_{t_2}), \text{  } \mathbb P-a.s.
\end{equation}
Thus, the {\rm 2BSDE} \reff{2bsde} has at most one solution in $\mathbb D_0^{p}(\F^{\Pc_0+}) \times \mathbb H^{p}_0(\F^{\Pc_0})\times \mathbb M^p_0((\F^{\P}_+)_{\P\in\Pc_0})\times \mathbb I^p_0((\F^{\P}_+)_{\P\in\Pc_0}).$
\end{Theorem}

\proof
We start by proving the representation \reff{eq:2bsderep} in three steps.

\vspace{0.5em}
$\mathrm{(i)}$ Fix some $\P\in\mathcal P_0$ and then some $\P^{'}\in\mathcal{P}_0(t_1,\mathbb{P},\F^{\P}_+)$. Since \reff{2bsde} holds $\P^{'}-a.s.$, we can see $Y$ as a supersolution of the {\rm BSDE} on $[t_1,t_2]$, under $\P^{'}$, with generator $\widehat f^{\P^{'}}$ and terminal condition $Y_{t_2}$. By the comparison principle of Lemma \ref{lemma:comp} (together with Lemma \ref{lemma:deuxbsde}), we deduce immediately that $Y_{t_1}\geq \mathcal Y^{\P^{'}}_{t_1}(t_2,Y_{t_2}),$ $\P^{'}-a.s.$ Then, since $\mathcal Y^{\P^{'}}_{t_1}(t_2,Y_{t_2})$ (or a $\P^{'}-$version of it) is $\mathcal F_{t_1}^+-$measurable and since $Y_{t_1} $ is $\Fc^{\Pc_0+}_{t_1}-$measurable, we deduce that the inequality also holds $\P-a.s.$, by definition of $\mathcal{P}_0(t_1,\mathbb{P},\F_{+})$ and the fact that measures extend uniquely to the completed $\sigma-$algebras. We deduce that
$$Y_{t_1}\geq\underset{ \mathbb{P}^{'} \in \mathcal{P}_0(t_1,\mathbb{P},\F_{+}) }{ \esup^{\mathbb P} } \mathcal Y_{t_1}^{\mathbb{P}^{'}}(t_2,Y_{t_2}), \text{  } \mathbb P-a.s.,$$
by arbitrariness of $\P^{'}$.

\vspace{0.5em}
$\mathrm{(ii)}$ We now show that
$$C_{t_1}^{\mathbb P}:=\underset{ \mathbb{P}^{'} \in \mathcal{P}_0(t_1,\mathbb{P},\F_{+}) }{\esup^\mathbb{P}}\mathbb E^{\mathbb{P}^{'}}\left[\left.\left(K_{t_2}^{\mathbb P^{'}}-K_{t_1}^{\mathbb P^{'}}\right)^p\right|\Fc_{t_1}^+\right]<+\infty,\text{ }\mathbb P-a.s.$$
First of all, we have by definition
\begin{align*}
\left(K_{t_2}^{\mathbb P^{'}}-K_{t_1}^{\mathbb P^{'}}\right)^p\leq &\ C\left(\underset{t_1\leq t\leq t_2}{\sup}\abs{Y_t}^p+\left(\int_{t_1}^{t_2}\abs{\widehat f^{\P^{'},0}}ds\right)^p+\left(\int_{t_1}^{t_2}\No{(\widehat a^{1/2}_s)^\top Z_s}ds\right)^p\right)\\
&+C\left(\abs{\int_{t_1}^{t_2}Z_s\cdot dX_s^{c,\P^{'}}}^p+\abs{\int_{t_1}^{t_2}dM_s^{\P^{'}}}^p\right),
\end{align*}
for some constant $C> 0$,
so that we obtain by BDG inequalities
\begin{align}\label{eq:estimK}
\E^{\P^{'}}\left[
\left(K_{t_2}^{\mathbb P^{'}}-K_{t_1}^{\mathbb P^{'}}\right)^p
\right]\leq &\ C\left(\phi^{p,\kappa}_f+\No{Y}_{\mathbb D^{p}_0}^p+\No{Z}_{\H^p_0}^p+\underset{\P\in\mathcal P_0}{\sup}\E^{\P}\left[[M^{\P}]_T^{\frac p2}\right]\right),
\end{align}
for some other constant $C>0$ and hence $C^{\P}_{t_1} < + \infty$, $\P-$a.s. Next, we claim that the family
$$\left\{\mathbb E^{\mathbb{P}^{'}}\left[\left.\left(K_{t_2}^{\mathbb P^{'}}-K_{t_1}^{\mathbb P^{'}}\right)^p\right|\Fc_{t_1}^+\right],\ \mathbb{P}^{'} \in \mathcal{P}_0(t_1,\mathbb{P},\F_{+})\right\},$$
is upward directed. Indeed, this can be proved exactly as in the proof of {Lemma} \ref{lemma:rep}. For the ease of the reader, we give the details again. Let us consider $(\P^1,\P^2)\in  \mathcal{P}_0(t_1,\mathbb{P},\F^{+})\times \mathcal{P}_0(t_1,\mathbb{P},\F^{+})$, and let us define the following subsets of $\Omega$
$$A_1:=\left\{\omega\in\Omega: \mathbb E^{\mathbb{P}^{1}}\left[\left.\left(K_{t_2}^{\mathbb P^{1}}-K_{t_1}^{\mathbb P^{1}}\right)^p\right|\Fc_{t_1}^+\right](\omega)>\mathbb E^{\mathbb{P}^{2}}\left[\left.\left(K_{t_2}^{\mathbb P^{2}}-K_{t_1}^{\mathbb P^{2}}\right)^p\right|\Fc_{t_1}^+\right](\omega)\right\},$$
 \vspace{-0.5em}
 $$A_2:=\Omega\backslash A_1.$$
 Then, $A_1,A_2\in\mathcal F_{t_1}^{\P+}$, and we can define the following probability measure on $(\Omega,\Fc_T)$
 $$\P^{1,2}(B):=\P^1(A_1\cap B)+\P^2(A_2\cap B),\ \text{for any $B\in\Fc_T$}.$$
 By Assumption \ref{assum:main}(v), we know that $\P^{1,2}\in\mathcal P_0$, and by definition, we further have $ \P^{1,2}\in\mathcal{P}_0(t_1,\mathbb{P},\F_{+})$ as well as, $\P-a.s.$,
 $$\mathbb E^{\mathbb{P}^{1,2}}\left[\left.\left(K_{t_2}^{\mathbb P^{1,2}}-K_{t_1}^{\mathbb P^{1,2}}\right)^p\right|\Fc_{t_1}^+\right]=\mathbb E^{\mathbb{P}^{1}}\left[\left.\left(K_{t_2}^{\mathbb P^{1}}-K_{t_1}^{\mathbb P^{1}}\right)^p\right|\Fc_{t_1}^+\right]\vee\mathbb E^{\mathbb{P}^{2}}\left[\left.\left(K_{t_2}^{\mathbb P^{2}}-K_{t_1}^{\mathbb P^{2}}\right)^p\right|\Fc_{t_1}^+\right],$$
 which proves the claim.
 
 \vspace{0.5em}
\noindent Therefore, by classical results for the essential supremum (see e.g. Neveu \cite{neveu1975discrete}), there exists a sequence $(\P^n)_{n\geq 0}\subset \mathcal{P}_0(t_1,\mathbb{P},\F_{+})$ such that
 $$\underset{ \mathbb{P}^{'} \in \mathcal{P}_0(t_1,\mathbb{P},\F_{+}) }{\esup^\mathbb{P}}\mathbb E^{\mathbb{P}^{'}}\left[\left.\left(K_{t_2}^{\mathbb P^{'}}-K_{t_1}^{\mathbb P^{'}}\right)^p\right|\Fc_{t_1}^+\right]=\underset{n\to \infty}{\lim}\uparrow \E^{\mathbb{P}_{n}}\left[\left.\left(K_{t_2}^{\mathbb P_{n}}-K_{t_1}^{\mathbb P_{n}}\right)^p\right|\Fc_{t_1}^+\right]$$
Then using \reff{eq:estimK} and the monotone convergence theorem under $\P$, we deduce that
 \begin{align*}
 \E^\P\left[C_{t_1}^\P\right]&\leq \underset{n\to \infty}{\lim}\uparrow\E^\P\left[\E^{\mathbb{P}_{n}}\left[\left.\left(K_{t_2}^{\mathbb P_{n}}-K_{t_1}^{\mathbb P_{n}}\right)^p\right|\Fc_{t_1}^+\right]\right]\\
 &\leq C\left(\phi^{p,\kappa}_f+\No{Y}_{\mathbb D^{p}_0}^p+\No{Z}_{\H^p_0}^p+\underset{\P\in\mathcal P_0}{\sup}\E^{\P}\left[[M^{\P}]_T^{\frac p2}\right]\right)<+\infty,
 \end{align*}  
 which provides the desired result.
 
 \vspace{0.5em}
$\mathrm{(iii)}$ We now prove the reverse inequality. Since we will use a linearization argument, we work on the enlarged space, remembering that this is without loss of generality by Lemma \ref{lemma:deuxbsde}. Fix $\mathbb P\in\mathcal P_{0}$.  For every $\mathbb P^{'}\in { \mathbb{P}^{'} \in \mathcal{P}_0(t_1,\mathbb{P},\F_{+}) }$, 
we extend the definition of $(Y,Z,(M^\P)_{\P\in\mathcal P_0},(K^\P)_{\P\in\mathcal P_0})$ on $\Omb$ as in \eqref{eq:extension_def}, 
and denote $
	\delta Y:=Y-\widetilde{\mathcal Y}^{\mathbb P^{'}\otimes\P_0},$ $\delta Z:=Z-\widetilde{\mathcal Z}^{\mathbb P^{'}\otimes\P_0}$ and  $\delta M^{\mathbb P^{'}}:=M^{\mathbb P^{'}}-\widetilde{\Mc}^{\mathbb P^{'}\otimes\P_0}.$
By Assumption \ref{assum:main}(i), there exist two bounded processes $\lambda^{\P^{'}}$ and $\eta^{\P^{'}}$ such that for all $t_1\leq t\leq t_2$, $\mathbb P^{'}\otimes\P_0-a.s$
{$$\delta Y_t=\int_t^{t_2}\left(\lambda^{\P^{'}}_s\delta Y_s+\eta^{\P^{'}}_s\cdot (\widehat{a}_s^{1/2})^\top \delta Z_s\right)ds-\int_t^{t_2}\delta Z_s\cdot \widehat{a}_s^{1/2}dW^{\P^{'}}_s-\int_t^{t_2}d\left(\delta M_{s}^{\mathbb P^{'}}- K_{s}^{\mathbb P^{'}}\right).$$}
Define for $t_1\leq t\leq t_2$ the following continuous process
\begin{equation} \label{eq:def_Delta}
	\Delta^{\P^{'}}_t:=\exp\left(\int_{t_1}^t\left(\lambda^{\P^{'}}_s-\frac12\No{\eta^{\P^{'}}_s}^2\right)ds-\int_{t_1}^t\eta^{\P^{'}}_s\cdot dW^{\P^{'}}_s\right),\text{ }\mathbb P^{'}\otimes\P_0-a.s.
\end{equation}
Note that since $\lambda^{\P^{'}}$ and $\eta^{\P^{'}}$ are bounded, we have for all $p\geq 1$, for some constant $C_p>0$, independent of $\P^{'}$
\begin{equation}
\label{truc.M}
\mathbb E^{\mathbb P^{'}\otimes\P_0}\left[\left.\underset{t_1\leq t\leq t_2}{\sup}(\Delta^{\P^{'}}_t)^p+\underset{t_1\leq t\leq t_2}{\sup}({\Delta^{\P^{'}}_t})^{-p}\right|\overline\Fc_{t_1}^{+}\right]\leq C_p,\text{ }\mathbb P^{'}\otimes\P_0-a.s.
\end{equation}
Then, by It\^o's formula, we obtain
\begin{equation}
\label{brubru}
\delta Y_{t_1}=\mathbb E^{\mathbb P^{'}\otimes\P_0} \left[\left.\int_{t_1}^{t_2}\Delta^{\P^{'}}_t d K_t^{\mathbb P^{'}}\right|\overline\Fc_{t_1}^{+}\right].
\end{equation}
because the martingale terms vanish when taking conditional expectation. We therefore deduce
\begin{align*}
\delta Y_{t_1}&\leq \left(\E^{\P^{'}\otimes\P_0}\left[\left.\underset{t_1\leq t\leq t_2}{\sup}\abs{\Delta^{\P^{'}}_t}^{\frac{p+1}{p-1}}\right|\overline\Fc_{t_1}^{+}\right]\right)^{\frac{p-1}{p+1}}\left(\E^{\P^{'}\otimes\P_0}\left[\left.\left(K^{\P^{'}}_{t_2}-K^{\P^{'}}_{t_1}\right)^{\frac{p+1}{2}}\right|\overline\Fc_{t_1}^{+}\right]\right)^{\frac{2}{p+1}}\\
&\leq C \left(C_{t_1}^{\P'}\right)^{\frac{1}{p+1}}\left(\E^{\P^{'}\otimes\P_0}\left[\left.K^{\P^{'}}_{t_2}-K^{\P^{'}}_{t_1}\right|\overline\Fc_{t_1}^{+}\right]\right)^{\frac{1}{p+1}}.
\end{align*}
Remember $Y$ and $K^{\P^{'}}$ are extended on $\Omb$ as in \eqref{eq:extension_def}, then it only depends on $X$ and not on $B$. 
Going back now to the canonical space $\Om$, it follows by Lemma \ref{lemma:deuxbsde} that
$$
	\delta Y_{t_1}' := Y_{t_1} - \widetilde \Yc^{\P'}_{t_1} \le C \left(C_{t_1}^{\P'}\right)^{\frac{1}{p+1}}\left(\E^{\P^{'}}\left[\left.K^{\P^{'}}_{t_2}-K^{\P^{'}}_{t_1}\right| \Fc_{t_1}^{+}\right]\right)^{\frac{1}{p+1}}.
$$
By arbitrariness of $\P^{'}$, we deduce thanks to \reff{2bsde.minK} that
$$Y_{t_1}-\underset{ \mathbb{P}^{'} \in \mathcal{P}_0(t_1,\mathbb{P},\F_{+}) }{ \esup^{\mathbb P} } \mathcal Y_{t_1}^{\mathbb{P}^{'}}(t_2,Y_{t_2})\leq 0,\ \P-a.s.$$
Finally, the uniqueness of $Y$ is immediate by the representation \reff{eq:2bsderep}. Then, since 
$${\langle Y,X \rangle_t}=\int_0^t\widehat a_s Z_sds,\ \P-a.s.,$$
$Z$ is also uniquely defined, $\widehat a_tdt\otimes\Pc_0-q.s.$ We therefore deduce that the processes $M^\P-K^\P$ are also uniquely defined for any $\P\in\mathcal P_0$. But, since they are $(\F^{\P}_+,\P)-$super--martingales, such that in addition $(K_t^\P,M_t^\P)\in \L^p_0(\F^{\P}_+,\P)\times \L^p_0(\F^{\P}_+,\P)$ for any $t\in[0,T]$, and since $K^\P$ is $\F^{\P}_+-$predictable, the uniqueness of $M^\P$ and $K^\P$ is a simple consequence of the uniqueness in the Doob--Meyer decomposition of these super--martingales.
\ep

\vspace{0.5em}
\noindent With the previous theorem in hand, the following comparison result is an immediate consequence of the corresponding one for BSDEs (see for instance Lemma \ref{lemma:comp} in the Appendix)
\begin{Theorem}
For $i=1,2$, let $f^i$ and $\xi^i$ be respectively a generator map and a terminal condition satisfying the required conditions in Theorem \ref{th:mainex}.
Let also $Y^i$ be the first component of the solution to the {\rm 2BSDE} with generator $f^i$ and terminal condition $\xi^i$. Suppose in addition that for any $\P\in\Pc_0$ we have
\begin{itemize}
\item[$(i)$] $\xi^1\leq \xi^2,\ \P-a.s.$
\item[$(ii)$] $\widehat{f}^{1,\P}_s(y_s^1, (\widehat a_s^{1/2})^\top z_s^1)\geq \widehat{f}^{2,\P}_s(y_s^2, (\widehat a_s^{1/2})^\top z_s^2),\ ds\times d\P-a.e., \text{ on }[0,T]\times\Omega$, where for $i=1,2$, $(y^i, z^i)$ are the first two components of the solution of the {\rm BSDE} under $\P$ with generator $\widehat{f}^{i,\P}$ and terminal condition $\xi^i$.
\end{itemize}
Then we have $Y^1_t\leq Y^2_t$, $t\in[0,T]$, $\Pc_0-q.s.$
\end{Theorem}
\subsection{A priori estimates and stability}
In this section, we give {\it a priori} estimates for 2BSDEs, which, as in the case of the classical BSDEs, play a very important role in the study of associated numerical schemes for instance. The proofs are actually based on the general results given very recently in \cite{bouchard2015unified}.
\vspace{0.5em}

\begin{Theorem}\label{estimatesref}
	Let Assumptions \ref{assum:main} and \ref{assum:integ} hold true,
	and $(Y,Z,(M^\P)_{\P\in\mathcal P_0},(K^\P)_{\P\in\mathcal P_0})$ be a solution to the {\rm 2BSDE} \reff{2bsde}. Then, there exists a constant $C_\kappa$ depending only on $p$, $\kappa$, $T$ and the Lipschitz constant of $f$ such that
$$\No{Y}^p_{\mathbb D^{p}_0}+\No{Z}^p_{\mathbb H^{p}_0}+\underset{\mathbb P\in \mathcal P_{0}}{\sup}\mathbb E^\mathbb P\left[(K_T^\mathbb P)^p\right]+\underset{\mathbb{P} \in \mathcal{P}_{0}}{\sup}\mathbb E^{\mathbb P}\left[[M^\P]_T^{\frac p2}\right]\leq C_\kappa \left(\No{\xi}^p_{\mathbb L_0^{p,\kappa}}+\phi_f^{p,\kappa}\right).$$
\end{Theorem}

\proof
First, by Theorem \ref{representation}, we have for any $\P\in\Pc_0$, $Y_t=\underset{\P^{'}\in\Pc_0(t,\P,\F_{+})}{\esup^\P}\Yc^{\P^{'}}_{t}(T,\xi),\ \P-a.s.$
Furthermore, by Lemma \ref{lemma:estimbsde} (together with Lemma \ref{lemma:deuxbsde}), we know that there exists a constant $C$ (which may change from line to line) depending only on $\kappa$, $T$ and the Lipschitz constant of $\widehat f$, such that for all $\mathbb P$
\begin{equation}
\label{estimref}
\abs{\mathcal Y_t^\mathbb P(T,\xi)}\leq C \left(\mathbb E^\mathbb P\left[\left.\abs{\xi}^\kappa+\int_t^T\abs{\widehat f^{\P,0}_s}^\kappa ds\right|\mathcal F_t^+\right]\right)^{\frac{1}{\kappa}},\ \P-a.s.
\end{equation}
Hence, we deduce immediately $\No{Y}^p_{\mathbb D^{p}_0}\leq C\left(\No{\xi}^p_{\mathbb L_0^{p,\kappa}}+\phi_f^{p,\kappa}\right).$ Now, by extending the definition of $(Y, Z, (M^{\P})_{\P \in \Pc_0}, (K^{\P})_{\P \in \Pc_0})$ on the enlarged space $\Omb$ (see \eqref{eq:extension_def}), one has for every $\P\in\Pc_0$,
$$Y_t=\xi-\int_t^T\widehat f^\P_s(Y_s, (\widehat a^{1/2}_s)^\top Z_s)ds-\int_t^TZ_s\cdot\widehat a_s^{1/2}dW^{\P}_s-\int_t^TdM_s^{\P}+\int_t^TdK^\P_s,\ \P\otimes\P_0-a.s.$$
Then for every $\P \in \Pc_0$, 
$(Y, Z, M^{\P}, K^{\P})$ can be interpreted as a super--solution of a {\rm BSDE} in the enlarged space $ \Omb$. We can therefore use Theorem 2.1 of \cite{bouchard2015unified} (notice that the constants appearing there do not depend on the underlying probability measure) to obtain the required estimates.
Noticing once again that the norms of $Z$, $K^\P$ and $M^\P$ are the same on the enlarged space $\Omb$ or on $\Omega$, it follows then
\begin{align*}
&\No{Z}^p_{\mathbb H^{p}_0}+\underset{\mathbb P\in \mathcal P_{0}}{\sup}\mathbb E^\mathbb P\left[(K_T^\mathbb P)^p\right]+\underset{\mathbb{P} \in \mathcal{P}_{0}}{\sup}\mathbb E^{\mathbb P}\left[[M^\P]_T^{\frac p2}\right]\\
&\leq C\left(\No{\xi}^p_{\mathbb L_0^{p,\kappa}}+\phi_f^{p,\kappa}+\No{\xi}^p_{\L^p_0}+\underset{\P\in\mathcal P_0}{\sup}\E^\P\left[\int_0^T\abs{\widehat f^{\P,0}_s}^pds\right]\right)\leq C\left(\No{\xi}^p_{\mathbb L_0^{p,\kappa}}+\phi_f^{p,\kappa}\right),
\end{align*}
for some constant $C>0$,
where we used the fact that by definition $\No{\xi}^p_{\L^p_0}\leq \No{\xi}^p_{\mathbb L_0^{p,\kappa}}$ and $\underset{\P\in\mathcal P_0}{\sup}\E^\P\left[\int_0^T\abs{\widehat f^{\P,0}_s}^pds\right]\leq \phi_f^{p,\kappa}.$
\ep

\vspace{0.5em}
\noindent Next, we also have the following estimates for the difference of two solutions of 2BSDEs, which plays a fundamental role for stability properties.

\begin{Theorem}\label{th:estimates2}
	We are given two generators $f^1$ and $f^2$ and terminal variables $\xi^1$ and $\xi^2$ satisfying Assumptions \ref{assum:main} and \ref{assum:integ}.
	Let $(Y^i, Z^i,(M^{i,\P})_{\P\in\Pc_0},(K^{i,\P})_{\P\in\Pc_0})$ be a solution to the {\rm 2BSDE} with generator $f^i$ and terminal condition $\xi^i$, for $i=1,2$. 
Define
\begin{align*}
\phi^{p,\kappa}_{f^1,f^2}&:=\underset{\P\in\mathcal P_0}{\sup}\E^\P\left[\underset{0\leq t\leq T}{\esup^\P}\ \E^{\P}\left[\left.\left(\int_0^T\abs{\widehat{f}^{1,\P}_s-\widehat{f}^{2,\P}_s}^\kappa(y_s^{1,\P}, (\widehat a_s^{1/2})^\top z_s^{1,\P})ds\right)^{\frac{p}{\kappa}}\right|\Fc^+_t\right]\right]\\
\psi^{p}_{f^1,f^2}&:=\underset{\P\in\mathcal P_0}{\sup}\E^\P\left[\int_0^T\abs{\widehat{f}^{1,\P}_s-\widehat{f}^{2,\P}_s}^p(Y^1_s, (\widehat a_s^{1/2})^\top Z_s^{1})ds\right].
\end{align*}
Then, there exists a constant $C_\kappa$ depending only on $\kappa$, $T$ and the Lipschitz constant of $f^1$ and $f^2$ such that
\begin{align*}
&\No{Y^1-Y^2}^p_{\mathbb D^{p}_0}\leq C_\kappa\left(\No{\xi^1-\xi^2}^p_{\L_0^{p,\kappa}}+\phi^{p,\kappa}_{f^1,f^2}\right)\\
&\No{Z^1-Z^2}^p_{\mathbb H^{p}_0}+\underset{\mathbb P\in \mathcal P_{0}}{\sup}\mathbb E^\mathbb P\left[[N^{1,\P}-N^{2,\P}]_T^{\frac p2}\right]\leq C_\kappa \left(\No{\xi^1-\xi^2}^p_{\mathbb L_0^{p,\kappa}}+\phi^{p,\kappa}_{f^1,f^2}+\No{\xi^1-\xi^2}^{\frac p2\wedge(p-1)}_{\mathbb L_0^{p,\kappa}}\right.\\
&\hspace{18.5em}\left.+\; \psi^{p}_{f^1,f^2}+\left(\phi^{p,\kappa}_{f^1,f^2}\right)^{\frac p2\wedge(p-1)}\right),
\end{align*}
where we have once more defined $N^{i,\P}:=M^{i,\P}-K^{i,\P}$ for any $\P\in\Pc_0$, $i=1,2$.
\end{Theorem}
\proof
First of all, by Lemma \ref{lemma:estimbsde} (together with Lemma \ref{lemma:deuxbsde}), we know that there exists a constant $C$ depending only on $\kappa$, $T$ and the Lipschitz constant of $\widehat f$, such that for all $\mathbb P\in\Pc_0$, we have $\P-a.s.$
\begin{equation}
\label{estimref2}
\abs{y_t^{1,\mathbb P}-y_t^{2,\mathbb P}}\leq C \left(\mathbb E^\mathbb P\left[\left.\abs{\xi^1-\xi^2}^\kappa+\int_t^T\abs{\widehat f^{1,\P}_s-\widehat f^{2,\P}_s}^\kappa(y_s^{1,\P}, (\widehat a_s^{1/2})^\top z_s^{1,\P}) ds\right|\mathcal F_t^+\right]\right)^{\frac{1}{\kappa}}.
\end{equation}
This immediately provides the estimate for $Y^{1}-Y^{2}$ by the representation formula $\reff{eq:2bsderep}$ and the definition of the norms and of $\phi^{p,\kappa}_{f^1,f^2}$. 
Next, we argue exactly as in the proof of Theorem \ref{estimatesref} by working on the enlarged space $\Omb$ and using now Theorem $2.2$ of \cite{bouchard2015unified} to obtain the required estimates.
\ep

\subsection{Existence through dynamic programming}


In this section, we will show that $\Ych^+$ defined in Section \ref{sec:1} is indeed a solution to the {\rm 2BSDE} \reff{2bsde}, thus completing the proof of Theorem \ref{th:mainex}.

\vspace{0.5em}
\noindent Recall that $\Ych^+$ is defined by \eqref{eq:def_Ychplus}, and one has processes
$(Z,(M^\P)_{\P\in\mathcal P_0},(K^\P)_{\P\in\mathcal P_0})\in \mathbb H^{p}_0(\F^{\Pc_0})\times \mathbb M^p_0((\F^{\P+})_{\P\in\Pc_0})\times \mathbb I^p_0((\F^{\P+})_{\P\in\Pc_0})$
given by Lemma \ref{lemma:semi},
so that the only thing left for us is to show that the family $(K^\P)_{\P\in\Pc_0}$ satisfies the minimality condition \reff{2bsde.minK}. 

\vspace{0.5em}
\noindent We again extend the definition of $(Y, Z,(M^\P)_{\P\in\mathcal P_0},(K^\P)_{\P\in\mathcal P_0})$
and $\widehat \Yc^+$, $\Yc^{\P'}(T, \xi)$ on $\Omb$ as in \eqref{eq:extension_def} (recall also Lemma \ref{lemma:deuxbsde}).
Then by \reff{brubru}, denoting $\delta \Ych^+:=\Ych^+-\Yc^{\P^{'}}(T,\xi)$, we have for any $t\in[0,T]$, for any $\P\in\Pc_0$ and any $\P^{'}\in\Pc_0(t,\P,\F_+)$
\begin{equation*}
\delta \Ych^+_{t}=\mathbb E^{\mathbb P^{'}\otimes\P_0} \left[\left.\int_{t}^{T}\Delta^{\P^{'}}_s d K_s^{\mathbb P^{'}}\right|\overline\Fc_{t}^+\right]\geq \mathbb E^{\mathbb P^{'}\otimes\P_0} \left[\left.\underset{t\leq s\leq T}{\inf}\Delta^{\P^{'}}_s \left(K_T^{\mathbb P^{'}}-K_t^{\mathbb P^{'}}\right)\right|\overline\Fc_{t}^+\right],\ \P-a.s.,
\end{equation*}
where $\Delta^{\P'}$ is defined in \eqref{eq:def_Delta}. Denote $\Kc^{\P^{'}}_t:=\mathbb E^{\mathbb P^{'}\otimes\P_0}[ K_T^{\mathbb P^{'}}-K_t^{\mathbb P^{'}}\|\overline\Fc_{t}^+].$
We therefore have
\begin{align*}
\Kc^{\P^{'}}_t\leq&\  \bigg(\mathbb E^{\mathbb P^{'}\otimes\P_0} \bigg[\underset{t\leq s\leq T}{\inf}\Delta^{\P^{'}}_s \Big(K_T^{\mathbb P^{'}}-K_t^{\mathbb P^{'}}\Big)\bigg|\overline{\Fc}_{t}^+\bigg]\bigg)^{\frac12}  \\
&\times \Big(\E^{\mathbb P^{'}\otimes\P_0} \Big[\Big(K_T^{\mathbb P^{'}}-K_t^{\mathbb P^{'}}\Big)^p\Big|\overline\Fc_{t}^+\Big]\Big)^{\frac{1}{2p}}\bigg(\E^{\mathbb P^{'}\otimes\P_0} \bigg[\Big(\underset{t\leq s\leq T}{\inf}\Delta^{\P^{'}}_s\Big)^{-q}\bigg|\overline\Fc_{t}^+\bigg]\bigg)^{\frac{1}{2q}}\\
\leq &\ C\left(C_t^{\P'} \right)^{\frac{1}{2p}}\left(\delta \Ych^+_t\right)^{\frac12},
\end{align*}
with $q>1$ such that $\frac{1}{p}+\frac{1}{q}=1$.

\vspace{0.5em}\noindent
Notice that by definition $K^{\P^{'}}$ (defined on $\Omb$) only depends on $X$ and not on $B$, 
so that we can go back to $\Om$ and obtain
$$
	\E^{\P'}  \left[\left. K_T^{\mathbb P^{'}}-K_t^{\mathbb P^{'}}\right| \Fc_{t}^+\right]
	~\leq~
	C\left(C_t^{\P'}\right)^{\frac{1}{2p}}\left(\delta \Ych^+_t\right)^{\frac12}.
$$
Then the result follows immediately thanks to Lemma \ref{lemma:rep}.

\begin{Remark}
{\rm For other classes of 2BSDEs with possibly non--Lipschitz generator, such as 2BSDEs under a monotonicity condition \cite{possamai2013second1}, quadratic 2BSDEs \cite{possamai2013second}, second--order reflected BSDEs \cite{matoussi2013second,matoussi2014second}, or doubly stochastic 2BSDEs \cite{matoussi2014probabilistic}, if a Doob--Meyer decomposition for the corresponding nonlinear super--martingales is available under any probability measure in the set $\Pc_0$, then together with Proposition \ref{others}, we can generalize the wellposedness result in Theorem \ref{th:mainex} to these classes of 2BSDEs when there is no regularity condition on the terminal condition and the generator. In particular, all probability measures in the non--dominated set considered in the articles above do satisfy this property, which means that our result extends directly to their context.}
\end{Remark}


\section{Nonlinear optional decomposition and super--hedging duality}\label{sec:4}
In this section, we show that under an additional assumption on the sets $\Pc_0$, basically stating that it is rich enough, we can give a different definition of second-order BSDEs, which is akin to a nonlinear optional decomposition theorem, as initiated by \cite{el1995dynamic,follmer1997optional,kramkov1996optional} in a dominated model framework, and more recently by \cite{nutz2015robust} for non--dominated models.

\subsection{Saturated 2BSDEs}
We introduce the following definition.

\begin{Definition}
The set $\Pc_0$ is said to be saturated if, when $\P\in\Pc_0$, we have $\Q\in\Pc_0$ for every probability measure $\Q$ on $(\Omega, \Fc)$ which is equivalent to $\P$ and under which $X$ is local martingale.
\end{Definition}

\begin{Remark}\label{rem.sat}\rm{
It can be readily checked that for instance the set $\overline \Pc_S$ of \cite{soner2012wellposedness}, whose measures only change the volatility of $X$ is saturated. We can also accommodate cases where $X$ has a drift, provided the latter belongs to the range of $\widehat a$, so that a Girsanov transformation can be applied. Similar examples can be obtained using weak or relaxed formulations, instead of the strong one used in $\overline \Pc_S$.}
\end{Remark}

\noindent We give now an alternative definition for 2BSDEs of the form
\begin{equation}
Y_t=\xi -\int_t^T\widehat{f}^\P_s(Y_s,  (\widehat{a}_s^{1/2})^\top Z_s)ds -  {\left(\int_t^TZ_s\cdot  dX_s^{c,\P}\right)^\P} + K^{\mathbb P}_T-K^{\mathbb P}_t,\ 0\leq t\leq T.
\label{2bsdesat}
\end{equation}

\begin{Definition}\label{def:2}
	We say $(Y,Z,(K^\P)_{\P\in\mathcal P_0})\in \mathbb D^{p}(\F^{\Pc_0+}) \times \mathbb H^{p}(\F^{\Pc_0})\times (\mathbb I^{o,p}((\F^{\P}_+))_{\P\in\Pc_0}$ is a saturated solution to ${\rm 2BSDE}$ \reff{2bsdesat} if equation \reff{2bsdesat} holds $\Pc_0-q.s.$ and if the family $\left\{K^{\mathbb P}, \mathbb P \in \mathcal P_{0}\right\}$ satisfies the minimality condition \reff{2bsde.minK}.
\end{Definition}

\begin{Remark}{\rm
In the above definition, two changes have occurred. First, the orthogonal martingales $M^\P$ have disappeared, and the non--decreasing processes $K^\P$ are assumed to be $\F^{\P}_+-$optional instead of being predictable.}
\end{Remark}

\noindent We then have the following result.
\begin{Theorem}
	Let Assumptions \ref{assum:main}, \ref{assum:integ} hold.
	Assume in addition that the set $\Pc_0$ is saturated. 
	Then there is a unique {\it saturated} solution of the {\rm 2BSDE} \reff{2bsdesat}.
\end{Theorem}

\proof
By Theorem \ref{th:mainex}, we know that the following {\rm 2BSDE} is well--posed, $ \Pc_0-q.s.$
$$Y_t=\xi -\int_t^T\widehat{f}^\P_s(Y_s, (\widehat{a}_s^{1/2})^\top Z_s)ds -\left(\int_t^T Z_s\cdot dX^{c,\P}_s\right)^\P-\int_t^T dM^{\mathbb P}_s + K^{\mathbb P}_T-K^{\mathbb P}_t,\ t\in[0,T].$$
In particular, this means that the process
$$Y_\cdot-\int_0^\cdot\widehat{f}^\P_s(Y_s, (\widehat{a}_s^{1/2})^\top Z_s)ds,$$
is a $(\F^{\P}_+,\P)-$supermartingale in $\D^p_0(\F^{\P}_+,\P)$ for every $\P\in\Pc_0$. 
Since $\Pc_0$ is saturated, it follows by Theorem 1 of \cite{follmer1997optional} (see also Theorem 3.1 of \cite{follmer1997optional2}), that there exists a $\F-$predictable process $\widetilde Z^{\P}$ such that
$$Y_\cdot-\int_0^\cdot\widehat{f}^\P_s(Y_s, (\widehat{a}_s^{1/2})^\top Z_s)ds-\int_0^\cdot \widetilde Z_s^{\P} \cdot d X^{c,\P}_s\text{ is non--increasing, $\P-a.s.$, for every $\P\in\Pc_0$.}$$
Hence, we can write
$$Y_t=\xi -\int_t^T\widehat{f}^\P_s(Y_s, (\widehat{a}_s^{1/2})^\top Z_s)ds -\left(\int_t^T \widetilde Z_s^{\P}\cdot dX^{c,\P}_s\right)^\P+ \tilde K^{\mathbb P}_T-\tilde K^{\mathbb P}_t,\ 0\leq t\leq T,\ \Pc_0-q.s.,$$
where for any $\P\in\Pc_0$, $\tilde K^\P$ is c\`adl\`ag, non--decreasing $\P-a.s.$ and $\F^{\P}_+-$optional. Moreover, by identification of the martingale parts, we deduce that we necessarily have $\widetilde Z^{\P}=Z$, $\widehat a_t dt\times\Pc_0-q.s.$ Finally, following the same arguments as in the proof of Theorem \ref{estimatesref}, we deduce that 
	$(\tilde K^\P)_{\P\in\Pc_0}\in (\mathbb I^{o,p}((\F^{\P}_+))_{\P\in\Pc_0}$, which ends the proof.
\ep

\subsection{A super--hedging duality in uncertain, incomplete and nonlinear markets}
\label{subsec:duality}

The result of the previous section finds an immediate application to the so-called problem of robust super--hedging. Before discussing the related results in the literature, let us explain exactly what the problem is.

\vspace{0.5em}
	\noindent We consider a standard financial market (possibly incomplete) consisting of a non--risky asset and $n$ risky assets whose dynamics are uncertain but given as solutions of controlled SDEs $S^{t,\om, \nu}$ in Remark \ref{remark:assum_meas_selec}.
	The collection of the law of these dynamics, denoted by $\Pc^{\Uc}_0$, satisfies Assumption \ref{assum:main} $\mathrm{(iii)-(v)}$.
	We assume in addition that $\Pc^{\Uc}_0$ is saturated (recall Remark \ref{rem.sat}).

\vspace{0.5em}
\noindent A portfolio strategy is then defined as a $\R^n-$valued and $\F^{\Pc^{\Uc}_0}-$predictable process $(Z_t)_{t\in[0,T]}$, such that $Z^i_t$ describes the number of units of asset $i$ in the portfolio of the investor at time $t$. It is well--known that under some constrained cases, the wealth $Y^{y_0,Z}$ associated to the strategy $Z$ and initial capital $y_0\in\R$ can be written, for every $\P\in\Pc^\Uc_0$, as
$$ {Y^{y_0,Z}_t:=y_0+\int_0^t\widehat f^\P_s(Y_s^{y_0,Z}, (\widehat{a}_s^{1/2})^\top Z_s)ds +\int_0^tZ_s\cdot dX_s^{c,\P},\ t\in[0,T],\ \P-a.s. }$$
For instance, the classical case corresponds to
\begin{equation}\label{eq:exemple}
\widehat f^\P_s(y,z)=r_sy+z\cdot\theta^\P_s,
\end{equation}
where $r_s$ is the risk--free rate of the market and $\theta^\P$ is the risk premium vector under $\P$, 
defined by $\theta^\P_s:= \big(\widehat a_s^{1/2}\big)^{\oplus} (b_s^\P-r_s{\bf 1}_n)$, where $\big(\widehat a_s^{1/2}\big)^{\oplus} $ denotes the Moore--Penrose pseudoinverse of $\widehat a^{1/2}_s$.

\vspace{0.5em}
\noindent The simplest example of a nonlinear $\widehat f^\P$ corresponds to the case where there are different lending and borrowing rates $\underline r_t\leq \overline r_t$, in which (see Example 1.1 in \cite{el1997backward})
$$\widehat f^\P_s(y,z)=\underline r_sy+z\cdot\theta^\P_s-(\overline r_s-\underline r_s)\left(y-z\cdot {\bf 1}_n\right)^-.$$
We will always assume that $\widehat f^\P$ satisfies our standing hypotheses in Assumptions \ref{assum:main} and \ref{assum:integ}.

\vspace{0.5em}
\noindent Let us now be given some Borel random variable $\xi\in \L^p(\Fc_T^{\Pc^{\Uc}_0})$. The problem of super--hedging $\xi$ corresponds to finding its super--replication price, defined as
$$P_{\rm sup}(\xi):=\inf\left\{y_0\in\R: \exists Z\in\mathcal H,\ Y_T^{y_0,Z}\geq \xi,\ \Pc^{\Uc}_0-q.s.\right\},$$
where the set of admissible trading strategies $\mathcal H$ is defined as the set of $\F^{\Pc_0^{\Uc}}-$predictable processes $Z$ such that in addition, $(Y^{y_0,Z}_t)_{t\in[0,T]}$ is a non--linear super--martingale under $\P$ for any $\P\in\Pc^{\Uc}_0$, in the sense that for any $0\leq s\leq t\leq T$
$$Y^{y_0,Z}_s\geq \Yc^\P_s(t,Y^{y_0,Z}_t),\ \P-a.s.$$
In the case where $\widehat f^\P$ corresponds to our first example \reff{eq:exemple} with $r=0$, and where the set of measures considered satisfy the predictable martingale representation property (that is the financial market is complete under any of the measures considered), this super--hedging price has been thoroughly studied in the recent literature, see among others \cite{avellaneda1995pricing,biagini2015robust,denis2006theoretical,dolinsky2014martingale,neufeld2013superreplication,nutz2012superhedging,nutz2015optimal,peng2012nonlinear,possamai2013robust,soner2011martingale,soner2013dual,song2011some}. The extension to possibly incomplete markets has been carried out notably by \cite{bouchard2015arbitrage} in discrete--time and more recently by \cite{nutz2015robust} in continuous time for models possibly incorporating jumps. Our result below extends all the results for continuous processes to markets with nonlinear portfolio dynamics. Of course, the same proof would go through for the more general jump case, provided that a {\rm 2BSDE} theory, extending that of \cite{kazi2015second,kazi2015second2}, is obtained in such a setting.

\begin{Theorem}
	Suppose that
	Assumptions \ref{assum:main}, \ref{assum:integ} hold and the set $\Pc_0$ is saturated. 
	Let $(Y,Z)$ be the first two components of the saturated solution of the {\rm 2BSDE} with generator $\widehat f^\P$ and terminal condition $\xi$. 
	Then
	$$P_{\rm sup}(\xi)=\underset{\P\in\Pc^{\Uc}_0}{\sup}\E^\P\left[Y_0\right],$$
	and $Z\in\Hc$ is a super--hedging strategy for $\xi$.
\end{Theorem}
\proof
First of all, assume that we have some $Z\in\Hc$ such that $Y_T^{y_0,Z}\geq \xi,\ \Pc^{\Uc}_0-q.s.$ Then, since $Y_T^{y_0,Z}$ is a non--linear super--martingale under $\P$ for any $\P\in\Pc^{\Uc}_0$, we have
$$y_0\geq \Yc^\P_0(T,Y^{y_0,Z}_T),\ \Pc_0^{\Uc}-q.s.$$
However, by the comparison result of Lemma \ref{lemma:comp} (together with Lemma \ref{lemma:deuxbsde}), we also have $\Yc^\P_0(T,Y^{y_0,Z}_T)\geq \Yc^\P_0(T,\xi)$, from which we deduce
$$y_0\geq \Yc^\P_0(T,\xi),\ \P-a.s.$$
In particular, for any $\P\in\Pc_0^{\Uc}$, we deduce that
$$y_0\geq \underset{\P^{'}\in\Pc_0^{\Uc}(0,\P,\F_+)}{\esup^\P}\Yc^\P_0(T,\xi)=Y_0,\ \P-a.s.,$$
where we have used Lemma \ref{lemma:rep}. It therefore directly implies, since $y_0$ is deterministic, that
$$y_0\geq \underset{\P\in\Pc^{\Uc}_0}{\sup}\E^\P\left[Y_0\right].$$
For the reverse inequality, let $(Y,Z,(K^\P)_{\P\in\mathcal P_0^{\Uc}})\in \mathbb D^{p}(\F^{U+,\Pc^{\Uc}_0}) \times \mathbb H^{p}(\F^{U,\Pc^{\Uc}_0})\times \mathbb I^{o,p}((\F^{\P}_+)_{\P\in\Pc^{\Uc}_0})$ be the unique saturated solution to the {\rm 2BSDE} with generator $\widehat f^\P$ and terminal condition $\xi$. Then, we have for any $\P\in\Pc_0^{\Uc}$
$$Y_0+\int_0^T\widehat f^\P_s(Y_s, (\widehat{a}_s^{1/2})^\top Z_s)ds+\int_0^TZ_s\cdot {dX^{c,\P}_s} =\xi+K_T^\P-K_t^\P\geq \xi,\ \P-a.s.$$
However, since $Y_0$ is only $\Fc_0^{\Pc^{\Uc}_0+}-$measurable, it is not, in general, deterministic, so that we cannot conclude directly. Let us nonetheless consider, for any $\P\in\Pc^{\Uc}_0$, $y_0^\P$ the smallest constant which dominates $Y_0$, $\P-a.s.$ We therefore want to show that for any $\P\in\Pc^{\Uc}_0$
$$y_0^\P\leq \underset{\P\in\Pc^{\Uc}_0}{\sup}\E^\P\left[Y_0\right],$$
which can be done by following exactly the same arguments as in the proof of Theorem 3.2 in \cite{nutz2015robust}. Finally, we do have $Z\in \Hc$, since by Lemma \ref{lemma:optsamp}, $Y$ is automatically a non--linear super--martingale for every $\P\in\Pc^{\Uc}_0$.
\ep

\section{Path--dependent PDEs}\label{sec:5}
\def\ch{\textsc{h}}

	In the context of stochastic control theory, using the dynamic programming principle,
	we can characterize the value function as a viscosity solution of PPDE.
	Recall that $\mu$, $\sigma$, $U$ as well as $\Uc$ are the same given in Section \ref{subsec:duality},
	we introduce a path--dependent PDE
	\begin{equation} \label{eq:PPDE}
		\partial_t v (t, \om)
		+
		G(t, \om, v(t,\om), \partial_{\om} v, \partial^2_{\om, \om}v)
		=
		0,
	\end{equation}
	where
	\begin{equation*}
		G(t, \om, y, z, \gamma)
		:=
		\sup_{u \in U} 
		\left\{
			f(t, \om, y , \sigma(\cdot) z, \mu, a \big)(t,\om,u) )
			+ \mu(t, \om, u) \cdot z
			+ \frac{1}{2} {\rm Tr}\big[a(t, \om, u) \gamma\big]
		\right\}.
	\end{equation*}
	As in the survey of Ren, Touzi and Zhang \cite{ren2014overview} (see also \cite{ren2014comparison}), one may define viscosity solutions of path dependent PDEs by using jets. For $\alpha \in \R,\ \beta \in \R^d,\ \gamma \in \S^d$, denote
	\begin{equation*}
		\phi^{\alpha, \beta, \gamma}(t,x):=\alpha t+ \beta \cdot x + \frac12{\rm Tr}\big[ \gamma(x x^T)\big]
		~\mbox{for all}~ 
		(t,x)\in \R^+\times \R^d,
	\end{equation*}
	where $A_1:A_2 := {\rm Tr}[A_1 A_2]$.
	Let $\mathrm{BUC}([0,T] \x \Om)$ denote the set of all bounded functions in $\Omega$ which are in addition uniformly continuous w.r.t. the metric $d$ defined by
	\begin{equation*}
		d \big( (t,\om), (t', \om') \big)
		:=
		\sqrt{ |t - t'|} + \| \om_{t \wedge \cdot} - \om'_{t' \wedge \cdot} \|_\infty.
	\end{equation*}
	Then define the semi--jets of a function $u\in \mbox{\rm BUC}([0,T] \x \Om)$ at the point $(t,\om)\in [0,T)\times\Om$ by
	\begin{align*}
		\underline \Jc u(t,\om) 
		~:=~ 
		\Big\{(\alpha,\beta,\gamma): u(t,\om)=\max_{\tau \in\Tc_{\ch_\delta}} \overline \Ec_{t,\om} [ u_\tau - \phi^{\alpha,\beta,\gamma,t,\om}_\tau],~\mbox{for some}~\delta>0\Big\}, \\
		\overline \Jc u(t,\om)
		~:=~
		\Big\{(\alpha, \beta, \gamma): u(t,\om)=\min_{\tau \in\Tc_{\ch_\delta}} \underline \Ec_{t,\om} [u_\tau - \phi^{\alpha,\beta,\gamma,t,\om}_\tau],~\mbox{for some}~\delta>0\Big\},
	\end{align*}
	where 
	$$
		\overline \Ec_{t,\om} [ \xi ] := \sup_{\P \in \Pc^{\Uc}_{t, \om}} \E^{\P}[ \xi],
		~
		\underline \Ec_{t,\om} [ \xi ] := \inf_{\P \in \Pc^{\Uc}_{t, \om}} \E^{\P}[ \xi],		
	$$
	and
	$\ch_\delta (\om'):=\delta \wedge \inf\{s\ge 0:|\om'_s|\ge \delta \}$ and  $\Tc_{\ch_\delta}$ denotes the collection of all $\F-$stopping times larger than $\ch_{\delta}$,
	and $\phi^{\alpha,\beta,\gamma, t,\om}_s : = \phi^{\alpha,\beta,\gamma}(s -t,X_s - \om_t)$.

	\begin{Definition}
		Let $u\in \mbox{\rm BUC}([0,T] \x \Om)$.

		\noindent {\rm (i)}~ $u$ is a $\Pc^{\Uc}-$viscosity sub--solution $($resp. super--solution$)$ of the path dependent PDE \eqref{eq:PPDE}, if at any point $(t,\om)\in [0,T)\times\Om$ it holds for all $(\alpha,\beta,\gamma)\in \underline \Jc u(t,\om)$ $($resp. $\overline\Jc u(t,\om) )$ that
		\begin{equation*}
			-\alpha - G(t,\om, u(t,\om),\beta,\gamma)\le \mbox{$($resp. $\ge)$}0.
		\end{equation*}

		\noindent {\rm (ii)}~ $u$ is a $\Pc^{\Uc}-$viscosity solution of PPDE \eqref{eq:PPDE}, if $u$ is both a $\Pc^{\Uc}-$viscosity sub--solution and a $\Pc^{\Uc}-$viscosity super--solution of \eqref{eq:PPDE}.
	\end{Definition}
\noindent Using the dynamic programming principle, and by exactly the same arguments as in \cite[Section 4.3]{ekren2016viscosity},
	we can characterize the value $Y$ as viscosity solution of the above PPDE.
	
	\begin{Theorem}
		Let Assumptions \ref{assum:main} and \ref{assum:integ} hold. 
		Suppose in addition that $\xi$ and $f$ are uniformly bounded,
		and $\om \longmapsto \xi(\om)$ and $\om \longmapsto f(t, \om, y, z, b, a)$ are uniformly continuous with respect to $\No{\cdot}_\infty$.
		Then the value function $v(t, \om) := Y_t(\om) \in \mathrm{BUC}([0,T] \x \Om)$ 
		and $v$ is a viscosity solution of PPDE \eqref{eq:PPDE}.
	\end{Theorem}

\noindent Of course, in order to have a complete characterization of the solution to a {\rm 2BSDE} as viscosity solution of the corresponding PPDE, the above result has to be complemented with a comparison theorem. In the case of fully nonlinear PPDE, such a result has been recently achieved by Ren, Touzi and Zhang \cite{ren2015comparison}. However, their main result Theorem 4.2 needs to consider viscosity sub--solutions and super--solutions in a smaller set than $\mathrm{BUC}([0,T] \x \Om)$. Namely, define for any $\ell>0$ the set $\mathrm{BUC}^\ell([0,T] \x \Om)$ of all bounded functions in $\Omega$ which are in addition uniformly continuous w.r.t. the metric $d^\ell$ defined by
	\begin{equation*}
		d^\ell \big( (t,\om), (t', \om') \big)
		:=
		\sqrt{ |t - t'|} +\No{\omega_{\cdot\wedge t}-\omega'_{\cdot\wedge t}}_{\ell},
	\end{equation*}
	where for any $\omega\in\Omega$
	$$\No{\omega}_{\ell}:=\left( \int_0^{T+1} \No{\om_{ s} }^\ell ds\right)^{\frac1\ell}.$$
 We then have the following result.
 	\begin{Theorem}
		Let Assumptions \ref{assum:main}, \ref{assum:integ} hold and let $G$ satisfy Assumption 4.1 of {\rm\cite{ren2015comparison}}. 
		Suppose in addition that $xi$ and $f$ are uniformly bounded 
		and $\om \longmapsto \xi(\om)$ and $\om \longmapsto f(t, \om, y, z, b, a)$ are uniformly continuous with respect to $\No{\cdot}_\ell$ for some $\ell\leq p$.
		Then the value function $Y_t(\om) $ 
		is the unique viscosity solution of PPDE \eqref{eq:PPDE} in $\mathrm{BUC}^\ell([0,T] \x \Om)$.
	\end{Theorem}
	
	\proof
	The only thing to prove here is that $Y$ does belong to $\mathrm{BUC}^\ell([0,T] \x \Om)$, since we can then apply immediately Theorem 4.2 of \cite{ren2015comparison}. However, this regularity can be obtained from classical {\it a priori} estimates for BSDEs, and arguments similar to the ones used in Example 7.1 of \cite{ren2015comparison}. 
	\ep
	
\begin{appendix}
\section{Technical results for BSDEs}
\noindent In this Appendix, we collect several results related to {\rm BSDE} theory which are used throughout the paper. We fix $r \in [0,T]$ and some $\P\in\Pc(r,\omega)$. A generator will here be a map $g:[r,T]\times\Omega\times\R\times\R^d\longrightarrow\R$ which is $\F_{+}-$progressively measurable and uniformly Lipschitz in $(y,z)$, satisfying
$$\E^\P\left[\int_{r}^T\abs{g_s(0,0)}^pds\right]<+\infty.$$
Similarly, a terminal condition will be a $\Fc_T-$measurable random variable in $\L^p_{r}(\Fc_T,\P)$. To state our results, we will actually need to work on the enlarged canonical space $\overline\Omega$, but we remind the reader that by Lemma \ref{lemma:deuxbsde}, it is purely a technical tool. Let $\overline\P:=\P\otimes\P_0$. We will then say that $(y,z,m)\in\D^p_{r}(\overline\F^{\overline\P+},\overline\P)\times\H^p_{r}(\overline\F^{\overline\P},\overline\P)\times\M^p_{r}(\overline\F^{\overline\P+},\overline\P)$ is a solution to the {\rm BSDE} with generator $g$ and terminal condition $\xi$ if
\begin{equation} \label{eq:BSDEp}
y_t=\xi(X_{\cdot}) -\int_t^Tg_s(y_s, (\widehat{a}_s^{1/2})^\top z_s)ds-\int_t^Tz_s\cdot \widehat a_s^{1/2}dW^{\P}_s-\int_t^Tdm_s,\ t\in[r,T],\ \overline\P-a.s.
\end{equation}
Similarly, if we are given a process $k\in\mathbb I^p_{r}(\overline\F^{\overline\P+},\overline\P)$, we call $(y,z,m,k)$ a super--solution of the {\rm BSDE} with generator $g$ and terminal condition $\xi$ if
\begin{equation}\label{eq:BSDE_supsol}
y_t=\xi-\int_t^Tg_s(y_s, (\widehat{a}_s^{1/2})^\top z_s)ds-\int_t^Tz_s\cdot \widehat a_s^{1/2}dW^{\P}_s-\int_t^Tdm_s+\int_t^Tdk_s,\ t\in[r,T],\ \overline\P-a.s.
\end{equation}

\begin{Lemma}[Estimates and stability]\label{lemma:estimbsde}
	Let the generator functions $g^i$ and the terminal condition $\xi^i$ for $i=1,2$ satisfy Assumption \ref{assum:main},
	we denote by $(y^{i},z^{i},m^{i})$ the solution of the {\rm BSDE} \eqref{eq:BSDEp} with $g^i$ and $\xi^i$.
	Then, for $\kappa\in(1,p]$, there exists some constant $C>0$ such that {for all $\F_+-$stopping time $\tau$ taking value in $[r, T]$,
\begin{align*}
\abs{y^i_{\tau}}&\leq C\left(\E^{\overline\P}\left[\left.\abs{\xi^i}^\kappa+\int_{r}^T\abs{g^i_s(0,0)}^\kappa ds\right|\overline\Fc_{\tau}^{+}\right]\right)^{\frac{1}{\kappa}},\ \overline\P-a.s.,
\end{align*}
}
and
$$\No{z^{i}}^p_{\H_{r}^p(\Pb)} +\No{m^i}^p_{\M_{r}^p(\Pb)}\leq C\left(\No{\xi^i}^p_{\L^p_{r}(\Pb)}+\E^\P\left[\int_{r}^T\abs{g^i_s(0,0)}^pds\right]\right).$$
Denoting $\delta \xi:=\xi^1-\xi^2$, $\delta y:=y^1-y^2$, $\delta z:=z^1-z^2$, $\delta m:=m^1-m^2$, $\delta g:= (g^1-g^2)(\cdot, y^1,z^1)$, we also have
{
\begin{align*}
\abs{\delta y_{\tau}}&\leq C\left(\E^{\overline\P}\left[\left.\abs{\delta\xi}^\kappa+\int_{r}^T\abs{\delta g_s}^\kappa ds\right|\overline\Fc_{\tau}^{+}\right]\right)^{\frac{1}{\kappa}},\ \overline\P-a.s.,
\end{align*}
}
and
$$\No{\delta z}^p_{\H_{r}^p(\Pb)} +\No{\delta m}^p_{\M_{r}^p(\Pb)}\leq C\left(\No{\delta\xi}^p_{\L^p_{r}(\Pb)}+\E^\Pb\left[\int_{r}^T\abs{\delta g_s}^pds\right]\right).$$
\end{Lemma}
\proof
See Section 4 of \cite{bouchard2015unified}.
\ep

\begin{Lemma}\label{lemma:bsde}
For any $\F-$stopping times $0\leq r \leq \rho\leq\tau\leq T$, any decreasing sequence of $\F-$stopping times $(\tau_n)_{n\geq 1}$ converging $\P-a.s.$ to $\tau$, and any $\F_+-$progressively measurable and right-continuous process $V \in\D^p_{r}(\F^{\P}_+,\P)$, if $y(\cdot,V_{\cdot})$ denotes the first component of the solution to the {\rm BSDE} \eqref{eq:BSDEp} on $[r,\cdot]$ with terminal condition $V_\cdot$ and some generator $g$, which satisfy Assumption \ref{assum:main}, we have
$$\E^{\Pb} \left[\abs{y_\rho(\tau,V_\tau)-y_\rho(\tau_n,V_{\tau_n})}\right]\underset{n\rightarrow +\infty}{\longrightarrow}0.$$
\end{Lemma}

\proof
First of all, by Lemma \reff{lemma:dppbsde}, we have
\begin{align*}
y_\rho(\tau,V_\tau)-y_\rho(\tau_n,V_{\tau_n})=y_\rho(\tau,V_\tau)-y_\rho(\tau,y_\tau(\tau_n,V_{\tau_n})).
\end{align*}
By Lemma \ref{lemma:estimbsde}, we therefore have for $\kappa\in(1,p]$
$$\E^\Pb\left[\abs{y_\rho(\tau,V_\tau)-y_\rho(\tau,y_\tau(\tau_n,V_{\tau_n}))}\right]\leq C\E^\Pb\left[\abs{V_\tau-y_\tau(\tau_n,V_{\tau_n})}^\kappa\right].$$
Next, again by a linearization argument, we can find bounded processes $\lambda$ and $\eta$ which are $\Fb-$progressively measurable such that
$$y_\tau(\tau_n,V_{\tau_n})=\E^{\P\otimes\P_0}\left[\left.\Ec\left(\int_{\tau}^{\tau_n}\eta_s\cdot dW_s^\P\right)\left(e^{\int_{\tau}^{\tau_n}\lambda_sds}V_{\tau_n}-\int_\tau^{\tau_n}e^{\int_\tau^s\lambda_udu}g_s(0,0)ds\right)\right|\overline\Fc^+_{\tau}\right].$$
Hence, choosing $\tilde p \in (\kappa, p)$
\begin{align*}
&\E^{\P\otimes\P_0}\left[\abs{y_\rho(\tau,V_\tau)-y_\rho(\tau,y_\tau(\tau_n,V_{\tau_n}))}\right]\\
&\leq  C\E^{\P\otimes\P_0}\left[\Ec\left(\int_{\tau}^{\tau_n}\eta_s\cdot dW_s^\P\right)^\kappa e^{\kappa\int_{\tau}^{\tau_n}\lambda_sds}\abs{V_{\tau_n}-V_\tau}^\kappa\right]\\
&\hspace{0.5em}+C\E^{\P\otimes\P_0}\left[\abs{1-\Ec\left(\int_{\tau}^{\tau_n}\eta_s\cdot dW_s^\P\right) e^{\int_{\tau}^{\tau_n}\lambda_sds}}^\kappa\abs{V_\tau}^\kappa\right]\\
&\hspace{0.5em}+C\E^{\P\otimes\P_0}\left[\Ec\left(\int_{\tau}^{\tau_n}\eta_s\cdot dW_s^\P\right)^\kappa \int_\tau^{\tau_n}e^{\kappa\int_{\tau}^{\tau_n}\lambda_sds}\abs{g_s(0,0)}^\kappa ds\right]\\
&\leq  \ C\left(\E^{\P\otimes\P_0}\left[\abs{V_{\tau_n}-V_\tau}^{\tilde p}\right]\right)^{\frac{\kappa}{\tilde p}}+C\left(\E^{\P\otimes\P_0}\left[\abs{1-\Ec\left(\int_{\tau}^{\tau_n}\eta_s\cdot dW_s^\P\right) e^{\int_{\tau}^{\tau_n}\lambda_sds}}^{\frac{p}{p-\kappa}}\right]\right)^{\frac{p-\kappa}{p}}\\
&\hspace{0.5em}+C\E^{\P\otimes\P_0}\left[ \int_\tau^{\tau_n}e^{\tilde p\int_{\tau}^{\tau_n}\lambda_sds}\abs{g_s(0,0)}^{\tilde p}ds\right],
\end{align*}
where we have used H\"{o}lder inequality, that $\lambda$ is bounded and the fact that since $\eta$ is also bounded, the Dol\'eans-Dade exponential appearing above has finite moments of any order. Now the terms inside the expectations on the right-hand side all converge in probability to $0$ and are clearly uniformly integrable by de la Vall\'ee--Poussin criterion since $V\in \D^p_{r}(\F^{\P}_+,\P)$ and $\tilde p<p$. We can therefore conclude by dominated convergence.
\ep

\begin{Lemma}[Comparison]\label{lemma:comp}
	Let Assumption \ref{assum:main} hold. Then, for $i=1,2$, let us denote by $(y^{i},z^{i},m^{i},k^i)$ the supersolution of the {\rm BSDE} \eqref{eq:BSDE_supsol} with generator $g^i$ and terminal condition $\xi^i$. If it holds $\mathbb P-a.s.$ that 
$$\xi_1\geq \xi_2,\ k^1-k^2 \text{ is non-decreasing and }g^1(s,y^1_s,z^1_s)\geq g^2(s,y^1_s,z^1_s),$$
 then we have for all $t\in [0,T]$
$$y_t^1\geq y_t^2,\ \mathbb P-a.s.$$
\end{Lemma}
\proof
We remind the reader that since $W^\P$ and $m^i$, $i=1,2$ are orthogonal and since $W^\P$ is actually continuous, we not only have $[W^\P,m^i]=0,\ \overline\P-a.s.$, but also
$$\langle W^\P, m^i\rangle = \langle W^\P, m^{i,c,\overline\P}\rangle=\langle W^\P, m^{i,d,\overline\P}\rangle=0,\ \overline\P-a.s.,$$
where $m^{i,c,\overline\P}$ (resp. $m^{i,d,\overline\P}$) is the continuous (resp. purely discontinuous) martingale part of $m^i$, under the measure $\overline\P$.

\vspace{0.5em}
\noindent Then, since the $g^i$ are uniformly Lipschitz, there exist two processes $\lambda$ and $\eta$ which are bounded, $\overline\P-a.s.$, and which are respectively $\overline\F^{\Pb}_+-$progressively measurable and $\overline\F^{\overline\P}-$predic-table, such that
$$ g^2(s,y_s^1,z_s^1)-g^2(s,y_s^2,z_s^2)=\lambda_t\left(y_s^1-y_s^2\right)+\eta_s\left(z_s^1-z_s^2\right),\ ds\times d\overline\P-a.e.$$
For any $0\leq t\leq s\leq T$, let us define the following continuous, positive and $\overline\F^{\Pb}_+-$progres-sively measurable process
$$A_{t,s}:=\exp\left(\int_t^s\lambda_udu-\int_t^s\eta_u\cdot dW^\P_u-\frac12\int_t^s\No{\eta_u}^2du\right).$$
By It\^o's formula, we deduce classically that
$$y_t^1-y_t^2=\mathbb E^{\overline{\mathbb P}}\left[\left. A_{t,T}(\xi^1-\xi^2)+\int_t^TA_{t,s}\left[(g^1-g^2)(s,y^1_s,z^1_s)ds+d(k^1_s-k^2_s)\right]\right|\overline\Fc_{t^+}\right],$$
from which we deduce immediately that $y^1_t\geq y^2_t,\ \overline\P-a.s.$
\ep

\end{appendix}

\section*{Acknowledgement}
{This work was started while Dylan Possama\"i and Xiaolu Tan were visiting the National University of Singapore, whose hospitality is kindly acknowledged. We are grateful to an anonymous referee as well as Nizar Touzi, Yiqing Lin, Zhenjie Ren and Junjian Yang for their helpful comments.}

\vspace{3mm}
Dylan Possama\"i gratefully acknowledges the financial support of the ANR Pacman ANR-16-CE05-0027.

\vspace{3mm}
Xiaolu Tan gratefully acknowledges the financial support of the ERC 321111 Rofirm, the ANR Isotace, and the Chairs Financial Risks (Risk Foundation, sponsored by Soci\'et\'e G\'en\'erale) and Finance and Sustainable Development (IEF sponsored by EDF and CA).

\vspace{3mm}
Chao Zhou gratefully acknowledges the financial support of Singapore MOE AcRF grants R-146-000-179-133 and R-146-000-219-112.

\end{document}